\documentclass[english,final,leqno,onefignum,onetabnum]{siamltex1213}
\usepackage[T1]{fontenc}
\usepackage[latin9]{inputenc}
\usepackage{array}
\usepackage{booktabs}
\usepackage{hyperref}
\usepackage{textcomp}
\usepackage{multirow}
\usepackage{amsmath}
\usepackage{algorithm2e}
\usepackage{amssymb}
\usepackage{graphicx}
\usepackage{esint}
\usepackage{babel}

\title{An Efficient Intrusive Uncertainty Propagation Method For Multi-Physics System With Random Inputs} 

\author{A.~Mittal\footnotemark[1] \footnotemark[3]
\and G. Iaccarino\footnotemark[2]}

\begin{document}
\maketitle
\slugger{sisc}{xxxx}{xx}{x}{x--x}%slugger should be set to mms, siap, sicomp, sicon, sidma, sima, simax, sinum, siopt, sisc, or sirev

\renewcommand{\thefootnote}{\fnsymbol{footnote}}
\footnotetext[1]{Institute for Computational and Mathematical Engineering, Stanford University, Stanford, CA 94305.}
\footnotetext[2]{Mechanical Engineering, Stanford University, Stanford, CA 94305.}
\footnotetext[3]{Corresponding author (Email: \email{mittal0@stanford.edu})}
\renewcommand{\thefootnote}{\arabic{footnote}}

\begin{abstract}
Coupled partial differential equation (PDE) systems, which often represent multi-physics models, are naturally suited for modular numerical solution methods. However, several challenges yet remain in extending the benefits of modularization practices to the task of uncertainty propagation. Since the cost of each deterministic PDE solve can be expected to be usually quite significant, statistical sampling based methods like Monte-Carlo (MC) are inefficient because they do not take advantage of the mathematical structure of the problem and suffer for poor convergence properties. On the other hand, even if each module contains a moderate number of uncertain parameters, implementing spectral methods on the combined high-dimensional parameter space can be prohibitively expensive due to the curse of dimensionality. In this work, we present a module-based and efficient intrusive spectral projection (ISP) method for uncertainty propagation. In our proposed method, each subproblem is separated and modularized via block Gauss-Seidel (BGS) techniques, such that each module only needs to tackle the local stochastic parameter space. Moreover, the computational costs are significantly mitigated by constructing reduced chaos approximations of the input data that enter each module. We demonstrate implementations of our proposed method and its computational gains over the standard ISP method using numerical examples.
\end{abstract}



\pagestyle{myheadings}
\thispagestyle{plain}
\markboth{A.~MITTAL AND G. IACCARINO}{INTRUSIVE UNCERTAINTY PROPAGATION}

\section{Introduction}
In recent years, Modeling and Simulation (M\&S) tools have become ubiquitous due to the rapid growth in computing power and significant developments in numerical techniques. Predictive science, which involves the utilization of verified and validated computer simulations in studying complex physical phenomena, has become especially useful in multi-physics applications involving heat transfer, fluid flow, structural dynamics, chemical kinetics and nuclear processes, due to the difficulty and costs associated with experimental testing. However, to validate numerical predictions against experimental observations, uncertainties (characterized as random quantities) must be included within the associated physical model. This is primarily because mathematical descriptions of physical processes are often unphysical idealizations of their target scenarios (in terms of material and geometric properties), and also because knowledge of the model parameters is usually limited. Uncertainty Quantification (UQ) has therefore emerged as a vital discipline in the certification of predictive simulations. Once the uncertain parameters of a physical system have been identified, a UQ study entails the following sequence of steps: (1) uncertainty characterization, (2) uncertainty propagation and (3) uncertainty analysis. 

In most cases where physical systems are formulated as partial differential equations (PDEs) with random inputs, uncertainty propagation becomes a compuationally demanding task. Monte-Carlo (MC) based sampling [\hyperlink{ref1}{1}, \hyperlink{ref2}{2}], wherein a set of input samples are generated and the corresponding output samples are collected via repeated simulation runs, is the simplest strategy to tackle this problem. However, obtaining a large enough sample size to perform accurate uncertainty analysis can be prohibitively expensive due to the slow convergence rate of MC. Alternative approaches, involving the construction of response surfaces of the quantities of interest using sparse regression methods [\hyperlink{ref3}{3}, \hyperlink{ref4}{4}] and subsequently sampling the surrogate model, have been proposed. However, in the context of multi-physics models, governed by coupled PDE systems, these methods ignore the coupling structures, which can be potentially exploited to reduce computational costs. Although spectral uncertainty propagation methods, for example, generalized polynomials chaos (gPC) [\hyperlink{ref5}{5}] based methods, exhibit geometric convergence properties, their monolithic implementations, would succumb to the curse of dimensionality. 

Multi-physics systems, are naturally suited for modular (partitioned) solution strategies [\hyperlink{ref6}{6}, \hyperlink{ref7}{7}], wherein specialized numerical methods and solvers for each constituent single-physics component can be reused and integrated with minimal effort, to construct solvers for the full system. These methods can therefore facilitate a natural division of M\&S expertise, and consequently a drastic reduction of developmental costs and maintanence overheads associated with multi-physics software. However, extending standard modularization practices for efficient uncertainty propagation yet remains challenging. 

In this work, we propose a module-based reduced intrusive spectral projection (ISP) method for gPC based uncertainty propagation in coupled PDE solvers. In addition to providing a practical framework for solver modularization, our proposed method mitigates the curse of dimensionality by constructing reduced approximations of the input data before being communicated to their respective computational modules. The reduced approximations are constructed via simple linear algebraic transformations, and corresponding optimal quadrature rules are constructed to significantly reduce the repeated number of module executions. Previous work has exploited relatively trivial coupling structures, for example, unidirectional coupling [\hyperlink{ref8}{8}] and linear coupling [\hyperlink{ref9}{9}, \hyperlink{ref10}{10}, \hyperlink{ref11}{11}], to reduce the costs of uncertainty propagation, while here, we consider a more general setup with bidirectional and nonlinear coupling structures in the model. Moreover, Constantine et. al [\hyperlink{ref12}{12}] successfully demonstrated a reduced approximation method for network (weakly) coupled systems and our proposed method extends the approach towards field (strongly) coupled models. 

Several components of our proposed approach have been motivated by the recent works of Arnst et. al [\hyperlink{ref13}{13}, \hyperlink{ref14}{14}, \hyperlink{ref15}{15}], wherein reduced chaos representations [\hyperlink{ref16}{16}] are proposed as the reduced approximations. A limitation of their approach however, is that computational gains can only be achieved in some of the modules, since the dimension and order reduction is implemented unidirectionally. We address this limitation in our proposed method by facilitating the construction of reduced approximations across all
the modules. Moreover, the respective tolerance values, which control the approximation errors, can be individually prescribed for each module.

The remainder of this article is organized as follows. In \hyperlink{sec2}{\S2}, the preliminary definitions and concepts concerning
ISP based uncertainty propagation for a stochastic multi-physics system are introduced along with a description of the 
standard ISP method. In \hyperlink{sec3}{\S3}, we describe the proposed
reduced ISP method by addressing each of its components individually. In \hyperlink{sec4}{\S4}, we
report and compare the performance and accuracy of the standard and reduced ISP method implementations on
two numerical examples.

\section{Preliminary definitions and concepts}

\subsection{Stochastic multi-physics model}

Without loss of generality, we consider to a multi-physics
setup wherin a two-component system of stochastic, steady-state
PDEs are bidirectionally coupled. The analysis that will be presented can be trivially extended to unsteady systems as well. A spatial discretization of each component of the system yields
the coupled stochastic algebraic system:$\forall\boldsymbol{\xi}=\left[\boldsymbol{\xi}_{1};\boldsymbol{\xi}_{2}\right]\in\Xi,$
\hypertarget{eq21}{}
\begin{align}
 & \boldsymbol{f}_{1}\left(\boldsymbol{u}_{1}\left(\boldsymbol{\xi}\right);\boldsymbol{u}_{2}\left(\boldsymbol{\xi}\right),\boldsymbol{\xi}_{1}\right)=\boldsymbol{0},\nonumber \\
 & \boldsymbol{f}_{2}\left(\boldsymbol{u}_{2}\left(\boldsymbol{\xi}\right);\boldsymbol{u}_{1}\left(\boldsymbol{\xi}\right),\boldsymbol{\xi}_{2}\right)=\boldsymbol{0},
\end{align}
where $\boldsymbol{f}_{1}\in\mathbb{R}^{n_{1}}$ and $\boldsymbol{f}_{2}\in\mathbb{R}^{n_{2}}$
denote the component residuals, $\boldsymbol{u}_{1}\in\mathbb{R}^{n_{1}}$
and $\boldsymbol{u}_{2}\in\mathbb{R}^{n_{2}}$ denote the component
solution vectors, which represent the respective finite-dimensional
discretizations of the spatially varying solution fields intrinsic
in each respective component PDE system. Moreover, $\boldsymbol{\xi}_{1}\in\Xi_{1}\subseteq\mathbb{R}^{s_{1}}$
and $\boldsymbol{\xi}_{2}\in\Xi_{2}\subseteq\mathbb{R}^{s_{2}}$ denote
the component random input parameters with probability density functions
$\mu_{1}:\Xi_{1}\rightarrow\mathbb{R}^{+}$ and $\mu_{2}:\Xi_{2}\rightarrow\mathbb{R}^{+}$
respectively. Furthermore, $\Xi\equiv\Xi_{1}\times\Xi{}_{2}$ denotes
the combined parameter space with dimension $s=s_{1}+s_{2}$ and joint
probability density function $\mu:\Xi\rightarrow\mathbb{R}^{+}:\forall\boldsymbol{\xi}\in\Xi,\mu\left(\boldsymbol{\xi}\right)=\mu_{1}\left(\boldsymbol{\xi}_{1}\left(\boldsymbol{\xi}\right)\right)\mu_{2}\left(\boldsymbol{\xi}_{2}\left(\boldsymbol{\xi}\right)\right),$
implying that $\boldsymbol{\xi}_{1}$ and $\boldsymbol{\xi}_{2}$
are statistically independent.

\subsection{Generalized Polynomial Chaos}

Let $\mathcal{L}^{2}\left(\Xi\right)$ denote the space of all $\mu$-weighted,
square-integrable scalar functions in $\Xi$. Moreover, $\forall i\in\left\{ 1,2\right\} $,
let $\boldsymbol{G}_{i}\in\mathbb{R}^{n_{1}\times n_{1}}$ denote
the symmetric-positive-definite Gramian matrix {[}\hyperlink{ref17}{17}{]} corresponding
to $\boldsymbol{u}_{i}$ and 
\begin{equation}
\mathcal{L}_{i}^{2}\left(\Xi\right)=\left\{ \boldsymbol{u}:\Xi\rightarrow\mathbb{R}^{n_{i}}:\int_{\Xi}\boldsymbol{u}\left(\boldsymbol{\xi}\right)^{\mathbf{T}}\boldsymbol{G}_{i}\boldsymbol{u}\left(\boldsymbol{\xi}\right)\mu\left(\boldsymbol{\xi}\right)d\boldsymbol{\xi}<\infty\right\} 
\end{equation}
denote the space of $\mu$-weighted, $\boldsymbol{G}_{i}$-square
integrable functions that map from $\Xi$ to $\mathbb{R}^{n_{i}}$.
If $\boldsymbol{u}_{i}\in\mathcal{L}_{i}^{2}\left(\Xi\right)$, then
it can be formulated as an infinite polynomial series as follows.
$\forall\boldsymbol{\xi}\in\Xi$,
\hypertarget{eq23}{}
\begin{equation}
\boldsymbol{u}_{i}\left(\boldsymbol{\xi}\right)=\sum_{j\geq0}\hat{\boldsymbol{u}}_{i,j}\psi_{j}\left(\boldsymbol{\xi}\right)
\end{equation}
where $\left\{ \psi_{j}:\Xi\rightarrow\mathbb{R}\right\} _{j\geq0}$
denotes the set of distinct $\mu$-orthonormal polynomials, and a
basis for $\mathcal{L}^{2}\left(\Xi\right)$. The coordinate directions
in $\Xi$ are assumed to be statistically independent and therefore,
each multivariate basis polynomial can be constructed as a product
of $s$ univariate orthonormal polynomials {[}\hyperlink{ref18}{18}{]}, which in turn
can be precomputed using the Golub-Welsch algorithm {[}\hyperlink{ref19}{19}{]}. Moreover,
the indexing in the polynomial series is assumed to follow a total
degree ordering, such that 
\begin{equation}
\mathrm{deg}\left(\psi_{j}\right)\geq\mathrm{deg}\left(\psi_{k}\right)\Leftrightarrow j\geq k\geq0,
\end{equation}
where $\forall j\geq0$, $\mathrm{deg}\left(\psi_{j}\right)$ denotes
the total degree of $\psi_{j}$.

Consequently, a gPC approximation {[}\hyperlink{ref5}{5}{]} $\boldsymbol{u}_{i}^{p}\approx\boldsymbol{u}_{i}$
of order $p\geq0$ can be formulated by truncating \hyperlink{eq23}{Eq. 2.3} as follows.
$\forall\boldsymbol{\xi}\in\Xi$, 
\begin{equation}
\boldsymbol{u}_{i}^{p}\left(\boldsymbol{\xi}\right)=\sum_{j=0}^{P}\hat{\boldsymbol{u}}_{i,j}\psi_{j}\left(\boldsymbol{\xi}\right)=\hat{\boldsymbol{U}}_{i}\boldsymbol{\psi}\left(\boldsymbol{\xi}\right),
\end{equation}
where $\hat{\boldsymbol{U}}_{i}\equiv\hat{\boldsymbol{U}}_{i}^{p}=\left[\begin{array}{ccc}
\hat{\boldsymbol{u}}_{i,0} & \cdots & \hat{\boldsymbol{u}}_{i,P}\end{array}\right]\in\mathbb{R}^{n_{i}\times\left(P+1\right)}$ denotes the gPC coefficient matrix and $\boldsymbol{\psi}\equiv\boldsymbol{\psi}^{p}=\left[\begin{array}{ccc}
\psi_{0} & \cdots & \psi_{P}\end{array}\right]^{\mathbf{T}}:\Xi\rightarrow\mathbb{R}^{P+1}$ denotes the basis vector. If the truncation is isotropic and based
on the total degree, then the number of coefficients $P+1={\displaystyle \frac{\left(p+s\right)!}{p!s!}}.$ 

The Cameron-Martin theorem {[}\hyperlink{ref20}{20}{]} states that if $\boldsymbol{u}_{i}$
is infinitely regular in $\Xi$, then as $p\rightarrow\infty$, the
gPC approximation $\boldsymbol{u}_{i}^{p}$ converges exponentially
to $\boldsymbol{u}_{i}$, in the mean-square sense. To state this
formally, $\exists\chi^{*}>0,\rho^{*}>1$, such that 
\begin{equation}
\sqrt{\int_{\Xi}\left\Vert \boldsymbol{u}_{i}^{p}\left(\boldsymbol{\xi}\right)-\boldsymbol{u}_{i}\left(\boldsymbol{\xi}\right)\right\Vert _{\boldsymbol{G}_{i}}^{2}\mu\left(\boldsymbol{\xi}\right)d\boldsymbol{\xi}}\leq\chi^{*}\rho^{-p}
\end{equation}
for some $\rho>\rho^{*}$. Moreover, a corollary to the theorem states
that if the regularity of $\boldsymbol{u}_{i}$ is $k$, then the
asymptotic rate of convergence is polynomial and the upper bound in
the mean-square error would be $\mathcal{O}\left(p^{-k}\right)$.

For computing statistical quantities of interest such as the moments
and related probability density functions of the solutions via exhaustive
sampling, the gPC approximation can be used as a significantly cheaper
surrogate model when compared to repeated executions of a costly numerical
PDE solver. In particular, approximations of the first two moments
can be computed directly from the gPC coefficients as follows. 
\begin{align}
 & \mathbb{E}\left(\boldsymbol{u}_{i}\right)\approx\mathbb{E}\left(\boldsymbol{u}_{i}^{p}\right)=\hat{\boldsymbol{u}}_{i,0},\nonumber \\
 & \mathrm{Cov}\left(\boldsymbol{u}_{i},\boldsymbol{u}_{i}\right)\approx\mathrm{Cov}\left(\boldsymbol{u}_{i}^{p},\boldsymbol{u}_{i}^{p}\right)=\hat{\boldsymbol{U}}_{i}\hat{\boldsymbol{U}}_{i}^{\mathbf{T}}-\hat{\boldsymbol{u}}_{i,0}\hat{\boldsymbol{u}}_{i,0}^{\mathbf{T}}.
\end{align}

Moreover, the probability distribution function of any related quantity
of interest can be computed with the kernel density estimation (KDE)
method {[}\hyperlink{ref21}{21}{]}, using a large number of surrogate solution samples.
Therefore, prior to any uncertainty analysis, the coefficient matrices
of the solutions must be computed via intrusive spectral projection
(ISP) or non-intrusive spectral projection (NISP) based uncertainty
propagation methods. The focus of this work is the stochastic Galerkin
method {[}\hyperlink{ref22}{22}{]}, which is a commonly used ISP based uncertainty
propagation method.

\subsection{Intrusive spectral projection}

Let $\mathcal{P}\left(\Xi\right)\equiv\mathcal{P}^{p}\left(\Xi\right)\equiv\left\{ \boldsymbol{v}^{\mathbf{T}}\boldsymbol{\psi}:\boldsymbol{v}\in\mathbb{R}^{P+1}\right\} $
denote the space of polynomials in $\Xi$, with total degree $\leq p$.
Subsequently, the Galerkin form of the coupled stochastic
system (\hyperlink{eq21}{Eq. 2.1}) can be expressed as follows. Find $\hat{\boldsymbol{U}}{}_{1}\in\mathbb{R}^{n_{1}\times\left(P+1\right)}$
and $\hat{\boldsymbol{U}}{}_{2}\in\mathbb{R}^{n_{2}\times\left(P+1\right)}$,
such that $\forall\omega_{1},\omega_{2}\in\mathcal{P}\left(\Xi\right)$,
\begin{align}
\hypertarget{eq29}{}
 & \int_{\Xi}\boldsymbol{f}_{1}\left(\hat{\boldsymbol{U}}{}_{1}\boldsymbol{\psi}\left(\boldsymbol{\xi}\right);\hat{\boldsymbol{U}}{}_{2}\boldsymbol{\psi}\left(\boldsymbol{\xi}\right),\hat{\boldsymbol{\Xi}}{}_{1}\boldsymbol{\psi}\left(\boldsymbol{\xi}\right)\right)\omega_{1}\left(\boldsymbol{\xi}\right)\mu\left(\boldsymbol{\xi}\right)d\boldsymbol{\xi}\nonumber \\
 & +\int_{\Xi}\boldsymbol{f}_{2}\left(\hat{\boldsymbol{U}}{}_{2}\boldsymbol{\psi}\left(\boldsymbol{\xi}\right);\hat{\boldsymbol{U}}{}_{1}\boldsymbol{\psi}\left(\boldsymbol{\xi}\right),\hat{\boldsymbol{\Xi}}{}_{2}\boldsymbol{\psi}\left(\boldsymbol{\xi}\right)\right)\omega_{2}\left(\boldsymbol{\xi}\right)\mu\left(\boldsymbol{\xi}\right)d\boldsymbol{\xi}=\boldsymbol{0}.\\
\Rightarrow & \hat{\boldsymbol{F}}_{1}\left(\hat{\boldsymbol{U}}{}_{1};\hat{\boldsymbol{U}}{}_{2}\right)=\int_{\Xi}\boldsymbol{f}_{1}\left(\hat{\boldsymbol{U}}{}_{1}\boldsymbol{\psi}\left(\boldsymbol{\xi}\right);\hat{\boldsymbol{U}}{}_{2}\boldsymbol{\psi}\left(\boldsymbol{\xi}\right),\hat{\boldsymbol{\Xi}}{}_{1}\boldsymbol{\psi}\left(\boldsymbol{\xi}\right)\right)\boldsymbol{\psi}\left(\boldsymbol{\xi}\right)^{\mathbf{T}}\mu\left(\boldsymbol{\xi}\right)d\boldsymbol{\xi}=\boldsymbol{0},\nonumber \\
 & \hat{\boldsymbol{F}}_{2}\left(\hat{\boldsymbol{U}}{}_{2};\hat{\boldsymbol{U}}{}_{1}\right)=\int_{\Xi}\boldsymbol{f}_{2}\left(\hat{\boldsymbol{U}}{}_{2}\boldsymbol{\psi}\left(\boldsymbol{\xi}\right);\hat{\boldsymbol{U}}{}_{1}\boldsymbol{\psi}\left(\boldsymbol{\xi}\right),\hat{\boldsymbol{\Xi}}{}_{2}\boldsymbol{\psi}\left(\boldsymbol{\xi}\right)\right)\boldsymbol{\psi}\left(\boldsymbol{\xi}\right)^{\mathbf{T}}\mu\left(\boldsymbol{\xi}\right)d\boldsymbol{\xi}=\boldsymbol{0}.
\end{align}
where $\forall i\in\left\{ 1,2\right\} $, $\hat{\boldsymbol{F}}_{i}\equiv\hat{\boldsymbol{F}}_{i}^{p}$
and $\hat{\boldsymbol{\Xi}}_{i}\equiv\hat{\boldsymbol{\Xi}}_{i}^{p}\in\mathbb{R}^{s_{i}\times\left(P+1\right)}:$
\begin{align}
\hat{\boldsymbol{\Xi}}_{i} & =\int_{\Xi}\boldsymbol{\xi}_{i}\left(\boldsymbol{\xi}\right)\boldsymbol{\psi}\left(\boldsymbol{\xi}\right)^{\mathbf{T}}\mu\left(\boldsymbol{\xi}\right)d\boldsymbol{\xi}
\end{align}
denote the gPC coefficient matrices of component residual $\boldsymbol{f}_{i}$
and random input vector $\boldsymbol{\xi}_{i}$ respectively. Alternatively,
the transposed Galerkin form {[}\hyperlink{ref23}{23}{]}, can be derived by projecting
the coupled PDE system in $\mathcal{P}\left(\Xi\right)$, and subsequently
applying the respective spatial discretization schemes on the deterministic
PDE system.

\subsection{Numerical solution methods}

\hyperlink{eq29}{Eq. 2.9} can be numerically solved using either monolithic (fully-coupled)
or modular (partitioned) methods. A monolithic approach, for example, Newton's method, would
require the solution of the fully-coupled linear system 
\hypertarget{eq211}{}
\begin{equation}
\left[\begin{array}{cc}
{\displaystyle \frac{\partial\underline{{\displaystyle \hat{\boldsymbol{F}}_{1}}}}{\partial\underline{\hat{\boldsymbol{U}}_{1}}}}\left(\hat{\boldsymbol{U}}_{1}^{\ell};\hat{\boldsymbol{U}}_{2}^{\ell}\right) & {\displaystyle \frac{\partial\underline{{\displaystyle \hat{\boldsymbol{F}}_{1}}}}{\partial\underline{\hat{\boldsymbol{U}}_{2}}}}\left(\hat{\boldsymbol{U}}_{1}^{\ell};\hat{\boldsymbol{U}}_{2}^{\ell}\right)\\
\\
{\displaystyle \frac{\partial\underline{{\displaystyle \hat{\boldsymbol{F}}_{2}}}}{\partial\underline{\hat{\boldsymbol{U}}_{1}}}}\left(\hat{\boldsymbol{U}}_{2}^{\ell};\hat{\boldsymbol{U}}_{1}^{\ell}\right) & {\displaystyle \frac{\partial\underline{{\displaystyle \hat{\boldsymbol{F}}_{2}}}}{\partial\underline{\hat{\boldsymbol{U}}_{2}}}}\left(\hat{\boldsymbol{U}}_{2}^{\ell};\hat{\boldsymbol{U}}_{1}^{\ell}\right)
\end{array}\right]\left[\begin{array}{c}
\underline{\hat{\boldsymbol{U}}_{1}^{\ell+1}}-\underline{\hat{\boldsymbol{U}}_{1}^{\ell}}\\
\\
\\
\underline{\hat{\boldsymbol{U}}_{2}^{\ell+1}}-\underline{\hat{\boldsymbol{U}}_{2}^{\ell}}
\end{array}\right]=-\left[\begin{array}{c}
\underline{{\displaystyle \hat{\boldsymbol{F}}_{1}}}\left(\hat{\boldsymbol{U}}_{1}^{\ell};\hat{\boldsymbol{U}}_{2}^{\ell}\right)\\
\\
\\
\underline{{\displaystyle \hat{\boldsymbol{F}}_{2}}}\left(\hat{\boldsymbol{U}}_{2}^{\ell};\hat{\boldsymbol{U}}_{1}^{\ell}\right)
\end{array}\right],
\end{equation}
to obtain the solution updates and converge to the solution of the
coupled Galerkin system (\hyperlink{eq29}{Eq. 2.9}). Here, $\underline{\ \cdot\ }$
simply denotes the vectorization operator. In general, developing
linear solvers (and associated preconditioners) for \hyperlink{eq211}{Eq. 2.11} is quite
challenging since each component of the coupled system may require
disparate numerical treatment. For example, fluid-structure interaction
models {[}\hyperlink{ref24}{24}{]} contains elliptic and hyperbolic PDEs, which are spatially
discretized using Finite Element and Finite Volume methods respectively.
Moreover, the quadratic convergence rate may not be guaranteed in
practice. Furthermore, variants of Newton's method, for instance,
Gauss-Newton {[}\hyperlink{ref25}{25}{]} and Levenberg-Marquardt {[}\hyperlink{ref26}{26}{]}, when implemented
in this monolithic fashion, would also be affected by these limitations. 

A modular solution method, for example, the block-Gauss-Seidel (BGS)
method {[}\hyperlink{ref6}{6}{]}, can address some of these limitations. Primarily,
separate solvers for each subproblem in \hyperlink{eq29}{Eq. 2.9}, possibly employing
disparate solution methods and numerical treatment, can be developed
with some degree of independence and coupled together as modules with
minimal modifications. In what follows, \hyperlink{alg1}{Algorithm 1} describes the standard ISP based
uncertainty propagation method, which is based on the BGS iterative
method. Here, $\hat{\boldsymbol{M}}_{1}\equiv\hat{\boldsymbol{M}}_{1}^{p}$
and $\hat{\boldsymbol{M}}_{2}\equiv\hat{\boldsymbol{M}}_{2}^{p}$
denote the module operators used in the respective solvers for ${\displaystyle \hat{\boldsymbol{F}}_{1}}=\boldsymbol{0}$
and ${\displaystyle \hat{\boldsymbol{F}}_{2}}=\boldsymbol{0}$.

A more general formulation of a two-component stochastic
multi-physics model can be defined as follows. $\forall\boldsymbol{\xi}=\left[\boldsymbol{\xi}_{1};\boldsymbol{\xi}_{2}\right]\in\Xi,$
\hypertarget{eq212}{}
\begin{align}
 & \boldsymbol{f}_{1}\left(\boldsymbol{u}_{1}\left(\boldsymbol{\xi}\right);\boldsymbol{v}_{2}\left(\boldsymbol{\xi}\right),\boldsymbol{\xi}_{1}\right)=\boldsymbol{0},\ \boldsymbol{v}_{1}\left(\boldsymbol{\xi}\right)=\boldsymbol{g}_{1}\left(\boldsymbol{u}_{1}\left(\boldsymbol{\xi}\right)\right),\nonumber \\
 & \boldsymbol{f}_{2}\left(\boldsymbol{u}_{2}\left(\boldsymbol{\xi}\right);\boldsymbol{v}_{1}\left(\boldsymbol{\xi}\right),\boldsymbol{\xi}_{2}\right)=\boldsymbol{0},\ \boldsymbol{v}_{2}\left(\boldsymbol{\xi}\right)=\boldsymbol{g}_{2}\left(\boldsymbol{u}_{2}\left(\boldsymbol{\xi}\right)\right),
\end{align}
where $\boldsymbol{g}_{1}\in\mathbb{R}^{m_{1}}$ and $\boldsymbol{g}_{2}\in\mathbb{R}^{m_{2}}$
denote the interface (coupling) functions. For this system, the BGS
method for ISP based uncertainty propagation can be characterized
as a slight modification of \hyperlink{alg1}{Algorithm 1}, wherein the iterations are
formulated as follows.
\begin{align}
 & \hat{\boldsymbol{U}}_{1}^{\ell+1}=\hat{\boldsymbol{M}}_{1}\left(\hat{\boldsymbol{U}}_{1}^{\ell};\hat{\boldsymbol{V}}_{2}^{\ell}\right),\ \hat{\boldsymbol{V}}_{1}^{\ell+1}=\hat{\boldsymbol{G}}_{1}\left(\hat{\boldsymbol{U}}_{1}^{\ell+1}\right),\nonumber \\
 & \hat{\boldsymbol{U}}_{2}^{\ell+1}=\hat{\boldsymbol{M}}_{2}\left(\hat{\boldsymbol{U}}_{2}^{\ell};\hat{\boldsymbol{V}}_{1}^{\ell}\right),\ \hat{\boldsymbol{V}}_{2}^{\ell+1}=\hat{\boldsymbol{G}}_{2}\left(\hat{\boldsymbol{U}}_{2}^{\ell+1}\right),
\end{align}
where $\forall i\in\left\{ 1,2\right\} $, $\boldsymbol{G}_{i}\equiv\boldsymbol{G}_{i}^{p}$
denotes the interface operator, such that $\forall\hat{\boldsymbol{U}}\in\mathbb{R}^{n_{i}\times\left(P+1\right)}$,
\begin{equation}
\boldsymbol{G}_{i}\left(\hat{\boldsymbol{U}}\right)=\int_{\Xi}\boldsymbol{g}_{i}\left(\hat{\boldsymbol{U}}\boldsymbol{\psi}\left(\boldsymbol{\xi}\right)\right)\boldsymbol{\psi}\left(\boldsymbol{\xi}\right)^{\mathbf{T}}\mu\left(\boldsymbol{\xi}\right)d\boldsymbol{\xi}.
\end{equation}

While modular solution methods are well suited for coupling disparate
numerical strategies for various components of the multi-physics model,
it is worth noting that the convergence rate is, in general, linear
{[}\hyperlink{ref27}{27}{]}.

\subsubsection{Computational cost and limitations}

Let $\bar{\mathcal{C}}_{1}\equiv\bar{\mathcal{C}}_{1}\left(n_{1}\right)$
and $\bar{\mathcal{C}}_{2}\equiv\bar{\mathcal{C}}_{2}\left(n_{2}\right)$
denote the respective average costs of solving $\boldsymbol{f}_{1}=\boldsymbol{0}$
and $\boldsymbol{f}_{2}=\boldsymbol{0}$ deterministically, given
a set of input realizations in $\Xi$. Therefore, the overall computational
cost of the standard ISP method 
\begin{equation}
\mathcal{C}_{s}\approx\mathcal{O}\left(\bar{\mathcal{C}}_{1}P^{\alpha_{1}}+\bar{\mathcal{C}}_{2}P^{\alpha_{2}}\right),
\end{equation}
where, in general, $\alpha_{1},\alpha_{2}>1$. Therefore, the $\mathcal{C}_{s}$
would grow exponentially with respect to the combined stochastic dimension
$s$ and order $p$. This undesirable complexity, known as the curse
of dimensionality, would still persist if a monolithic solution strategy
were used instead. 

Moreover, since the Galerkin form of each subproblem is formulated
using polynomials in the global (combined) stochastic parameter space
$\Xi$, modifying the characteristics of any particular module's stochastic
inputs would have to be reflected in all other modules. This requirement
limits the developmental independence of modules and therefore, limits
the practical benefits of modularization. 

In the next section, we describe the proposed reduced ISP based uncertainty
propagation method, which mitigates the curse of dimensionality via
reduced approximation methods, and facilitates a desirable module-based
independence by eliminating the influence of external stochastic inputs
in formulating the Galerkin form of each subproblem. 

\RestyleAlgo{boxruled}

\begin{algorithm}
\hypertarget{alg1}{}
\caption{Standard ISP based uncertainty propagation for a two-module multi-physics
system}
\SetAlgoLined

\SetKwInOut{Input}{inputs}\SetKwInOut{Output}{outputs}
\DontPrintSemicolon
\Input{order $p\geq0$, $\hat{\boldsymbol{U}}_{1}^{0}$, $\hat{\boldsymbol{U}}_{2}^{0}$}
\Output{$\hat{\boldsymbol{U}}_{1}$, $\hat{\boldsymbol{U}}_{2}$}

\textbf{$\ell\leftarrow0$}

\Repeat{$\hat{\boldsymbol{U}}_{1}^{\ell}$, $\hat{\boldsymbol{U}}_{2}^{\ell}$ $\mathrm{converge}$}{

$\hat{\boldsymbol{U}}_{1}^{\ell+1}\leftarrow\hat{\boldsymbol{M}}_{1}\left(\hat{\boldsymbol{U}}_{1}^{\ell};\hat{\boldsymbol{U}}_{2}^{\ell}\right)$

$\hat{\boldsymbol{U}}_{2}^{\ell+1}\leftarrow\hat{\boldsymbol{M}}_{2}\left(\hat{\boldsymbol{U}}_{2}^{\ell};\hat{\boldsymbol{U}}_{1}^{\ell+1}\right)$

\textbf{$\ell\leftarrow \ell+1$}

}
\end{algorithm}

\section{Reduced ISP based uncertainty propagation}

We will now individually describe the components of the proposed algorithm.

\subsection{Modular gPC approximation}

Since the component inputs are assumed to be statistically independent,
then $\forall i\in\left\{ 1,2\right\} $, the gPC approximation of
component solution $\boldsymbol{u}_{i}$ can be rewritten
as follows. $\forall\boldsymbol{\xi}_{1}\in\Xi_{1},\boldsymbol{\xi}_{2}\in\Xi_{2}$,
\hypertarget{eq31}{}
\begin{align}
\boldsymbol{u}_{i}^{p}\left(\boldsymbol{\xi}_{1},\boldsymbol{\xi}_{2}\right) & =\sum_{j=0}^{P}\hat{\boldsymbol{u}}_{i,j}\psi_{j}\left(\boldsymbol{\xi}_{1},\boldsymbol{\xi}_{2}\right)\nonumber \\
 & =\sum_{j=0}^{P}\hat{\boldsymbol{u}}_{i,\jmath_{1}\left(j\right),\jmath_{2}\left(j\right)}\psi_{1,\jmath_{1}\left(j\right)}\left(\boldsymbol{\xi}_{1}\right)\psi_{2,\jmath_{2}\left(j\right)}\left(\boldsymbol{\xi}_{2}\right),
\end{align}
where $\forall i\in\left\{ 1,2\right\} $, $\left\{ \psi_{i,j}:\Xi_{i}\rightarrow\mathbb{R}\right\} _{j\geq0}$
denotes the set of distinct $\mu_{i}-$orthonormal polynomials ordered
according to their total degree and $\jmath_{i}\equiv\jmath_{i}^{p}:\left\{ 0\leq j\leq P\right\} \rightarrow\left\{ 0\leq j\leq P_{i}\right\} $
denotes the map from the global index set to the modular index set.
Since the total degree of the polynomial expansion must be $\leq p$,
in the isotropic case, the number of distinct polynomials needed is
$P_{i}+1={\displaystyle \frac{\left(p+s_{i}\right)}{p!s_{i}!}}$ .

$\forall i\in\left\{ 1,2\right\} $, let $\boldsymbol{\psi}_{i}\equiv\boldsymbol{\psi}_{i}^{p}=\left[\begin{array}{ccc}
\psi_{i,0} & \cdots & \psi_{i,P_{i}}\end{array}\right]^{\mathbf{T}}:\Xi_{i}\rightarrow\mathbb{R}^{P_{i}+1}$ denote the modular basis vector. Therefore, from \hyperlink{eq31}{Eq. 3.1}, the global
basis vector $\boldsymbol{\psi}$ can be formulated in terms of $\boldsymbol{\psi}_{i}$
as follows. $\forall\boldsymbol{\xi}_{1}\in\Xi_{1},\boldsymbol{\xi}_{2}\in\Xi_{2}$,
\hypertarget{eq32}{}
\begin{equation}
\boldsymbol{\psi}\left(\boldsymbol{\xi}_{1},\boldsymbol{\xi}_{2}\right)=\boldsymbol{\Pi}_{1}\left(\boldsymbol{\xi}_{2}\right)\boldsymbol{\psi}_{1}\left(\boldsymbol{\xi}_{1}\right)=\boldsymbol{\Pi}_{2}\left(\boldsymbol{\xi}_{1}\right)\boldsymbol{\psi}_{2}\left(\boldsymbol{\xi}_{2}\right),
\end{equation}
where $\boldsymbol{\Pi}_{1}:\Xi_{2}\rightarrow\mathbb{R}^{\left(P+1\right)\times\left(P_{1}+1\right)}$
and $\boldsymbol{\Pi}_{2}:\Xi_{1}\rightarrow\mathbb{R}^{\left(P+1\right)\times\left(P_{2}+1\right)}$
denote sparse polynomial matrices with at-most $P+1$ non-zero elements.
$\forall i\in\left\{ 1,2\right\} $, $0\leq j\leq P,0\leq k\leq P_{i}$,
let $\Pi_{i,j,k}$ denote the $\left(j+1,k+1\right)-th$ element in
$\boldsymbol{\Pi}_{i}$. Therefore, $\forall\boldsymbol{\xi}_{1}\in\Xi_{1},\boldsymbol{\xi}_{2}\in\Xi_{2}$,
\hypertarget{eq33}{}
\begin{align}
\Pi_{1,j,k}\left(\boldsymbol{\xi}_{2}\right) & =\begin{cases}
\psi_{2,\jmath_{2}\left(j\right)}\left(\boldsymbol{\xi}_{2}\right) & k=\jmath_{1}\left(j\right)\\
0 & k\neq\jmath_{1}\left(j\right)
\end{cases},\nonumber \\
\Pi_{2,j,k}\left(\boldsymbol{\xi}_{1}\right) & =\begin{cases}
\psi_{1,\jmath_{1}\left(j\right)}\left(\boldsymbol{\xi}_{1}\right) & k=\jmath_{2}\left(j\right)\\
0 & k\neq\jmath_{2}\left(j\right)
\end{cases}.
\end{align}

Subsequently, the gPC approximations of $\boldsymbol{u}_{1}$ and
$\boldsymbol{u}_{2}$ can be formulated as follows. $\forall\boldsymbol{\xi}_{1}\in\Xi_{1},\boldsymbol{\xi}_{2}\in\Xi_{2},$
\begin{equation}
\boldsymbol{u}_{1}^{p}\left(\boldsymbol{\xi}_{1},\boldsymbol{\xi}_{2}\right)=\tilde{\boldsymbol{U}}_{1}\left(\boldsymbol{\xi}_{2}\right)\boldsymbol{\psi}_{1}\left(\boldsymbol{\xi}_{1}\right),\ \boldsymbol{u}_{2}^{p}\left(\boldsymbol{\xi}_{1},\boldsymbol{\xi}_{2}\right)=\tilde{\boldsymbol{U}}_{2}\left(\boldsymbol{\xi}_{1}\right)\boldsymbol{\psi}_{2}\left(\boldsymbol{\xi}_{2}\right),
\end{equation}
where $\tilde{\boldsymbol{U}}_{1}\equiv\tilde{\boldsymbol{U}}_{1}^{p}:\Xi_{2}\rightarrow\mathbb{R}^{n_{1}\times\left(P_{1}+1\right)}$
and $\tilde{\boldsymbol{U}}_{2}\equiv\tilde{\boldsymbol{U}}_{2}^{p}:\Xi_{1}\rightarrow\mathbb{R}^{n_{2}\times\left(P_{2}+1\right)}$
denote the modular gPC coefficient matrices of $\boldsymbol{u}_{1}$
and $\boldsymbol{u}_{2}$ respectively. From \hyperlink{eq32}{Eq. 3.2}, these matrices
are related to the respective global gPC coefficient matrices $\hat{\boldsymbol{U}}_{1}$
and $\hat{\boldsymbol{U}}_{2}$ as follows. $\forall\boldsymbol{\xi}_{1}\in\Xi_{1},\boldsymbol{\xi}_{2}\in\Xi_{2}$,

\hypertarget{eq35}{}
\begin{equation}
\tilde{\boldsymbol{U}}_{1}\left(\boldsymbol{\xi}_{2}\right)=\hat{\boldsymbol{U}}_{1}\boldsymbol{\Pi}_{1}\left(\boldsymbol{\xi}_{2}\right),\ \tilde{\boldsymbol{U}}_{2}\left(\boldsymbol{\xi}_{1}\right)=\hat{\boldsymbol{U}}_{2}\boldsymbol{\Pi}_{2}\left(\boldsymbol{\xi}_{1}\right).
\end{equation}

Moreover, from their definition in \hyperlink{eq33}{Eq. 3.3}, $\boldsymbol{\Pi}_{1}$
and $\boldsymbol{\Pi}_{2}$ satisfy
\hypertarget{eq36}{}
\begin{equation}
\int_{\Xi_{2}}\boldsymbol{\Pi}_{1}\left(\boldsymbol{\xi}_{2}\right)\boldsymbol{\Pi}_{1}\left(\boldsymbol{\xi}_{2}\right)^{\mathbf{T}}\mu_{2}\left(\boldsymbol{\xi}_{2}\right)d\boldsymbol{\xi}_{2}=\int_{\Xi_{1}}\boldsymbol{\Pi}_{2}\left(\boldsymbol{\xi}_{1}\right)\boldsymbol{\Pi}_{2}\left(\boldsymbol{\xi}_{1}\right)^{\mathbf{T}}\mu_{1}\left(\boldsymbol{\xi}_{1}\right)d\boldsymbol{\xi}_{1}=\boldsymbol{I}_{P+1}.
\end{equation}

Therefore, from \hyperlink{eq35}{Eq. 3.5} and \hyperlink{eq36}{Eq. 3.6}, we have

\begin{align}
\hat{\boldsymbol{U}}_{1} & =\int_{\Xi_{2}}\tilde{\boldsymbol{U}}_{1}\left(\boldsymbol{\xi}_{2}\right)\boldsymbol{\Pi}_{1}\left(\boldsymbol{\xi}_{2}\right)^{\mathbf{T}}\mu_{2}\left(\boldsymbol{\xi}_{2}\right)d\boldsymbol{\xi}_{2},\nonumber \\
\hat{\boldsymbol{U}}_{2} & =\int_{\Xi_{1}}\tilde{\boldsymbol{U}}_{2}\left(\boldsymbol{\xi}_{2}\right)\boldsymbol{\Pi}_{2}\left(\boldsymbol{\xi}_{1}\right)^{\mathbf{T}}\mu_{1}\left(\boldsymbol{\xi}_{1}\right)d\boldsymbol{\xi}_{1}.
\end{align}

\subsection{Modular Galerkin form}

$\forall i\in\left\{ 1,2\right\} $, let $\mathcal{P}_{i}\left(\Xi_{i}\right)\equiv\mathcal{P}_{i}^{p}\left(\Xi_{i}\right)\equiv\left\{ \boldsymbol{v}^{\mathbf{T}}\boldsymbol{\psi}_{i}:\boldsymbol{v}\in\mathbb{R}^{P_{i}+1}\right\} $
denote the space of polynomials in $\Xi_{i}$ with total degree $\leq p$.
Subsequently, the modular Galerkin form of each component
of the coupled stochastic system (\hyperlink{eq21}{Eq. 2.1}) can be expressed as follows.
In module 1, find $\tilde{\boldsymbol{U}}{}_{1}:\Xi_{2}\rightarrow\mathbb{R}^{n_{1}\times\left(P_{1}+1\right)}$,
such that $\forall\hat{\boldsymbol{U}}_{2}\in\mathbb{R}^{n_{2}\times\left(P+1\right)}$,
$\boldsymbol{\xi}_{2}\in\Xi_{2}$, $\omega_{1}\in\mathcal{P}_{1}\left(\Xi_{1}\right)$,
\begin{align}
 & \int_{\Xi_{1}}\boldsymbol{f}_{1}\left(\tilde{\boldsymbol{U}}{}_{1}\left(\boldsymbol{\xi}_{2}\right)\boldsymbol{\psi}_{1}\left(\boldsymbol{\xi}_{1}\right);\hat{\boldsymbol{U}}_{2}\boldsymbol{\Pi}_{1}\left(\boldsymbol{\xi}_{2}\right)\boldsymbol{\psi}_{1}\left(\boldsymbol{\xi}_{1}\right),\tilde{\boldsymbol{\Xi}}{}_{1}\boldsymbol{\psi}_{1}\left(\boldsymbol{\xi}_{1}\right)\right)\omega_{1}\left(\boldsymbol{\xi}_{1}\right)\mu_{1}\left(\boldsymbol{\xi}_{1}\right)d\boldsymbol{\xi}_{1}\nonumber \\
 & =\boldsymbol{0}\\
\Rightarrow & \tilde{\boldsymbol{F}}_{1}\left(\tilde{\boldsymbol{U}}{}_{1};\hat{\boldsymbol{U}}{}_{2},\boldsymbol{\xi}_{2}\right)=\int_{\Xi_{1}}\left(\boldsymbol{f}_{1}\left(\tilde{\boldsymbol{U}}{}_{1}\left(\boldsymbol{\xi}_{2}\right)\boldsymbol{\psi}_{1}\left(\boldsymbol{\xi}_{1}\right);\hat{\boldsymbol{U}}_{2}\boldsymbol{\Pi}_{1}\left(\boldsymbol{\xi}_{2}\right)\boldsymbol{\psi}_{1}\left(\boldsymbol{\xi}_{1}\right),\tilde{\boldsymbol{\Xi}}{}_{1}\boldsymbol{\psi}_{1}\left(\boldsymbol{\xi}_{1}\right)\right)\right.\nonumber \\
 & \times\left.\boldsymbol{\psi}_{1}\left(\boldsymbol{\xi}_{1}\right)^{\mathbf{T}}\right)\mu_{1}\left(\boldsymbol{\xi}_{1}\right)d\boldsymbol{\xi}_{1}=\boldsymbol{0}.
\end{align}
Similarly, in module 2, find $\tilde{\boldsymbol{U}}{}_{2}:\Xi_{1}\rightarrow\mathbb{R}^{n_{2}\times\left(P_{2}+1\right)}$,
such that $\forall\hat{\boldsymbol{U}}_{1}\in\mathbb{R}^{n_{1}\times\left(P+1\right)}$,
$\boldsymbol{\xi}_{1}\in\Xi_{1}$, $\omega_{2}\in\mathcal{P}_{2}\left(\Xi_{2}\right)$,

\begin{align}
 & \int_{\Xi_{2}}\boldsymbol{f}_{2}\left(\tilde{\boldsymbol{U}}{}_{2}\left(\boldsymbol{\xi}_{1}\right)\boldsymbol{\psi}_{2}\left(\boldsymbol{\xi}_{2}\right);\hat{\boldsymbol{U}}_{1}\boldsymbol{\Pi}_{2}\left(\boldsymbol{\xi}_{1}\right)\boldsymbol{\psi}_{2}\left(\boldsymbol{\xi}_{2}\right),\tilde{\boldsymbol{\Xi}}{}_{2}\boldsymbol{\psi}_{2}\left(\boldsymbol{\xi}_{2}\right)\right)\omega_{2}\left(\boldsymbol{\xi}_{2}\right)\mu_{2}\left(\boldsymbol{\xi}_{2}\right)d\boldsymbol{\xi}_{2}\nonumber \\
 & =\boldsymbol{0}\\
\Rightarrow & \tilde{\boldsymbol{F}}_{2}\left(\tilde{\boldsymbol{U}}{}_{2};\hat{\boldsymbol{U}}{}_{1},\boldsymbol{\xi}_{1}\right)=\int_{\Xi_{2}}\left(\boldsymbol{f}_{2}\left(\tilde{\boldsymbol{U}}{}_{2}\left(\boldsymbol{\xi}_{1}\right)\boldsymbol{\psi}_{2}\left(\boldsymbol{\xi}_{2}\right);\hat{\boldsymbol{U}}_{1}\boldsymbol{\Pi}_{2}\left(\boldsymbol{\xi}_{1}\right)\boldsymbol{\psi}_{2}\left(\boldsymbol{\xi}_{2}\right),\tilde{\boldsymbol{\Xi}}{}_{2}\boldsymbol{\psi}_{2}\left(\boldsymbol{\xi}_{2}\right)\right)\right.\nonumber \\
 & \times\left.\boldsymbol{\psi}_{2}\left(\boldsymbol{\xi}_{2}\right)^{\mathbf{T}}\right)\mu_{2}\left(\boldsymbol{\xi}_{2}\right)d\boldsymbol{\xi}_{2}=\boldsymbol{0}.
\end{align}
Here, $\tilde{\boldsymbol{F}}_{i}\equiv\tilde{\boldsymbol{F}}_{i}^{p}$
and $\tilde{\boldsymbol{\Xi}}_{i}\equiv\tilde{\boldsymbol{\Xi}}_{i}^{p}:$
\begin{align}
\tilde{\boldsymbol{\Xi}}_{i} & =\int_{\Xi}\boldsymbol{\xi}_{i}\boldsymbol{\psi}_{i}\left(\boldsymbol{\xi}_{i}\right)^{\mathbf{T}}\mu_{i}\left(\boldsymbol{\xi}_{i}\right)d\boldsymbol{\xi}_{i}
\end{align}
denote the modular gPC coefficient matrices of component residual
$\boldsymbol{f}_{i}$ and random input vector $\boldsymbol{\xi}_{i}$
respectively. Alternatively, a transposed modular Galerkin form of
component $i$ in the coupled PDE system can be derived by projecting
it in $\mathcal{P}_{i}\left(\Xi_{i}\right)$, and subsequently applying
the respective spatial discretization scheme.

Let $\tilde{\boldsymbol{M}}_{1}\equiv\tilde{\boldsymbol{M}}_{1}^{p}$
and $\tilde{\boldsymbol{M}}_{2}\equiv\tilde{\boldsymbol{M}}_{2}^{p}$
denote the respective module operators used in the solvers for $\tilde{\boldsymbol{F}}_{1}=\boldsymbol{0}$
and $\tilde{\boldsymbol{F}}_{2}=\boldsymbol{0}$. Therefore, $\forall\hat{\boldsymbol{U}}_{i}\in\mathbb{R}^{n_{i}\times\left(P+1\right)},\boldsymbol{\xi}_{i}\in\Xi_{i}$,
the iterations in the respective solvers can be formulated as follows.
\begin{align}
\tilde{\boldsymbol{U}}{}_{1}^{\ell+1}\left(\boldsymbol{\xi}{}_{2}\right) & =\tilde{\boldsymbol{M}}_{1}\left(\tilde{\boldsymbol{U}}{}_{1}^{\ell}\left(\boldsymbol{\xi}{}_{2}\right);\hat{\boldsymbol{U}}_{2}^{\ell}\boldsymbol{\Pi}_{1}\left(\boldsymbol{\xi}_{2}\right)\right),\nonumber \\
\tilde{\boldsymbol{U}}{}_{2}^{\ell+1}\left(\boldsymbol{\xi}{}_{1}\right) & =\tilde{\boldsymbol{M}}_{2}\left(\tilde{\boldsymbol{U}}_{2}^{\ell}\left(\boldsymbol{\xi}{}_{1}\right);\hat{\boldsymbol{U}}_{1}^{\ell+1}\boldsymbol{\Pi}_{2}\left(\boldsymbol{\xi}_{1}\right)\right).
\end{align}

From \hyperlink{eq36}{Eq. 3.6}, the global gPC coefficient matrices $\hat{\boldsymbol{U}}_{1}$
and $\hat{\boldsymbol{U}}_{2}$ can be computed iteratively using
a BGS method wrapped around the module operators, as follows.
\hypertarget{eq314}{}
\begin{align}
\hat{\boldsymbol{U}}_{1}^{\ell+1} & =\int_{\Xi_{2}}\tilde{\boldsymbol{M}}_{1}\left(\hat{\boldsymbol{U}}_{1}^{\ell}\boldsymbol{\Pi}_{1}\left(\boldsymbol{\xi}_{2}\right);\hat{\boldsymbol{U}}_{2}^{\ell}\boldsymbol{\Pi}_{1}\left(\boldsymbol{\xi}_{2}\right)\right)\boldsymbol{\Pi}_{1}\left(\boldsymbol{\xi}_{2}\right)^{\mathbf{T}}\mu_{2}\left(\boldsymbol{\xi}_{2}\right)d\boldsymbol{\xi}_{2},\nonumber \\
\hat{\boldsymbol{U}}_{2}^{\ell+1} & =\int_{\Xi_{1}}\tilde{\boldsymbol{M}}_{2}\left(\hat{\boldsymbol{U}}_{2}^{\ell}\mathbf{\boldsymbol{\Pi}}_{2}\left(\boldsymbol{\xi}_{1}\right);\hat{\boldsymbol{U}}_{1}^{\ell+1}\boldsymbol{\Pi}{}_{2}\left(\boldsymbol{\xi}_{1}\right)\right)\boldsymbol{\Pi}_{2}\left(\boldsymbol{\xi}_{1}\right)^{\mathbf{T}}\mu_{1}\left(\boldsymbol{\xi}_{1}\right)d\boldsymbol{\xi}_{1}.
\end{align}

Let $\left\{ \left(\boldsymbol{\xi}_{i}^{\left(j\right)},w_{i}^{\left(j\right)}\right)\right\} _{j=1}^{Q_{i}}$
denote a quadrature rule with level $q\geq p$, in $\Xi_{i}$,
such that all polynomial functions with degree $\leq2q+1$ can be
numerically integrated up to machine precision. Subsequently, the
integrals in \hyperlink{eq314}{Eq. 3.14} can be approximated as follows.
\begin{align}
\hat{\boldsymbol{U}}_{1}^{\ell+1} & \approx\sum_{j=1}^{Q_{2}}w_{2}^{\left(j\right)}\tilde{\boldsymbol{M}}_{1}\left(\hat{\boldsymbol{U}}_{1}^{\ell}\boldsymbol{\Pi}{}_{1}\left(\boldsymbol{\xi}_{2}^{\left(j\right)}\right);\hat{\boldsymbol{U}}_{2}^{\ell}\boldsymbol{\Pi}_{1}\left(\boldsymbol{\xi}_{2}^{\left(j\right)}\right)\right)\boldsymbol{\Pi}_{1}\left(\boldsymbol{\xi}_{2}^{\left(j\right)}\right){}^{\mathbf{T}},\nonumber \\
\hat{\boldsymbol{U}}_{2}^{\ell+1} & \approx\sum_{j=1}^{Q_{1}}w_{1}^{\left(j\right)}\tilde{\boldsymbol{M}}_{2}\left(\hat{\boldsymbol{U}}_{2}^{\ell}\boldsymbol{\Pi}_{2}\left(\boldsymbol{\xi}_{1}^{\left(j\right)}\right);\hat{\boldsymbol{U}}_{1}^{\ell+1}\boldsymbol{\Pi}_{2}\left(\boldsymbol{\xi}_{1}^{\left(j\right)}\right)\right)\boldsymbol{\Pi}_{2}\left(\boldsymbol{\xi}_{1}^{\left(j\right)}\right)^{\mathbf{T}}.
\end{align}

As opposed to the standard ISP method, which is based on solving the
global Galerkin form, these module operators can be developed independently
and modified without affecting each other. Moreover, since the coordinate
directions in $\Xi_{i}$ are assumed to be independent, the quadrature
rules can be constructed by tensorization of univariate quadrature
rules in each direction, such that $Q_{i}=\left(q+1\right)^{s_{i}}\geq\left(p+1\right)^{s_{i}}.$Alternatively,
sparse grid tensorization methods {[}\hyperlink{ref28}{28}{]} exhibit a significantly
slower, albeit exponential growth in the size of the quadrature rules
with respect to the dimension $s_{i}$ and level $q$. 

Assuming that the repeated execution of module operators $\tilde{\boldsymbol{M}}_{1}$
and $\tilde{\boldsymbol{M}}_{2}$ dominate the computational costs,
the overall cost of this method would therefore be $\approx\mathcal{O}\left(\bar{\mathcal{C}}_{1}P_{1}^{\alpha_{1}}Q_{2}+\bar{\mathcal{C}}_{2}P_{2}^{\alpha_{2}}Q_{1}\right)$.
Since $\alpha_{1},\alpha_{2}>1$, and $P_{1},P_{2}<P$, the overall cost
would indeed be lesser than the cost of the standard ISP method.
However, the curse of dimensionality would still persist here, since
the overall costs would grow exponentially with respect to the combined
dimension $s$ and order $p$. To address this limitation, we propose
that reduced approximations of the data that enter their respective
modules be constructed at every iteration. We will now describe two
main components of this construction: (1) the dimension reduction
routine and (2) the order reduction routine.

\subsection{Dimension reduction}

At each iteration $\ell$, let $\boldsymbol{y}_{1}^{\ell}\equiv\left[\boldsymbol{u}_{1}^{\ell};\boldsymbol{u}_{2}^{\ell}\right]:\Xi\rightarrow\mathbb{R}^{n}$
and $\boldsymbol{y}_{2}^{\ell}\equiv\left[\boldsymbol{u}_{2}^{\ell};\boldsymbol{u}_{1}^{\ell+1}\right]:\Xi\rightarrow\mathbb{R}^{n}$
denote the input data for module 1 and 2 respectively, where $n=n_{1}+n_{2}$.
Consequently, the gPC coefficient matrices $\hat{\boldsymbol{Y}}_{1}^{\ell}\equiv\left[\hat{\boldsymbol{U}}_{1}^{\ell};\hat{\boldsymbol{U}}_{2}^{\ell}\right]=\left[\begin{array}{ccc}
\hat{\boldsymbol{y}}_{1,\jmath_{1}\left(0\right),\jmath_{2}\left(0\right)}^{\ell} & \cdots & \hat{\boldsymbol{y}}_{1,\jmath_{1}\left(P\right),\jmath_{2}\left(P\right)}^{\ell}\end{array}\right]\in\mathbb{R}^{n\times\left(P+1\right)}$ and $\hat{\boldsymbol{Y}}_{2}^{\ell}$$\equiv\left[\hat{\boldsymbol{U}}_{2}^{\ell};\hat{\boldsymbol{U}}_{1}^{\ell+1}\right]=\left[\begin{array}{ccc}
\hat{\boldsymbol{y}}_{2,\jmath_{1}\left(0\right),\jmath_{2}\left(0\right)}^{\ell} & \cdots & \hat{\boldsymbol{y}}_{2,\jmath_{1}\left(P\right),\jmath_{2}\left(P\right)}^{\ell}\end{array}\right]\in\mathbb{R}^{n\times\left(P+1\right)}$ define the
respective input data for the module operators $\tilde{\boldsymbol{M}}_{1}$
and $\tilde{\boldsymbol{M}}_{2}$, Moreover, let $\boldsymbol{\Gamma}_{1},\boldsymbol{\Gamma}_{2}:$
\begin{equation}
\boldsymbol{\Gamma}_{1}=\left[\begin{array}{cc}
\boldsymbol{G}_{1} & \boldsymbol{0}\\
\boldsymbol{0} & \boldsymbol{G}_{2}
\end{array}\right]\in\mathbb{R}^{n\times n},\ \boldsymbol{\Gamma}_{2}=\left[\begin{array}{cc}
\boldsymbol{G}_{2} & \boldsymbol{0}\\
\boldsymbol{0} & \boldsymbol{G}_{1}
\end{array}\right]\in\mathbb{R}^{n\times n}
\end{equation}
denote the Gramian matrices corresponding to $\boldsymbol{y}_{1}^{\ell}$
and $\boldsymbol{y}_{2}^{\ell}$ respectively, where $\forall i\in\left\{ 1,2\right\} $,
$P_{i}^{\prime}+1={\displaystyle \frac{\left(p-1+s_{i}\right)!}{\left(p-1\right)!s_{i}!}}$. 

Therefore, using $\hat{\boldsymbol{Y}}_{1}^{\ell}$ and $\hat{\boldsymbol{Y}}_{2}^{\ell}$,
we construct $\boldsymbol{Z}_{1}^{\ell}$ and $\boldsymbol{Z}_{2}^{\ell}$:
\begin{align}
\boldsymbol{Z}_{1}^{\ell} & =\left(\boldsymbol{I}_{P_{1}^{\prime}+1}\otimes\boldsymbol{\Gamma}_{1}^{\frac{1}{2}}\right)\left[\begin{array}{ccc}
\boldsymbol{\zeta}_{1,0,1}^{\ell} & \cdots & \boldsymbol{\zeta}_{1,0,P_{2}}^{\ell}\\
\vdots &  & \vdots\\
\boldsymbol{\zeta}_{1,P_{1}^{\prime},1}^{\ell} & \cdots & \boldsymbol{\zeta}_{1,P_{1}^{\prime},P_{2}}^{\ell}
\end{array}\right]\in\mathbb{R}^{\left(P_{1}^{\prime}+1\right)n\times P_{2}},\nonumber \\
\boldsymbol{Z}_{2}^{\ell} & =\left(\boldsymbol{I}_{P_{2}^{\prime}+1}\otimes\boldsymbol{\Gamma}_{2}^{\frac{1}{2}}\right)\left[\begin{array}{ccc}
\boldsymbol{\zeta}_{2,0,1}^{\ell} & \cdots & \boldsymbol{\zeta}_{2,0,P_{1}}^{\ell}\\
\vdots &  & \vdots\\
\boldsymbol{\zeta}_{2,P_{2}^{\prime},1}^{\ell} & \cdots & \boldsymbol{\zeta}_{2,P_{2}^{\prime},P_{1}}^{\ell}
\end{array}\right]\in\mathbb{R}^{\left(P_{2}^{\prime}+1\right)n\times P_{1}},
\end{align}
where $\forall0\leq j\leq P_{1}^{\prime}$, $1\leq k\leq P_{2}$,
$\boldsymbol{\zeta}_{1,j,k}^{\ell}\in\mathbb{R}^{n}$:
\begin{equation}
\boldsymbol{\zeta}_{1,j,k}^{\ell}=\begin{cases}
\hat{\boldsymbol{y}}_{1,j,k}^{\ell} & \mathrm{deg}\left(\psi_{1,j}\right)+\mathrm{deg}\left(\psi_{2,k}\right)\leq p\\
\boldsymbol{0} & \mathrm{deg}\left(\psi_{1,j}\right)+\mathrm{deg}\left(\psi_{2,k}\right)>p
\end{cases},
\end{equation}
and $\forall0\leq j\leq P_{2}^{\prime}$, $1\leq k\leq P_{1}$, $\boldsymbol{\zeta}_{2,j,k}^{\ell}\in\mathbb{R}^{n}$:

\begin{equation}
\boldsymbol{\zeta}_{2,j,k}^{\ell}=\begin{cases}
\hat{\boldsymbol{y}}_{2,j,k}^{\ell} & \mathrm{deg}\left(\psi_{2,j}\right)+\mathrm{deg}\left(\psi_{1,k}\right)\leq p\\
\boldsymbol{0} & \mathrm{deg}\left(\psi_{2,j}\right)+\mathrm{deg}\left(\psi_{1,k}\right)>p
\end{cases}.
\end{equation}

Subsequently, we compute the singular value decomposition (SVD) of
$\boldsymbol{Z}_{1}^{\ell}$ and $\boldsymbol{Z}_{2}^{\ell}$:
\begin{equation}
\boldsymbol{Z}_{1}^{\ell}=\tilde{\boldsymbol{\Upsilon}}_{1}^{\ell}\boldsymbol{\Sigma}_{1}^{\ell}\left(\tilde{\boldsymbol{\Theta}}_{1}^{\ell}\right)^{\mathbf{T}},\ \boldsymbol{Z}_{2}^{\ell}=\tilde{\boldsymbol{\Upsilon}}_{2}^{\ell}\boldsymbol{\Sigma}_{2}^{\ell}\left(\tilde{\boldsymbol{\Theta}}_{2}^{\ell}\right)^{\mathbf{T}},
\end{equation}
where the decomposition matrices have the following structure and
dimensions.

\begin{align}
\tilde{\boldsymbol{\Upsilon}}_{1}^{\ell} & =\left[\begin{array}{ccc}
\tilde{\boldsymbol{\upsilon}}_{1,0,1}^{\ell} & \cdots & \tilde{\boldsymbol{\upsilon}}_{1,0,P_{2}}^{\ell}\\
\vdots &  & \vdots\\
\tilde{\boldsymbol{\upsilon}}_{1,P_{1}^{\prime},1}^{\ell} & \cdots & \tilde{\boldsymbol{\upsilon}}_{1,P_{1}^{\prime},P_{2}}^{\ell}
\end{array}\right]\in\mathbb{R}^{\left(P_{1}^{\prime}+1\right)n\times P_{2}},\nonumber \\
\tilde{\boldsymbol{\Upsilon}}_{2}^{\ell} & =\left[\begin{array}{ccc}
\tilde{\boldsymbol{\upsilon}}_{2,0,1}^{\ell} & \cdots & \tilde{\boldsymbol{\upsilon}}_{2,0,P_{1}}^{\ell}\\
\vdots &  & \vdots\\
\tilde{\boldsymbol{\upsilon}}_{2,P_{2}^{\prime},1}^{\ell} & \cdots & \tilde{\boldsymbol{\upsilon}}_{2,P_{2}^{\prime},P_{1}}^{\ell}
\end{array}\right]\in\mathbb{R}^{\left(P_{2}^{\prime}+1\right)n\times P_{1}},
\end{align}
\begin{equation}
\boldsymbol{\Sigma}_{1}^{\ell}=\left[\begin{array}{ccc}
\sigma_{1,1}^{\ell}\\
 & \ddots\\
 &  & \sigma_{1,P_{2}}^{\ell}
\end{array}\right]\in\mathbb{R}^{P_{2}\times P_{2}},\ \boldsymbol{\Sigma}_{2}^{\ell}=\left[\begin{array}{ccc}
\sigma_{2,1}^{\ell}\\
 & \ddots\\
 &  & \sigma_{2,P_{1}}^{\ell}
\end{array}\right]\in\mathbb{R}^{P_{1}\times P_{1}},
\end{equation}
\begin{align}
\tilde{\boldsymbol{\Theta}}_{1}^{\ell} & =\left[\begin{array}{ccc}
\hat{\theta}_{1,1,1} & \cdots & \hat{\theta}_{1,1,P_{2}}\\
\vdots &  & \vdots\\
\hat{\theta}_{1,P_{2},1} & \cdots & \hat{\theta}_{1,P_{2},P_{2}}
\end{array}\right]\in\mathbb{R}^{P_{2}\times P_{2}},\nonumber \\
\tilde{\boldsymbol{\Theta}}_{2}^{\ell} & =\left[\begin{array}{ccc}
\hat{\theta}_{2,1,1} & \cdots & \hat{\theta}_{2,1,P_{1}}\\
\vdots &  & \vdots\\
\hat{\theta}_{2,P_{1},1} & \cdots & \hat{\theta}_{2,P_{1},P_{1}}
\end{array}\right]\in\mathbb{R}^{P_{1}\times P_{1}}.
\end{align}

We then construct $\hat{\boldsymbol{\Upsilon}}_{1}^{\ell}\in\mathbb{R}^{\left(P_{1}+1\right)n\times P_{2}}$,
$\hat{\boldsymbol{\Upsilon}}_{2}^{\ell}\in\mathbb{R}^{\left(P_{2}+1\right)n\times P_{1}}$:
\begin{align}
\hat{\boldsymbol{\Upsilon}}_{1}^{\ell} & =\left(\boldsymbol{I}_{P_{1}+1}\otimes\boldsymbol{\Gamma}_{1}^{-\frac{1}{2}}\right)\left[\begin{array}{c}
\tilde{\boldsymbol{\Upsilon}}_{1}^{\ell}\\
\boldsymbol{0}
\end{array}\right]=\left[\begin{array}{ccc}
\hat{\boldsymbol{\upsilon}}_{1,0,1}^{\ell} & \cdots & \hat{\boldsymbol{\upsilon}}_{1,0,P_{2}}^{\ell}\\
\vdots &  & \vdots\\
\hat{\boldsymbol{\upsilon}}_{1,P_{1},1}^{\ell} & \cdots & \hat{\boldsymbol{\upsilon}}_{1,P_{1},P_{2}}^{\ell}
\end{array}\right],\nonumber \\
\hat{\boldsymbol{\Upsilon}}_{2}^{\ell} & =\left(\boldsymbol{I}_{P_{2}^{\prime}+1}\otimes\boldsymbol{\Gamma}_{2}^{-\frac{1}{2}}\right)\left[\begin{array}{c}
\tilde{\boldsymbol{\Upsilon}}_{2}^{\ell}\\
\boldsymbol{0}
\end{array}\right]=\left[\begin{array}{ccc}
\hat{\boldsymbol{\upsilon}}_{2,0,1}^{\ell} & \cdots & \hat{\boldsymbol{\upsilon}}_{2,0,P_{1}}^{\ell}\\
\vdots &  & \vdots\\
\hat{\boldsymbol{\upsilon}}_{2,P_{2},1}^{\ell} & \cdots & \hat{\boldsymbol{\upsilon}}_{2,P_{2},P_{1}}^{\ell}
\end{array}\right],
\end{align}
and $\hat{\boldsymbol{\Theta}}_{1}^{\ell}\in\mathbb{R}^{\left(P_{2}+1\right)\times P_{2}}$,
$\hat{\boldsymbol{\Theta}}_{2}^{\ell}\in\mathbb{R}^{\left(P_{1}+1\right)\times P_{1}}$:
\begin{align}
\hat{\boldsymbol{\Theta}}_{1}^{\ell} & =\left[\begin{array}{ccc}
0 & \cdots & 0\\
 & \tilde{\boldsymbol{\Theta}}_{1}^{\ell}
\end{array}\right]=\left[\begin{array}{ccc}
\hat{\theta}_{1,0,1}^{\ell} & \cdots & \hat{\theta}_{1,0,P_{2}}^{\ell}\\
\vdots &  & \vdots\\
\hat{\theta}_{1,P_{2},1}^{\ell} & \cdots & \hat{\theta}_{1,P_{2},P_{2}}^{\ell}
\end{array}\right]=\left[\begin{array}{c}
\hat{\boldsymbol{\theta}}_{1,1}\\
\vdots\\
\hat{\boldsymbol{\theta}}_{1,P_{2}}
\end{array}\right]^{\mathbf{T}},\nonumber \\
\hat{\boldsymbol{\Theta}}_{2}^{\ell} & =\left[\begin{array}{ccc}
0 & \cdots & 0\\
 & \tilde{\boldsymbol{\Theta}}_{2}^{\ell}
\end{array}\right]=\left[\begin{array}{ccc}
\hat{\theta}_{2,0,1}^{\ell} & \cdots & \hat{\theta}_{2,0,P_{1}}^{\ell}\\
\vdots &  & \vdots\\
\hat{\theta}_{2,P_{1},1}^{\ell} & \cdots & \hat{\theta}_{2,P_{1},P_{1}}^{\ell}
\end{array}\right]=\left[\begin{array}{c}
\hat{\boldsymbol{\theta}}_{2,1}\\
\vdots\\
\hat{\boldsymbol{\theta}}_{2,P_{1}}
\end{array}\right]^{\mathbf{T}}.
\end{align}

Subsequently, the gPC approximations $\boldsymbol{y}_{1}^{\ell,p}$
and $\boldsymbol{y}_{2}^{\ell,p}$ can be rewritten as follows. $\forall\boldsymbol{\xi}_{1}\in\Xi_{1}$,
$\boldsymbol{\xi}_{2}\in\Xi_{2}$,
\hypertarget{eq326}{}
\begin{align}
\boldsymbol{y}_{1}^{\ell,p}\left(\boldsymbol{\xi}_{1},\boldsymbol{\xi}_{2}\right) & =\bar{\boldsymbol{y}}_{1}^{\ell}\left(\boldsymbol{\xi}_{1}\right)+\sum_{j=1}^{P_{2}}\sigma_{1,j}^{\ell}\boldsymbol{\upsilon}_{1,j}^{\ell}\left(\boldsymbol{\xi}_{1}\right)\theta_{1,j}^{\ell}\left(\boldsymbol{\xi}_{2}\right),\nonumber \\
\boldsymbol{y}_{2}^{\ell,p}\left(\boldsymbol{\xi}_{1},\boldsymbol{\xi}_{2}\right) & =\bar{\boldsymbol{y}}_{2}^{\ell}\left(\boldsymbol{\xi}_{2}\right)+\sum_{j=1}^{P_{1}}\sigma_{2,j}^{\ell}\boldsymbol{\upsilon}_{2,j}^{\ell}\left(\boldsymbol{\xi}_{2}\right)\theta_{2,j}^{\ell}\left(\boldsymbol{\xi}_{1}\right),
\end{align}
where
\begin{equation}
\bar{\boldsymbol{y}}_{1}^{\ell}\left(\boldsymbol{\xi}_{1}\right)=\sum_{j=0}^{P_{1}}\hat{\boldsymbol{y}}_{1,j,0}^{\ell}\psi_{1,j}\left(\boldsymbol{\xi}_{1}\right),\ \bar{\boldsymbol{y}}_{2}^{\ell}\left(\boldsymbol{\xi}_{2}\right)=\sum_{j=0}^{P_{2}}\hat{\boldsymbol{y}}_{2,0,j}^{\ell}\psi_{2,j}\left(\boldsymbol{\xi}_{2}\right).
\end{equation}
Moreover, $\forall1\leq j\leq P_{2}$, 
\begin{equation}
\boldsymbol{\upsilon}_{1,j}^{\ell}\left(\boldsymbol{\xi}_{1}\right)=\sum_{k=0}^{P_{1}}\hat{\boldsymbol{\upsilon}}_{1,k,j}^{\ell}\psi_{1,k}\left(\boldsymbol{\xi}_{1}\right),\ \theta_{1,j}^{\ell}\left(\boldsymbol{\xi}_{2}\right)=\sum_{k=0}^{P_{2}}\hat{\theta}_{1,k,j}^{\ell}\psi_{2,k}\left(\boldsymbol{\xi}_{2}\right)=\hat{\boldsymbol{\theta}}_{1,j}\boldsymbol{\psi}_{2}\left(\boldsymbol{\xi}_{2}\right),
\end{equation}
and $\forall1\leq j\leq P_{1}$, 
\begin{equation}
\boldsymbol{\upsilon}_{2,j}^{\ell}\left(\boldsymbol{\xi}_{2}\right)=\sum_{k=0}^{P_{2}}\hat{\boldsymbol{\upsilon}}_{2,k,j}^{\ell}\psi_{2,k}\left(\boldsymbol{\xi}_{2}\right),\ \theta_{2,j}^{\ell}\left(\boldsymbol{\xi}_{1}\right)=\sum_{k=0}^{P_{1}}\hat{\theta}_{2,k,j}^{\ell}\psi_{1,k}\left(\boldsymbol{\xi}_{1}\right)==\hat{\boldsymbol{\theta}}_{2,j}\boldsymbol{\psi}_{1}\left(\boldsymbol{\xi}_{1}\right),
\end{equation}

Each functional representation in \hyperlink{eq326}{Eq. 3.26} is a finite-dimensional
and separable variant of the widely used Karhunen-Loeve (KL) expansion
{[}\hyperlink{ref29}{29}{]}. While the coefficients in a standard KL expansion are deterministic,
the coefficients in the separable KL expansion are random. 

Moreover, the modular
gPC coefficient matrices of $\boldsymbol{y}_{1}^{\ell}$ and $\boldsymbol{y}_{2}^{\ell}$
can be formulated as follows. $\forall\boldsymbol{\xi}_{1}\in\Xi_{1}$,
$\boldsymbol{\xi}_{2}\in\Xi_{2}$, 
\begin{align}
\tilde{\boldsymbol{Y}}_{1}^{\ell}\left(\boldsymbol{\xi}_{2}\right) & =\tilde{\boldsymbol{Y}}_{1,0}^{\ell}+\sum_{j=1}^{P_{2}}\sigma_{1,j}^{\ell}\tilde{\boldsymbol{Y}}_{1,j}^{\ell}\theta_{1,j}^{\ell}\left(\boldsymbol{\xi}_{2}\right),\nonumber \\
\tilde{\boldsymbol{Y}}_{2}^{\ell}\left(\boldsymbol{\xi}_{1}\right) & =\tilde{\boldsymbol{Y}}_{2,0}^{\ell}+\sum_{j=1}^{P_{1}}\sigma_{2,j}^{\ell}\tilde{\boldsymbol{Y}}_{2,j}^{\ell}\theta_{2,j}^{\ell}\left(\boldsymbol{\xi}_{1}\right),
\end{align}
where $\tilde{\boldsymbol{Y}}_{1,0}^{\ell}=\left[\begin{array}{ccc}
\hat{\boldsymbol{y}}_{1,0,0}^{\ell} & \cdots & \hat{\boldsymbol{y}}_{1,P_{1},0}^{\ell}\end{array}\right],$ $\tilde{\boldsymbol{Y}}_{2,0}^{\ell}=\left[\begin{array}{ccc}
\hat{\boldsymbol{y}}_{2,0,0}^{\ell} & \cdots & \hat{\boldsymbol{y}}_{1,0,P_{2}}^{\ell}\end{array}\right]$, $\forall1\leq j\leq P_{2},$ $\tilde{\boldsymbol{Y}}_{1,j}^{\ell}=\left[\begin{array}{ccc}
\hat{\boldsymbol{\upsilon}}_{1,0,j}^{\ell} & \cdots & \hat{\boldsymbol{\upsilon}}_{1,P_{1},j}^{\ell}\end{array}\right],$ and $\forall1\leq j\leq P_{1},$ $\tilde{\boldsymbol{Y}}_{2,j}^{\ell}=\left[\begin{array}{ccc}
\hat{\boldsymbol{\upsilon}}_{2,0,j}^{\ell} & \cdots & \hat{\boldsymbol{\upsilon}}_{2,P_{2},j}^{\ell}\end{array}\right]$.

\subsubsection{Reduced chaos expansions}

At each iteration $\ell$, the expansions in \hyperlink{eq326}{Eq. 3.26} can be truncated
by retaining $d_{1}\equiv d_{1}^{\ell}$ and $d_{2}\equiv d_{2}^{\ell}$
terms, to define the respective reduced chaos expansions {[}\hyperlink{ref16}{16}{]}
$\boldsymbol{y}_{1}^{\ell,d_{1}}\equiv\boldsymbol{y}_{1}^{\ell,p,d_{1}}$
and $\boldsymbol{y}_{2}^{\ell,d_{2}}\equiv\boldsymbol{y}_{2}^{\ell,p,d_{2}}$:
$\forall\boldsymbol{\xi}_{1}\in\Xi_{1}$, $\boldsymbol{\xi}_{2}\in\Xi_{2}$,
\begin{align}
\boldsymbol{y}_{1}^{\ell,p}\left(\boldsymbol{\xi}_{1},\boldsymbol{\xi}_{2}\right)\approx\boldsymbol{y}_{1}^{\ell,d_{1}}\left(\boldsymbol{\xi}_{1},\boldsymbol{\xi}_{2}\right) & =\bar{\boldsymbol{y}}_{1}^{\ell}\left(\boldsymbol{\xi}_{1}\right)+\sum_{j=1}^{d_{1}}\sigma_{1,j}^{\ell}\boldsymbol{\upsilon}_{1,j}^{\ell}\left(\boldsymbol{\xi}_{1}\right)\theta_{1,j}^{\ell}\left(\boldsymbol{\xi}_{2}\right),\nonumber \\
\boldsymbol{y}_{2}^{\ell,p}\left(\boldsymbol{\xi}_{1},\boldsymbol{\xi}_{2}\right)\approx\boldsymbol{y}_{2}^{\ell,d_{2}}\left(\boldsymbol{\xi}_{1},\boldsymbol{\xi}_{2}\right) & =\bar{\boldsymbol{y}}_{2}^{\ell}\left(\boldsymbol{\xi}_{2}\right)+\sum_{j=1}^{d_{2}}\sigma_{2,j}^{\ell}\boldsymbol{\upsilon}_{2,j}^{\ell}\left(\boldsymbol{\xi}_{2}\right)\theta_{2,j}^{\ell}\left(\boldsymbol{\xi}_{1}\right).
\end{align}

Similarly, reduced dimensional approximations of the modular gPC coefficient
matrices $\tilde{\boldsymbol{Y}}_{1}^{\ell}$ and $\tilde{\boldsymbol{Y}}_{2}^{\ell}$
can be defined as follows. $\forall\boldsymbol{\xi}_{1}\in\Xi_{1}$,
$\boldsymbol{\xi}_{2}\in\Xi_{2}$,

\begin{align}
\tilde{\boldsymbol{Y}}_{1}^{\ell}\left(\boldsymbol{\xi}_{2}\right)\approx\tilde{\boldsymbol{Y}}_{1}^{\ell,d_{1}}\left(\boldsymbol{\xi}_{2}\right) & =\tilde{\boldsymbol{Y}}_{1,0}^{\ell}+\sum_{j=1}^{d_{1}}\sigma_{1,j}^{\ell}\tilde{\boldsymbol{Y}}_{1,j}^{\ell}\theta_{1,j}^{\ell}\left(\boldsymbol{\xi}_{2}\right),\nonumber \\
\tilde{\boldsymbol{Y}}_{2}^{\ell}\left(\boldsymbol{\xi}_{1}\right)\approx\tilde{\boldsymbol{Y}}_{2}^{\ell,d_{2}}\left(\boldsymbol{\xi}_{1}\right) & =\tilde{\boldsymbol{Y}}_{2,0}^{\ell}+\sum_{j=1}^{d_{2}}\sigma_{2,j}^{\ell}\tilde{\boldsymbol{Y}}_{2,j}^{\ell}\theta_{2,j}^{\ell}\left(\boldsymbol{\xi}_{1}\right),
\end{align}

Let $\boldsymbol{\theta}_{1}^{\ell}\equiv\left[\begin{array}{ccc}
\theta_{1,1}^{\ell} & \cdots & \theta_{1,d_{1}}^{\ell}\end{array}\right]^{\mathbf{T}}:\Xi_{2}\rightarrow\Theta_{1}^{\ell}\subseteq\mathbb{R}^{d_{1}}$ and $\boldsymbol{\theta}_{2}^{\ell}\equiv\left[\begin{array}{cc}
\theta_{2,1}^{\ell} & \cdots\end{array}\right.$ $\left.\theta_{2,d_{2}}^{\ell}\ \right]^{\mathbf{T}}:\Xi_{1}\rightarrow\Theta_{2}^{\ell}\subseteq\mathbb{R}^{d_{2}}$
denote the reduced dimensional random vectors that define $\boldsymbol{y}_{1}^{\ell,d_{1}}$
and $\boldsymbol{y}_{2}^{\ell,d_{2}}$ respectively, with probability
density functions $\nu_{1}^{\ell}:\Theta_{1}^{\ell}\rightarrow\mathbb{R}^{+}$
and $\nu_{2}^{\ell}:\Theta_{2}^{\ell}\rightarrow\mathbb{R}^{+}$.
Therefore, $\forall i\in\left\{ 1,2\right\} $, there exists an affine
map $\tilde{\boldsymbol{Z}}_{i}^{\ell}:\Theta_{i}^{\ell}\rightarrow\mathbb{R}^{n\times\left(P_{i}+1\right)}$
which defines the composition $\tilde{\boldsymbol{Y}}_{i}^{\ell}=\tilde{\boldsymbol{Z}}_{i}^{\ell}\circ\boldsymbol{\theta}_{i}^{\ell}$,
such that $\forall\boldsymbol{\theta}=\left[\begin{array}{ccc}
\theta_{1} & \cdots & \theta_{d_{i}}\end{array}\right]^{\mathbf{T}}\in\Theta_{i}^{\ell}$, 
\[
\tilde{\boldsymbol{Z}}_{i}^{\ell}\left(\boldsymbol{\theta}\right)=\tilde{\boldsymbol{Z}}_{i}^{\ell}\left(\theta_{1},\ldots,\theta_{d_{i}}\right)=\tilde{\boldsymbol{Z}}_{i,0}^{\ell}+\sum_{j=1}^{d_{i}}\tilde{\boldsymbol{Z}}_{i,j}^{\ell}\theta_{j},
\]
where $\tilde{\boldsymbol{Z}}_{i,0}^{\ell}=\tilde{\boldsymbol{Y}}_{i,0}^{\ell}$
and $\forall1\leq j\leq d_{i}$, $\tilde{\boldsymbol{Z}}_{i,j}^{\ell}=\sigma_{i,j}^{\ell}\tilde{\boldsymbol{Y}}_{i,j}^{\ell}$. 

Subsequently, we extract the following submatrix blocks from $\tilde{\boldsymbol{Z}}_{1,0}^{\ell},\ldots,\tilde{\boldsymbol{Z}}_{1,d_{1}}^{\ell}$
and $\tilde{\boldsymbol{Z}}_{2,0}^{\ell},\ldots,\tilde{\boldsymbol{Z}}_{2,d_{2}}^{\ell}$.
$\forall0\leq j\leq d_{1}$
\begin{equation}
\tilde{\boldsymbol{Z}}_{1,j}^{\ell}=\left[\begin{array}{c}
\tilde{\boldsymbol{U}}_{1,1,j}^{\ell}\\
\tilde{\boldsymbol{U}}_{2,1,j}^{\ell}
\end{array}\right],
\end{equation}
and $\forall0\leq j\leq d_{2}$
\begin{equation}
\tilde{\boldsymbol{Z}}_{2,j}^{\ell}=\left[\begin{array}{c}
\tilde{\boldsymbol{U}}_{2,2,j}^{\ell}\\
\tilde{\boldsymbol{U}}_{1,2,j}^{\ell}
\end{array}\right],
\end{equation}
where $\forall i,j\in\left\{ 1,2\right\} $, $0\leq k\leq d_{i}$,
$\tilde{\boldsymbol{U}}_{i,j,k}^{\ell}\in\mathbb{R}^{n_{i}\times\left(P_{j}+1\right)}$.

\subsubsection{Selecting the reduced dimensions}

We prescribe tolerances $\epsilon_{1,\mathrm{dim}}>0$ and $\epsilon_{2,\mathrm{dim}}>0$,
such that at each iteration $\ell$, $d_{1}$ and $d_{2}$ are respectively
selected as the minimum $k_{1}\in\mathbb{N}$ and $k_{2}\in\mathbb{N}$,
which satisfy the following inequalities.
\hypertarget{eq335}{}
\begin{align}
\sqrt{\sum_{j=k_{1}+1}^{P_{2}}\left(\sigma_{1,j}^{\ell}\right)^{2}} & \leq\epsilon_{1,\mathrm{dim}}\sqrt{\sum_{j=1}^{P_{2}}\left(\sigma_{1,j}^{\ell}\right)^{2}},\nonumber \\
\sqrt{\sum_{j=k_{2}+1}^{P_{1}}\left(\sigma_{2,j}^{\ell}\right)^{2}} & \leq\epsilon_{2,\mathrm{dim}}\sqrt{\sum_{j=1}^{P_{1}}\left(\sigma_{2,j}^{\ell}\right)^{2}}.
\end{align}
If $d_{1}<s_{2}$, then a reduced dimensional approximation of the
input data that enters module operator $\tilde{\boldsymbol{M}}_{1}$
exists, with an approximation error $\mathcal{O}\left(\epsilon_{1,\mathrm{dim}}\right)$.
Similarly, if $d_{2}<s_{1}$, then a reduced dimensional approximation
of the input data that enters module operator $\tilde{\boldsymbol{M}}_{2}$
exists, with an approximation error $\mathcal{O}\left(\epsilon_{2,\mathrm{dim}}\right)$.
\hyperlink{thm1}{Theorem 1} proves these error bounds. 
\subsubsection*{Theorem 1}
\hypertarget{thm1}{}
\emph{$\forall i\in\left\{ 1,2\right\} $ and iteration $\ell$, the approximation
$\boldsymbol{y}_{i}^{\ell,d_{i}}$ satisfies the inequality}
\hypertarget{eq336}{}
\begin{equation}
\sqrt{\int_{\Xi}\left\Vert \boldsymbol{y}_{i}^{\ell}\left(\boldsymbol{\xi}\right)-\boldsymbol{y}_{i}^{\ell,d_{i}}\left(\boldsymbol{\xi}\right)\right\Vert _{\boldsymbol{\Gamma}_{i}}^{2}\mu\left(\boldsymbol{\xi}\right)d\boldsymbol{\xi}}\leq\epsilon_{i,\mathrm{dim}}\sqrt{\int_{\Xi}\left\Vert \boldsymbol{y}_{i}^{\ell}\left(\boldsymbol{\xi}\right)-\tilde{\boldsymbol{y}}_{i}^{\ell}\left(\boldsymbol{\xi}_{i}\left(\boldsymbol{\xi}\right)\right)\right\Vert _{\boldsymbol{\Gamma}_{i}}^{2}\mu\left(\boldsymbol{\xi}\right)d\boldsymbol{\xi}}
\end{equation}
\subsubsection{Proof}
For $i=1$, the square of the left hand side expression in \hyperlink{eq336}{Eq. 3.36}
can be written as follows.

\begin{align}
\hypertarget{eq336}{}
 & \int_{\Xi_{1}\times\Xi_{2}}\left\Vert \boldsymbol{y}_{1}^{\ell}\left(\boldsymbol{\xi}_{1},\boldsymbol{\xi}_{2}\right)-\boldsymbol{y}_{1}^{\ell,d_{i}}\left(\boldsymbol{\xi}_{1},\boldsymbol{\xi}_{2}\right)\right\Vert _{\boldsymbol{\Gamma}_{1}}^{2}\mu_{1}\left(\boldsymbol{\xi}_{1}\right)\mu_{2}\left(\boldsymbol{\xi}_{2}\right)d\boldsymbol{\xi}_{1}d\boldsymbol{\xi}_{2}\nonumber\\
 & =\int_{\Xi_{1}\times\Xi_{2}}\left\Vert \sum_{j=d_{1}+1}^{P_{2}}\sigma_{1,j}^{\ell}\boldsymbol{\upsilon}_{1,j}^{\ell}\left(\boldsymbol{\xi}_{1}\right)\theta_{1,j}^{\ell}\left(\boldsymbol{\xi}_{2}\right)\right\Vert _{\boldsymbol{\Gamma}_{1}}^{2}\mu_{1}\left(\boldsymbol{\xi}_{1}\right)\mu_{2}\left(\boldsymbol{\xi}_{2}\right)d\boldsymbol{\xi}_{1}d\boldsymbol{\xi}_{2}\nonumber\\
 & =\int_{\Xi_{1}\times\Xi_{2}}\left\Vert \sum_{j=d_{1}+1}^{P_{2}}\sigma_{1,j}^{\ell}\tilde{\boldsymbol{Y}}_{1,j}^{\ell}\boldsymbol{\psi}_{1}\left(\boldsymbol{\xi}_{1}\right)\hat{\boldsymbol{\theta}}_{1,j}^{\ell}\boldsymbol{\psi}_{2}\left(\boldsymbol{\xi}_{2}\right)\right\Vert _{\boldsymbol{\Gamma}_{1}}^{2}\mu_{1}\left(\boldsymbol{\xi}_{1}\right)\mu_{2}\left(\boldsymbol{\xi}_{2}\right)d\boldsymbol{\xi}_{1}d\boldsymbol{\xi}_{2}\nonumber
\end{align}
\begin{align}
 & =\int_{\Xi_{1}\times\Xi_{2}}\left(\boldsymbol{\psi}_{2}\left(\boldsymbol{\xi}_{2}\right)^{\mathbf{T}}\left[\begin{array}{c}
\hat{\boldsymbol{\theta}}_{1,d_{1}+1}\\
\vdots\\
\hat{\boldsymbol{\theta}}_{1,P_{2}}
\end{array}\right]^{\mathbf{T}}\left[\begin{array}{ccc}
\sigma_{1,d_{1}+1}^{\ell}\nonumber\\
 & \ddots\\
 &  & \sigma_{1,P_{2}}^{\ell}
\end{array}\right]\right.\nonumber\\
 & \left.\times\left[\begin{array}{ccc}
\tilde{\boldsymbol{\upsilon}}_{1,0,1}^{\ell} & \cdots & \tilde{\boldsymbol{\upsilon}}_{1,0,P_{2}}^{\ell}\\
\vdots &  & \vdots\\
\tilde{\boldsymbol{\upsilon}}_{1,P_{1}^{\prime},1}^{\ell} & \cdots & \tilde{\boldsymbol{\upsilon}}_{1,P_{1}^{\prime},P_{2}}^{\ell}\\
 & \boldsymbol{0}
\end{array}\right]^{\mathbf{T}}\left(\boldsymbol{\psi}_{1}\left(\boldsymbol{\xi}_{1}\right)\otimes\boldsymbol{\Gamma}_{1}^{-\frac{1}{2}}\right)\right)\boldsymbol{\Gamma}_{1}\left(\left(\boldsymbol{\psi}_{1}\left(\boldsymbol{\xi}_{1}\right)^{\mathbf{T}}\otimes\boldsymbol{\Gamma}_{1}^{-\frac{1}{2}}\right)\right.\nonumber\\
 & \left.\times\left[\begin{array}{ccc}
\tilde{\boldsymbol{\upsilon}}_{1,0,1}^{\ell} & \cdots & \tilde{\boldsymbol{\upsilon}}_{1,0,P_{2}}^{\ell}\\
\vdots &  & \vdots\\
\tilde{\boldsymbol{\upsilon}}_{1,P_{1}^{\prime},1}^{\ell} & \cdots & \tilde{\boldsymbol{\upsilon}}_{1,P_{1}^{\prime},P_{2}}^{\ell}\\
 & \boldsymbol{0}
\end{array}\right]\left[\begin{array}{ccc}
\sigma_{1,d_{1}+1}^{\ell}\\
 & \ddots\\
 &  & \sigma_{1,P_{2}}^{\ell}
\end{array}\right]\left[\begin{array}{c}
\hat{\boldsymbol{\theta}}_{1,d_{1}+1}\\
\vdots\\
\hat{\boldsymbol{\theta}}_{1,P_{2}}
\end{array}\right]\right.\nonumber\\
 & \left.\times\boldsymbol{\psi}_{2}\left(\boldsymbol{\xi}_{2}\right)\right)\mu_{1}\left(\boldsymbol{\xi}_{1}\right)\mu_{2}\left(\boldsymbol{\xi}_{2}\right)d\boldsymbol{\xi}_{1}d\boldsymbol{\xi}_{2}\nonumber\\
 & =\mathrm{trace}\left(\left(\int_{\Xi_{2}}\boldsymbol{\psi}_{2}\left(\boldsymbol{\xi}_{2}\right)\boldsymbol{\psi}_{2}\left(\boldsymbol{\xi}_{2}\right)^{\mathbf{T}}\mu_{2}\left(\boldsymbol{\xi}_{2}\right)d\boldsymbol{\xi}_{2}\right)\left[\begin{array}{c}
\hat{\boldsymbol{\theta}}_{1,d_{1}+1}\\
\vdots\\
\hat{\boldsymbol{\theta}}_{1,P_{2}}
\end{array}\right]^{\mathbf{T}}\right.\nonumber\\
 & \left.\times\left[\begin{array}{ccc}
\sigma_{1,d_{1}+1}^{\ell}\\
 & \ddots\\
 &  & \sigma_{1,P_{2}}^{\ell}
\end{array}\right]\left[\begin{array}{ccc}
\tilde{\boldsymbol{\upsilon}}_{1,0,1}^{\ell} & \cdots & \tilde{\boldsymbol{\upsilon}}_{1,0,P_{2}}^{\ell}\\
\vdots &  & \vdots\\
\tilde{\boldsymbol{\upsilon}}_{1,P_{1}^{\prime},1}^{\ell} & \cdots & \tilde{\boldsymbol{\upsilon}}_{1,P_{1}^{\prime},P_{2}}^{\ell}\\
 & \boldsymbol{0}
\end{array}\right]^{\mathbf{T}}\right.\nonumber\\
 & \left.\times\left(\left(\int_{\Xi_{1}}\boldsymbol{\psi}_{1}\left(\boldsymbol{\xi}_{1}\right)\boldsymbol{\psi}_{1}\left(\boldsymbol{\xi}_{1}\right)^{\mathbf{T}}\mu_{1}\left(\boldsymbol{\xi}_{1}\right)d\boldsymbol{\xi}_{1}\right)\otimes\boldsymbol{I}_{n}\right)\left[\begin{array}{ccc}
\tilde{\boldsymbol{\upsilon}}_{1,0,1}^{\ell} & \cdots & \tilde{\boldsymbol{\upsilon}}_{1,0,P_{2}}^{\ell}\\
\vdots &  & \vdots\\
\tilde{\boldsymbol{\upsilon}}_{1,P_{1}^{\prime},1}^{\ell} & \cdots & \tilde{\boldsymbol{\upsilon}}_{1,P_{1}^{\prime},P_{2}}^{\ell}\\
 & \boldsymbol{0}
\end{array}\right]\right.\nonumber\\
 & \left.\times\left[\begin{array}{ccc}
\sigma_{1,d_{1}+1}^{\ell}\\
 & \ddots\\
 &  & \sigma_{1,P_{2}}^{\ell}
\end{array}\right]\left[\begin{array}{c}
\hat{\boldsymbol{\theta}}_{1,d_{1}+1}\\
\vdots\\
\hat{\boldsymbol{\theta}}_{1,P_{2}}
\end{array}\right]\right) =\sum_{j=d_{1}+1}^{P_{2}}\left(\sigma_{1,j}^{\ell}\right)^{2}.
\end{align}

Similarly, for the right hand side expression, we can show that 
\hypertarget{eq338}{}
\begin{equation}
\int_{\Xi_{1}\times\Xi_{2}}\left\Vert \boldsymbol{y}_{1}^{\ell}\left(\boldsymbol{\xi}_{1},\boldsymbol{\xi}_{2}\right)-\tilde{\boldsymbol{y}}_{1}^{\ell}\left(\boldsymbol{\xi}_{1}\right)\right\Vert _{\boldsymbol{\Gamma}_{1}}^{2}\mu_{1}\left(\boldsymbol{\xi}_{1}\right)\mu_{2}\left(\boldsymbol{\xi}_{2}\right)d\boldsymbol{\xi}_{1}d\boldsymbol{\xi}_{2}=\sum_{j=1}^{P_{2}}\left(\sigma_{1,j}^{\ell}\right)^{2}
\end{equation}

By substituting \hyperlink{eq337}{Eq. 3.37} and \hyperlink{eq338}{Eq. 3.38} into \hyperlink{eq335}{Eq. 3.35}, we can arrive
at \hyperlink{eq336}{Eq. 3.36}. Similarly, the inequality can be proven for $i=2$. $\square$

\subsubsection*{Remarks}
For the alternative (general) formulation
of the coupled stochastic system given in \hyperlink{eq212}{Eq. 2.12} is considered,
the dimension reduction routine would be exactly the same as described
here, with the exception that the input data must be respectively
formulated as $\boldsymbol{y}_{1}^{\ell}\equiv\left[\boldsymbol{u}_{1}^{\ell};\boldsymbol{v}_{2}^{\ell}\right]$
and $\boldsymbol{y}_{2}^{\ell}\equiv\left[\boldsymbol{u}_{2}^{\ell};\boldsymbol{v}_{1}^{\ell+1}\right]$.
Moreover, if any of the right hand side summations in \hyperlink{eq335}{Eq. 3.35} are
zero, then a corresponding zero dimensional reduced approximation
of the input data is constructed. We will now describe the order reduction
procedure.

\subsection{Order reduction}

$\forall i\in\left\{ 1,2\right\} $ and iteration $\ell$, an approximation
of the input data in module operator $\tilde{\boldsymbol{M}}_{i}$
can be constructed in the reduced dimensional stochastic space $\Theta_{i}^{\ell}$,
as described in \S3.3.

Consequently, for any square-integrable, vector valued functions $\boldsymbol{u}:\Theta_{i}^{\ell}\rightarrow\mathbb{R}^{n}$,
a reduced gPC approximation $\boldsymbol{u}^{\tilde{p}_{i}}\equiv\boldsymbol{u}^{p,d_{i},\tilde{p}_{i}}$
of order $\tilde{p}_{i}\equiv\tilde{p}_{i}^{\ell}\geq0$ can be formulated
as follows. $\forall\boldsymbol{\theta}\in\Theta_{i}^{\ell}$, 
\hypertarget{eq339}{}
\begin{equation}
\boldsymbol{u}\left(\boldsymbol{\theta}\right)\approx\boldsymbol{u}^{\tilde{p}_{i}}\left(\boldsymbol{\theta}\right)=\sum_{j=0}^{\tilde{P}_{i}}\tilde{\boldsymbol{u}}_{j}\phi_{i,j}^{\ell}\left(\boldsymbol{\theta}\right)=\tilde{\boldsymbol{U}}^{\tilde{p}_{i}}\boldsymbol{\phi}_{i}^{\ell,\tilde{p}_{i}}\left(\boldsymbol{\theta}\right),
\end{equation}
where $\left\{ \phi_{i,j}^{\ell}\equiv\phi_{i,j}^{\ell,p}:\Theta_{i}^{\ell}\rightarrow\mathbb{R}\right\} _{j\geq0}$
denotes the set of distinct $\nu_{i}^{\ell}-$orthonormal polynomials,
$\tilde{\boldsymbol{U}}^{\tilde{p}_{i}}\equiv\tilde{\boldsymbol{U}}^{p,\tilde{p}_{i}}=\left[\begin{array}{ccc}
\tilde{\boldsymbol{u}}_{0} & \cdots & \tilde{\boldsymbol{u}}_{\tilde{P}_{i}}\end{array}\right]\in\mathbb{R}^{n\times\left(\tilde{P}_{i}+1\right)}$ denotes the reduced order gPC coefficient matrix and $\boldsymbol{\phi}_{i}^{\ell,\tilde{p}_{i}}\equiv\boldsymbol{\phi}_{i}^{\ell,p,\tilde{p}_{i}}=\left[\begin{array}{ccc}
\phi_{i,0}^{\ell} & \cdots & \phi_{i,\tilde{P}_{i}}^{\ell}\end{array}\right]^{\mathbf{T}}:\Theta_{i}^{\ell}\rightarrow\mathbb{R}^{\tilde{P}_{i}+1}$ denotes the reduced order basis vector. If the formulation of the
gPC approximation $\boldsymbol{u}^{\tilde{p}_{i}}$ is based on an
isotropic, total-degree truncation of the infinite polynomial series,
then $\tilde{P}_{i}={\displaystyle \frac{\left(\tilde{p}_{i}+d_{i}\right)}{\tilde{p}_{i}!d_{i}!}}$.
Moreover, we assume that $\tilde{p}_{1},\tilde{p}_{2}\leq p$, which
implies that $\tilde{P}_{1}\leq P_{2}$ and $\tilde{P}_{2}\leq P_{1}$. 

Since the coordinate directions in $\Theta_{i}^{\ell}$ are not necessarily
statistically independent, we cannot simply compute the elements of
$\boldsymbol{\phi}_{i}^{\ell,\tilde{p}_{i}}$ as products of univariate
polynomials. Instead, we propose another SVD based numerical construction
method.

\subsubsection{Reduced order basis construction}

$\forall i\in\left\{ 1,2\right\} $, iteration $\ell$ and degree
index $\boldsymbol{\alpha}=\left[\begin{array}{ccc}
\alpha_{1} & \cdots & \alpha_{d_{i}}\end{array}\right]^{\mathbf{T}}\in\mathbb{N}_{0}^{d_{i}}$, let $m_{i,\boldsymbol{\alpha}}^{\ell}:\Theta_{i}^{\ell}\rightarrow\mathbb{R}$
denote the monomial function with $\mathrm{deg}\left(m_{i,\boldsymbol{\alpha}}^{\ell}\right)=\left|\boldsymbol{\alpha}\right|=\alpha_{1}+\cdots+\alpha_{d_{i}}$,
such that $\forall\boldsymbol{\theta}=\left[\begin{array}{ccc}
\theta_{1} & \cdots & \theta_{d_{i}}\end{array}\right]\in\Theta_{i}^{\ell}$, 
\begin{equation}
m_{i,\boldsymbol{\alpha}}^{\ell}=\prod_{j=1}^{d_{i}}\theta_{j}^{\alpha_{j}}.
\end{equation}
The number of such monomials with total degree $\leq\tilde{p}_{i}$
is equal to $\tilde{P}_{i}+1.$ 

Let $\left\{ \boldsymbol{\alpha}_{j}:\left|\boldsymbol{\alpha}_{j}\right|\leq\tilde{p}_{i}\right\} _{j=0}^{\tilde{P}_{i}}$
denote the corresponding set of indices and $\boldsymbol{m}_{i}^{\ell,\tilde{p}_{i}}\equiv\boldsymbol{\pi}_{i}^{\ell,p,\tilde{p}_{i}}=\left[\begin{array}{ccc}
m_{i,\boldsymbol{\alpha}_{0}} & \cdots & m_{i,\boldsymbol{\alpha}_{\tilde{P}_{i}}}\end{array}\right]^{\mathbf{T}}:\Theta_{i}^{\ell}\rightarrow\mathbb{R}^{\tilde{P}_{i}+1}$ denote the monomial vector. In general, the respective basis vector
$\boldsymbol{\phi}_{i}^{\ell,\tilde{p}_{1}}$ can be computed as follows.
$\forall\boldsymbol{\theta}\in\Theta_{i}^{\ell}$, 
\[
\boldsymbol{\phi}_{i}^{\ell,\tilde{p}_{i}}\left(\boldsymbol{\theta}\right)=\left(\boldsymbol{L}_{i}^{\ell,\tilde{p}_{i}}\right)^{-1}\boldsymbol{m}_{i}^{\ell,\tilde{p}_{i}}\left(\boldsymbol{\theta}\right),
\]
where $\boldsymbol{L}_{i}^{\ell,\tilde{p}_{i}}\in\mathbb{R}^{\tilde{P}_{i}+1}$
is the lower triangular matrix which defines the Cholesky factorization
of the Hankel matrix $\boldsymbol{H}_{i}^{\ell,\tilde{p}_{i}}\equiv\boldsymbol{H}_{i}^{\ell,p,\tilde{p}_{i}}=\boldsymbol{L}_{i}^{\ell,\tilde{p}_{i}}\left(\boldsymbol{L}_{i}^{\ell,\tilde{p}_{i}}\right)^{\mathbf{T}}:$
\begin{align}
\boldsymbol{H}_{i}^{\ell,\tilde{p}_{i}} & =\int_{\Theta_{i}^{\ell}}\boldsymbol{m}_{i}^{\ell,\tilde{p}_{i}}\left(\boldsymbol{\theta}\right)\boldsymbol{m}_{i}^{\ell,\tilde{p}_{i}}\left(\boldsymbol{\theta}\right)^{\mathbf{T}}\nu_{i}^{\ell}\left(\boldsymbol{\theta}\right)d\boldsymbol{\theta}=\nonumber \\
 & =\begin{cases}
\int_{\Xi_{2}}\boldsymbol{m}_{1}^{\ell,\tilde{p}_{1}}\left(\boldsymbol{\theta}_{1}^{\ell}\left(\boldsymbol{\xi}_{2}\right)\right)\boldsymbol{m}_{1}^{\ell,\tilde{p}_{1}}\left(\boldsymbol{\theta}_{1}^{\ell}\left(\boldsymbol{\xi}_{2}\right)\right)^{\mathbf{T}}\mu_{2}\left(\boldsymbol{\xi}_{2}\right)d\boldsymbol{\xi}_{2} & i=1\\
\int_{\Xi_{1}}\boldsymbol{m}_{2}^{\ell,\tilde{p}_{2}}\left(\boldsymbol{\theta}_{2}^{\ell}\left(\boldsymbol{\xi}_{1}\right)\right)\boldsymbol{m}_{2}^{\ell,\tilde{p}_{2}}\left(\boldsymbol{\theta}_{2}^{\ell}\left(\boldsymbol{\xi}_{1}\right)\right)^{\mathbf{T}}\mu_{1}\left(\boldsymbol{\xi}_{1}\right)d\boldsymbol{\xi}_{1} & i=2
\end{cases},
\end{align}
which in turn can be approximated using the respective quadrature
rule as follows. $\tilde{\boldsymbol{H}}_{i}^{\ell,\tilde{p}_{i}}\approx\boldsymbol{H}_{i}^{\ell,\tilde{p}_{i}}:$
\begin{equation}
\tilde{\boldsymbol{H}}_{i}^{\ell,\tilde{p}_{i}}=\begin{cases}
{\displaystyle \sum_{j=1}^{Q_{2}}}w_{2}^{\left(j\right)}\boldsymbol{m}_{1}^{\ell,\tilde{p}_{1}}\left(\boldsymbol{\theta}_{1}^{\ell}\left(\boldsymbol{\xi}_{2}^{\left(j\right)}\right)\right)\boldsymbol{m}_{1}^{\ell,\tilde{p}_{1}}\left(\boldsymbol{\theta}_{1}^{\ell}\left(\boldsymbol{\xi}_{2}^{\left(j\right)}\right)\right)^{\mathbf{T}} & i=1\\
{\displaystyle \sum_{j=1}^{Q_{1}}}w_{1}^{\left(j\right)}\boldsymbol{m}_{2}^{\ell,\tilde{p}_{2}}\left(\boldsymbol{\theta}_{2}^{\ell}\left(\boldsymbol{\xi}_{1}^{\left(j\right)}\right)\right)\boldsymbol{m}_{2}^{\ell,\tilde{p}_{2}}\left(\boldsymbol{\theta}_{2}^{\ell}\left(\boldsymbol{\xi}_{1}^{\left(j\right)}\right)\right)^{\mathbf{T}} & i=2
\end{cases}.
\end{equation}
However, the possibility of negative weights in the quadrature rules can lead to negative or zero
eigenvalues in the approximation $\tilde{\boldsymbol{H}}_{i}^{\ell,\tilde{p}_{i}}$, and consequently preclude the Cholesky factorization. Therefore, we instead compute the rank-reduced
SVD of $\tilde{\boldsymbol{H}}_{i}^{\ell,\tilde{p}_{i}}:$
\begin{equation}
\tilde{\boldsymbol{H}}_{i}^{\ell,\tilde{p}_{i}}=\tilde{\boldsymbol{V}}_{i}^{\ell,\tilde{p}_{i}}\tilde{\boldsymbol{S}}_{i}^{\ell,\tilde{p}_{i}}\tilde{\boldsymbol{\Sigma}}_{i}^{\ell,\tilde{p}_{i}}\left(\tilde{\boldsymbol{V}}_{i}^{\ell,\tilde{p}_{i}}\right)^{\mathbf{T}},
\end{equation}
where $\tilde{\boldsymbol{V}}_{i}^{\ell,\tilde{p}_{i}},\tilde{\boldsymbol{\Sigma}}_{i}^{\ell,\tilde{p}_{i}}$
denote the usual decomposition matrices and $\tilde{\boldsymbol{S}}_{i}^{\ell,\tilde{p}_{i}}$
denotes a diagonal matrix with $\pm1$ as its diagonal elements. Since
the approximation $\tilde{\boldsymbol{H}}_{i}^{\ell,\tilde{p}_{i}}$
is symmetric, such a decomposition will always exist {[}\hyperlink{ref30}{30}{]}. 

Subsequently, the basis vector $\boldsymbol{\phi}_{i}^{\ell,\tilde{p}_{1}}$
is formulated as follows. $\forall\boldsymbol{\theta}\in\Theta_{i}^{\ell}$,
\begin{equation}
\boldsymbol{\phi}_{i}^{\ell,\tilde{p}_{1}}\left(\boldsymbol{\theta}\right)=\left(\tilde{\boldsymbol{\Sigma}}_{i}^{\ell,\tilde{p}_{i}}\right)^{-\frac{1}{2}}\left(\tilde{\boldsymbol{V}}_{i}^{\ell,\tilde{p}_{i}}\right)^{\mathbf{T}}\boldsymbol{m}_{i}^{\ell,\tilde{p}_{i}}\left(\boldsymbol{\theta}\right).
\end{equation}
Here, $\boldsymbol{\phi}_{1}^{\ell,\tilde{p}_{1}}$ and $\boldsymbol{\phi}_{2}^{\ell,\tilde{p}_{2}}$
satisfy the discrete orthogonality condition
\begin{align*}
{\displaystyle \sum_{j=1}^{Q_{2}}}w_{2}^{\left(j\right)}\boldsymbol{\phi}_{1}^{\ell,\tilde{p}_{1}}\left(\boldsymbol{\theta}_{1}^{\ell}\left(\boldsymbol{\xi}_{2}^{\left(j\right)}\right)\right)\boldsymbol{\phi}_{1}^{\ell,\tilde{p}_{1}}\left(\boldsymbol{\theta}_{1}^{\ell}\left(\boldsymbol{\xi}_{2}^{\left(j\right)}\right)\right)^{\mathbf{T}} & =\boldsymbol{S}_{1}^{\ell,\tilde{p}_{1}},\\
{\displaystyle \sum_{j=1}^{Q_{1}}}w_{1}^{\left(j\right)}\boldsymbol{\pi}_{2}^{\ell,\tilde{p}_{2}}\left(\boldsymbol{\theta}_{2}^{\ell}\left(\boldsymbol{\xi}_{1}^{\left(j\right)}\right)\right)\boldsymbol{\pi}_{2}^{\ell,\tilde{p}_{2}}\left(\boldsymbol{\theta}_{2}^{\ell}\left(\boldsymbol{\xi}_{1}^{\left(j\right)}\right)\right)^{\mathbf{T}} & =\boldsymbol{S}_{2}^{\ell,\tilde{p}_{2}}.
\end{align*}
Therefore, the reduced order gPC coefficient matrix can be approximated
as follows.
\hypertarget{eq345}{}
\begin{align}
\tilde{\boldsymbol{U}}^{\tilde{p}_{i}} & =\int_{\Theta_{i}^{\ell}}\boldsymbol{u}\left(\boldsymbol{\theta}\right)\boldsymbol{\phi}_{i}^{\ell,\tilde{p}_{1}}\left(\boldsymbol{\theta}\right)^{\mathbf{T}}\boldsymbol{S}_{i}^{\ell,\tilde{p}_{1}}\nu_{i}^{\ell}\left(\boldsymbol{\theta}\right)d\boldsymbol{\theta}\nonumber \\
 & =\begin{cases}
\int_{\Xi_{2}}\boldsymbol{u}\left(\boldsymbol{\theta}_{1}^{\ell}\left(\boldsymbol{\xi}_{2}\right)\right)\boldsymbol{\phi}_{1}^{\ell,\tilde{p}_{1}}\left(\boldsymbol{\theta}_{1}^{\ell}\left(\boldsymbol{\xi}_{2}\right)\right)^{\mathbf{T}}\boldsymbol{S}_{1}^{\ell,\tilde{p}_{1}}\mu_{2}\left(\boldsymbol{\xi}_{2}\right)d\boldsymbol{\xi}_{2} & i=1\\
\int_{\Xi_{1}}\boldsymbol{u}\left(\boldsymbol{\theta}_{2}^{\ell}\left(\boldsymbol{\xi}_{1}\right)\right)\boldsymbol{\phi}_{2}^{\ell,\tilde{p}_{2}}\left(\boldsymbol{\theta}_{2}^{\ell}\left(\boldsymbol{\xi}_{1}\right)\right)^{\mathbf{T}}\boldsymbol{S}_{2}^{\ell,\tilde{p}_{2}}\mu_{1}\left(\boldsymbol{\xi}_{1}\right)d\boldsymbol{\xi}_{1} & i=2
\end{cases}\nonumber \\
 & \approx\begin{cases}
{\displaystyle \sum_{j=1}^{Q_{2}}}w_{2}^{\left(j\right)}\boldsymbol{u}\left(\boldsymbol{\theta}_{1}^{\ell}\left(\boldsymbol{\xi}_{2}^{\left(j\right)}\right)\right)\boldsymbol{\phi}_{1}^{\ell,\tilde{p}_{1}}\left(\boldsymbol{\theta}_{1}^{\ell}\left(\boldsymbol{\xi}_{2}^{\left(j\right)}\right)\right)^{\mathbf{T}}\boldsymbol{S}_{1}^{\ell,\tilde{p}_{1}} & i=1\\
{\displaystyle \sum_{j=1}^{Q_{1}}}w_{1}^{\left(j\right)}\boldsymbol{u}\left(\boldsymbol{\theta}_{2}^{\ell}\left(\boldsymbol{\xi}_{1}^{\left(j\right)}\right)\right)\boldsymbol{\phi}_{2}^{\ell,\tilde{p}_{2}}\left(\boldsymbol{\theta}_{2}^{\ell}\left(\boldsymbol{\xi}_{1}^{\left(j\right)}\right)\right)^{\mathbf{T}}\boldsymbol{S}_{2}^{\ell,\tilde{p}_{2}} & i=2
\end{cases}.
\end{align}

\hyperlink{eq345}{Eq. 3.45} implies that $\left\{ \left(\boldsymbol{\theta}_{1}^{\left(j\right)}=\boldsymbol{\theta}_{1}^{\ell}\left(\boldsymbol{\xi}_{2}^{\left(j\right)}\right),w_{2}^{\left(j\right)}\right)\right\} _{j=1}^{Q_{2}}$
is a quadrature rule in $\Theta_{1}^{\ell}$ with level $\geq\tilde{p}_{1}$
and $\left\{ \left(\boldsymbol{\theta}_{2}^{\left(j\right)}=\boldsymbol{\theta}_{2}^{\ell}\left(\boldsymbol{\xi}_{1}^{\left(j\right)}\right),w_{1}^{\left(j\right)}\right)\right\} _{j=1}^{Q_{1}}$
is a quadrature rule in $\Theta_{2}^{\ell}$ with level $\geq\tilde{p}_{2}$.
However, since $Q_{1}$ grows exponentially with respect to $s_{1}$
and $p$, and $Q_{2}$ grows exponentially with respect to $s_{2}$
and $p$, these quadrature rules are not optimal with respect to the
reduced dimensions $d_{1},d_{2}$ and reduced orders $\tilde{p}_{1},\tilde{p}_{2}$.
Therefore, we propose a QR factorization based method to construct
the computationally optimal quadrature rules. $ $

\subsubsection{Optimal quadrature rule construction}

At each iteration $\ell$, the respective optimally sparse quadrature
rules in $\Theta_{1}^{\ell}$ and $\Theta_{2}^{\ell}$, denoted as
$\left\{ \left(\boldsymbol{\theta}_{1}^{\left(j\right)},\tilde{w}_{2}^{\left(j\right)}\right)\right\} _{j=1}^{Q_{2}}$
and $\left\{ \left(\boldsymbol{\theta}_{2}^{\left(j\right)},\tilde{w}_{1}^{\left(j\right)}\right)\right\} _{j=1}^{Q_{1}}$,
with levels $\tilde{p}_{1}$ and $\tilde{p}_{2}$, would contain the
minimum possible number of non-zero weights and still be able to numerically
integrate all polynomials with total degree $\leq2\tilde{p}_{1}$
and $\leq2\tilde{p}_{2}$, up to machine precision. This requirement
is dictated by the projection formula in \hyperlink{eq345}{Eq. 3.45}. 

$\forall i\in\left\{ 1,2\right\} $, let $\boldsymbol{w}_{i}=\left[\begin{array}{ccc}
w_{i}^{\left(1\right)} & \cdots & w_{i}^{\left(Q_{i}\right)}\end{array}\right]^{\mathbf{T}}\in\mathbb{R}^{Q_{i}}$ denote the dense weight vector and and $\tilde{\boldsymbol{w}}_{i}=\left[\begin{array}{ccc}
w_{i}^{\left(1\right)} & \cdots & w_{i}^{\left(Q_{i}\right)}\end{array}\right]^{\mathbf{T}}\in\mathbb{R}^{Q_{i}}$ denote the optimally sparse weight vector. Therefore, $\tilde{\boldsymbol{w}}_{1}$
and $\tilde{\boldsymbol{w}}_{2}$ each solve an $\ell_{0}$-minimization
problem as follows.
\begin{align}
\tilde{\boldsymbol{w}}_{1} & =\arg\min_{\boldsymbol{\omega}\in\mathbb{R}^{Q_{1}}}\left\Vert \boldsymbol{\omega}\right\Vert _{0}:\boldsymbol{M}_{2}^{\ell,2\tilde{p}_{2}}\boldsymbol{\omega}=\boldsymbol{M}_{2}^{\ell,2\tilde{p}_{2}}\boldsymbol{w}_{1},\nonumber \\
\tilde{\boldsymbol{w}}_{2} & =\arg\min_{\boldsymbol{\omega}\in\mathbb{R}^{Q_{2}}}\left\Vert \boldsymbol{\omega}\right\Vert _{0}:\boldsymbol{M}_{1}^{\ell,2\tilde{p}_{1}}\boldsymbol{\omega}=\boldsymbol{M}_{1}^{\ell,2\tilde{p}_{1}}\boldsymbol{w}_{2},
\end{align}
where 
\begin{align}
\boldsymbol{M}_{1}^{\ell,2\tilde{p}_{1}}\equiv\boldsymbol{M}_{1}^{\ell,p,2\tilde{p}_{1}} & =\left[\begin{array}{ccc}
\boldsymbol{m}_{1}^{\ell,2\tilde{p}_{1}}\left(\boldsymbol{\theta}_{1}^{\left(1\right)}\right) & \cdots & \boldsymbol{m}_{1}^{\ell,2\tilde{p}_{1}}\left(\boldsymbol{\theta}_{1}^{\left(Q_{2}\right)}\right)\end{array}\right]\in\mathbb{R}^{\left(\tilde{N}_{1}+1\right)\times Q_{2}},\nonumber \\
\boldsymbol{M}_{2}^{\ell,2\tilde{p}_{2}}\equiv\boldsymbol{M}_{2}^{\ell,p,2\tilde{p}_{2}} & =\left[\begin{array}{ccc}
\boldsymbol{m}_{2}^{\ell,2\tilde{p}_{2}}\left(\boldsymbol{\theta}_{2}^{\left(1\right)}\right) & \cdots & \boldsymbol{m}_{2}^{\ell,2\tilde{p}_{1}}\left(\boldsymbol{\theta}_{2}^{\left(Q_{1}\right)}\right)\end{array}\right]\in\mathbb{R}^{\left(\tilde{N}_{2}+1\right)\times Q_{1}}
\end{align}
denote the corresponding respective Vandermonde matrices with $\tilde{N}_{1}+1={\displaystyle \frac{\left(2\tilde{p}_{1}+d_{1}\right)}{\left(2\tilde{p}_{1}\right)!d_{1}!}}$
and $\tilde{N}_{2}+1={\displaystyle \frac{\left(2\tilde{p}_{2}+d_{2}\right)}{\left(2\tilde{p}_{2}\right)!d_{2}!}}$
rows respectively. In general, these matrices may not be full rank
and therefore, we can setup numerically stable and equivalent $\ell_{0}$-minimization
problems for $\tilde{\boldsymbol{w}}_{1}$ and $\tilde{\boldsymbol{w}}_{2}$
as follows.
\begin{align}
\tilde{\boldsymbol{w}}_{1} & =\arg\min_{\boldsymbol{\omega}\in\mathbb{R}^{Q_{1}}}\left\Vert \boldsymbol{\omega}\right\Vert _{0}:\left(\boldsymbol{Q}_{2,r_{2}}^{\ell,2\tilde{p}_{2}}\right)^{\mathbf{T}}\boldsymbol{\omega}=\left(\boldsymbol{Q}_{2,r_{2}}^{\ell,2\tilde{p}_{2}}\right)^{\mathbf{T}}\boldsymbol{w}_{1},\nonumber \\
\tilde{\boldsymbol{w}}_{2} & =\arg\min_{\boldsymbol{\omega}\in\mathbb{R}^{Q_{2}}}\left\Vert \boldsymbol{\omega}\right\Vert _{0}:\left(\boldsymbol{Q}_{1,r_{1}}^{\ell,2\tilde{p}_{1}}\right)^{\mathbf{T}}\boldsymbol{\omega}=\left(\boldsymbol{Q}_{1,r_{1}}^{\ell,2\tilde{p}_{1}}\right)^{\mathbf{T}}\boldsymbol{w}_{2},
\end{align}
where $\forall i\in\left\{ 1,2\right\} $, $r_{i}$ denotes the rank
of $\boldsymbol{M}_{i}^{\ell,2\tilde{p}_{i}}$ and $\boldsymbol{Q}_{i,r_{i}}^{\ell,2\tilde{p}_{i}}$
defines the pivoted-QR factorization of $\left(\boldsymbol{M}_{i}^{\ell,2\tilde{p}_{i}}\right)^{\mathbf{T}}$
:
\begin{equation}
\left(\boldsymbol{M}_{i}^{\ell,2\tilde{p}_{i}}\right)^{\mathbf{T}}\boldsymbol{\Pi}_{i}^{\ell,2\tilde{p}_{i}}=\boldsymbol{Q}_{i,r_{i}}^{\ell,2\tilde{p}_{i}}\boldsymbol{R}_{i}^{\ell,2\tilde{p}_{i}}.
\end{equation}

Instead of solving the NP-hard $\ell_{0}$-minimization problems,
we employ a direct approach to compute 'weakly' optimal quadrature
rules, with at most $r_{2}$ and $r_{1}$ non-zero weights respectively.
Therefore, $\forall i\in\left\{ 1,2\right\} $, we compute the pivoted-QR
factorization of $\left(\boldsymbol{Q}_{i,r_{i}}^{\ell,2\tilde{p}_{i}}\right)^{\mathbf{T}}$
:
\begin{equation}
\left(\boldsymbol{Q}_{i,r_{i}}^{\ell,2\tilde{p}_{i}}\right)^{\mathbf{T}}\tilde{\boldsymbol{\Pi}}_{i}^{\ell,2\tilde{p}_{i}}=\tilde{\boldsymbol{Q}}_{i}^{\ell,2\tilde{p}_{i}}\tilde{\boldsymbol{R}}_{i}^{\ell,2\tilde{p}_{i}},
\end{equation}
construct the upper triangular square matrix $\tilde{\boldsymbol{R}}_{i,r_{i}}^{\ell,2\tilde{p}_{i}}$
using the first $r_{i}$ columns of $\tilde{\boldsymbol{R}}_{i}^{\ell,2\tilde{p}_{i}}$,
and compute the sparse weight vectors as follows. 
\begin{align}
\tilde{\boldsymbol{w}}_{1} & =\tilde{\boldsymbol{\Pi}}_{2}^{\ell,2\tilde{p}_{2}}\left[\begin{array}{c}
\left(\tilde{\boldsymbol{R}}_{2,r_{2}}^{\ell,2\tilde{p}_{2}}\right)^{-1}\tilde{\boldsymbol{R}}_{2}^{\ell,2\tilde{p}_{2}}\left(\tilde{\boldsymbol{\Pi}}_{2}^{\ell,2\tilde{p}_{2}}\right)^{\mathbf{T}}\boldsymbol{w}_{1}\\
\boldsymbol{0}
\end{array}\right],\nonumber \\
\tilde{\boldsymbol{w}}_{2} & =\tilde{\boldsymbol{\Pi}}_{1}^{\ell,2\tilde{p}_{1}}\left[\begin{array}{c}
\left(\tilde{\boldsymbol{R}}_{1,r_{1}}^{\ell,2\tilde{p}_{1}}\right)^{-1}\tilde{\boldsymbol{R}}_{1}^{\ell,2\tilde{p}_{1}}\left(\tilde{\boldsymbol{\Pi}}_{1}^{\ell,2\tilde{p}_{1}}\right)^{\mathbf{T}}\boldsymbol{w}_{2}\\
\boldsymbol{0}
\end{array}\right].
\end{align}

$\forall i\in\left\{ 1,2\right\} $, if $\mathcal{Z}_{i}=\left\{ 1\leq j\leq Q_{i}:\tilde{w}_{i}^{\left(i\right)}\neq0\right\} $
denotes the index set corresponding to the non-zero elements in $\tilde{\boldsymbol{w}}_{i}$,
then an efficient reduced projection formula for computing the reduced
order gPC coefficient matrix $\tilde{\boldsymbol{U}}$ can be formulated
as follows. 
\hypertarget{eq352}{}
\begin{equation}
\tilde{\boldsymbol{U}}^{\tilde{p}_{i}}\approx\begin{cases}
{\displaystyle \sum_{j\in\mathcal{Z}_{2}}}\tilde{w}_{2}^{\left(j\right)}\boldsymbol{u}\left(\boldsymbol{\theta}_{1}^{\left(j\right)}\right)\boldsymbol{\phi}_{1}^{\ell,\tilde{p}_{1}}\left(\boldsymbol{\theta}_{1}^{\left(j\right)}\right)^{\mathbf{T}}\boldsymbol{S}_{1}^{\ell,\tilde{p}_{1}} & i=1\\
{\displaystyle \sum_{j\in\mathcal{Z}_{1}}}\tilde{w}_{1}^{\left(j\right)}\boldsymbol{u}\left(\boldsymbol{\theta}_{2}^{\left(j\right)}\right)\boldsymbol{\phi}_{2}^{\ell,\tilde{p}_{2}}\left(\boldsymbol{\theta}_{2}^{\left(j\right)}\right)^{\mathbf{T}}\boldsymbol{S}_{2}^{\ell,\tilde{p}_{2}} & i=2
\end{cases}.
\end{equation}

\subsubsection{Approximating modular and global gPC coefficients}

$\forall i\in\left\{ 1,2\right\} $ and iteration $\ell$, let $\tilde{\boldsymbol{U}}^{d_{i}}\equiv\tilde{\boldsymbol{U}}^{p,d_{i}}:\Theta_{i}^{\ell}\rightarrow\mathbb{R}^{n\times\left(P_{i}+1\right)}$
denote the reduced dimensional modular gPC coefficient matrix of $\boldsymbol{u}$,
corresponding to module $i$. Consequently, for the original modular
gPC coefficient matrix $\tilde{\boldsymbol{U}}^{p}$, defined \S3.1,
a composition $\tilde{\boldsymbol{U}}^{p}=\tilde{\boldsymbol{U}}^{d_{i}}\circ\boldsymbol{\theta}_{i}^{\ell}$
exists. Therefore, using the reduced projection formula (\hyperlink{eq352}{Eq. 3.52}),
the reduced dimensional global gPC coefficient matrix $\hat{\boldsymbol{U}}^{d_{i}}\equiv\hat{\boldsymbol{U}}^{p,d_{i}}$
can be computed as follows.
\hypertarget{eq353}{}
\begin{equation}
\hat{\boldsymbol{U}}^{d_{i}}\approx\begin{cases}
{\displaystyle \sum_{j\in\mathcal{Z}_{2}}}\tilde{w}_{2}^{\left(j\right)}\tilde{\boldsymbol{U}}^{d_{1}}\left(\boldsymbol{\theta}_{1}^{\left(j\right)}\right)\left(\left(\boldsymbol{\phi}_{1}^{\ell,\tilde{p}_{1}}\left(\boldsymbol{\theta}_{1}^{\left(j\right)}\right)^{\mathbf{T}}\boldsymbol{S}_{1}^{\ell,\tilde{p}_{1}}\right)\otimes\boldsymbol{I}_{P_{1}+1}\right) & i=1\\
{\displaystyle \sum_{j\in\mathcal{Z}_{1}}}\tilde{w}_{1}^{\left(j\right)}\tilde{\boldsymbol{U}}^{d_{2}}\left(\boldsymbol{\theta}_{2}^{\left(j\right)}\right)\left(\left(\boldsymbol{\phi}_{2}^{\ell,\tilde{p}_{2}}\left(\boldsymbol{\theta}_{2}^{\left(j\right)}\right)^{\mathbf{T}}\boldsymbol{S}_{2}^{\ell,\tilde{p}_{2}}\right)\otimes\boldsymbol{I}_{P_{2}+1}\right) & i=2
\end{cases}.
\end{equation}

Subsequently, using \hyperlink{eq353}{Eq. 3.53}, the global gPC coefficient matrix $\hat{\boldsymbol{U}}^{p,d_{i},\tilde{p}_{i}}$
can be computed as follows.
\begin{equation}
\hat{\boldsymbol{U}}^{p,d_{i},\tilde{p}_{i}}\approx\begin{cases}
{\displaystyle \sum_{j=1}^{Q_{2}}}w_{2}^{\left(j\right)}\hat{\boldsymbol{U}}^{d_{1}}\left(\boldsymbol{\phi}_{1}^{\ell,\tilde{p}_{1}}\otimes\boldsymbol{I}_{P_{1}+1}\right)\boldsymbol{\Pi}_{1}\left(\boldsymbol{\xi}_{2}^{\left(j\right)}\right)^{\mathbf{T}} & i=1\\
{\displaystyle \sum_{j=1}^{Q_{1}}}w_{1}^{\left(j\right)}\hat{\boldsymbol{U}}^{d_{2}}\left(\boldsymbol{\phi}_{2}^{\ell,\tilde{p}_{2}}\otimes\boldsymbol{I}_{P_{2}+1}\right)\boldsymbol{\Pi}_{2}\left(\boldsymbol{\xi}_{1}^{\left(j\right)}\right)^{\mathbf{T}} & i=2
\end{cases}.
\end{equation}

In general, the level and characteristics of the quadrature rule used
in approximating $\hat{\boldsymbol{U}}^{p,d_{i},\tilde{p}_{i}}$ can
differ from the level and characteristics of the quadrature rule used
in constructing the optimally sparse quadrature rules. Moreover, a
strictly positive weight vector with $r_{i}\leq\tilde{N}_{i}+1$ non-zeros
can also be computed using the Nelder-Mead method [\hyperlink{ref31}{31}], in which case,
the upper bound of $\tilde{N}_{i}+1$ on the number of non-zero weights
would coincide with the upper bound proved by Tchakaloff {[}\hyperlink{ref32}{32}{]}.
Furthermore, the sparsity can also be explicitly controlled by prescribing
a threshold on the diagonal elements of $\boldsymbol{R}_{i}^{\ell,2\tilde{p}_{i}}$.

\subsubsection{Selecting the reduced representation}

$\forall i\in\left\{ 1,2\right\} $, we prescribe a tolerance $\epsilon_{i,\mathrm{ord}}>0$,
such that at each iteration $\ell$, the reduced order $\tilde{p}_{i}$
is selected as the smallest $k\in\mathbb{N}$, which satisfies 
\hypertarget{eq355}{}
\begin{equation}
\left\Vert \boldsymbol{U}_{i}^{\ell,p,d_{i},k+1}-\boldsymbol{U}_{i}^{\ell,p,d_{i},k}\right\Vert _{\boldsymbol{G}_{i}}\leq\epsilon_{i,\mathrm{ord}}\left\Vert \boldsymbol{U}_{i}^{\ell,p,d_{i},k+1}\right\Vert _{\boldsymbol{G}_{i}},
\end{equation}
where $\left\Vert \cdot\right\Vert _{\boldsymbol{G}_{i}}$ denotes
the $\boldsymbol{G}_{i}$-weighted Frobenius norm, such that $\forall\hat{\boldsymbol{U}}\in\mathbb{R}^{n\times\left(P+1\right)}$,
$\left\Vert \hat{\boldsymbol{U}}\right\Vert _{\boldsymbol{G}_{i}}=\sqrt{\mathrm{trace}\left(\hat{\boldsymbol{U}}^{\mathbf{T}}\boldsymbol{G}_{i}\hat{\boldsymbol{U}}\right)}$. 

We propose the following heuristic to select $\tilde{p}_{i}$. $\tilde{p}_{i}$
is initialized to $0$ at $\ell=0$, and subsequently incremented
by $1$ at any $\ell>0$, if $\left\Vert \boldsymbol{U}_{i}^{\ell,p,d_{i},\tilde{p}_{i}+1}-\boldsymbol{U}_{i}^{\ell,p,d_{i},\tilde{p}_{i}}\right\Vert _{\boldsymbol{G}_{i}}>\epsilon_{i,\mathrm{ord}}\left\Vert \boldsymbol{U}_{i}^{\ell,p,d_{i},\tilde{p}_{i}+1}\right\Vert _{\boldsymbol{G}_{i}}$.
Therefore, by choosing an appropriate value for the tolerance $\epsilon_{i,\mathrm{ord}}$,
we can guarantee that $\tilde{p}_{i}<p$, and therefore, a reduction
of the gPC approximation order. Moreover, the requirement of computing
$\boldsymbol{U}_{i}^{\ell,p,d_{i},\tilde{p}_{i}+1}$ implies that
the level of the constructed optimal quadrature rule must be $\geq\tilde{p}_{i}+1$. 

\hyperlink{thm2}{Theorem 2} proves an important relationship between the tolerance $\epsilon_{i,\mathrm{ord}}$
and the upper bound on the error incurred by the reduced order approximation
defined in \hyperlink{eq339}{Eq. 3.39}.
\subsubsection*{Theorem 2}
\hypertarget{thm2}{}
\emph{$\forall i\in\left\{ 1,2\right\} $ and iteration $\ell$, let $\boldsymbol{u},\tilde{\boldsymbol{u}}^{\tilde{p}_{i}}:\Theta_{i}^{\ell}\rightarrow\mathbb{R}^{n}$
denote an infinitely regular vector valued function in $\Theta_{i}^{\ell}$
and its reduced gPC approximation of order $\tilde{p}_{i}\geq0$ respectively.
Given a positive-definite Gramian $\boldsymbol{G}\in\mathbb{R}^{n\times n}$,
if $\boldsymbol{u}^{d_{i}}\equiv\boldsymbol{u}^{p,d_{i}}$ and $\boldsymbol{u}^{d_{i},\tilde{p}_{i}}\equiv\boldsymbol{u}^{p,d_{i},\tilde{p}_{i}}$
denote the respective global gPC approximations of the composite functions
$\boldsymbol{u}\circ\boldsymbol{\theta}_{i}^{\ell}$ and $\boldsymbol{u}^{\tilde{p}_{i}}\circ\boldsymbol{\theta}_{i}^{\ell}$,
then $\exists\chi,\chi^{*}>0,\rho^{*}>1$, such that $\forall p,\tilde{p}_{i}\geq0$,}
\hypertarget{eq356}{}
\begin{equation}
\sqrt{\int_{\Xi}\left\Vert \boldsymbol{u}^{d_{i}}\left(\boldsymbol{\xi}\right)-\boldsymbol{u}^{d_{i},\tilde{p}_{i}}\left(\boldsymbol{\xi}\right)\right\Vert _{\boldsymbol{G}}^{2}\mu\left(\boldsymbol{\xi}\right)d\boldsymbol{\xi}}\leq\chi\epsilon_{i,\mathrm{ord}}\sqrt{\int_{\Xi}\left\Vert \boldsymbol{u}^{d_{i}}\left(\boldsymbol{\xi}\right)\right\Vert _{\boldsymbol{G}}^{2}\mu\left(\boldsymbol{\xi}\right)d\boldsymbol{\xi}}+\chi^{*}\rho^{-p}
\end{equation}
\emph{for some $\rho\geq\rho^{*}$.}
\subsubsection*{Proof}
Since $\boldsymbol{u}$ is infinitely regular in $\Theta_{i}^{\ell}$,
$\boldsymbol{u}\circ\boldsymbol{\theta}_{i}^{\ell}$ is infinitely
regular in $\Xi$. Therefore, from the Cameron-Martin theorem, $\exists\chi^{*},\tilde{\chi}>0,\rho^{*}>1$,
such that for any positive-definite $\boldsymbol{G}\in\mathbb{R}^{n\times n}$,
the approximation error between $\boldsymbol{u}^{d_{i}}$ and $\boldsymbol{u}^{d_{i},\tilde{p}_{i}}$
has the following upper bound.
\begin{align}
\hypertarget{eq357}{}
 & \sqrt{\int_{\Xi}\left\Vert \boldsymbol{u}^{d_{i}}\left(\boldsymbol{\xi}\right)-\boldsymbol{u}^{d_{i},\tilde{p}_{i}}\left(\boldsymbol{\xi}\right)\right\Vert _{\boldsymbol{G}}^{2}\mu\left(\boldsymbol{\xi}\right)d\boldsymbol{\xi}}\nonumber \\
 & \leq\tilde{\chi}\sqrt{\int_{\Xi}\left\Vert \boldsymbol{u}^{d_{i},\tilde{p}_{i}+1}\left(\boldsymbol{\xi}\right)-\boldsymbol{u}^{d_{i},\tilde{p}_{i}}\left(\boldsymbol{\xi}\right)\right\Vert _{\boldsymbol{G}}^{2}\mu\left(\boldsymbol{\xi}\right)d\boldsymbol{\xi}}+\chi^{*}\rho^{-p}\nonumber \\
 & =\tilde{\chi}\sqrt{\int_{\Xi}\left\Vert \left(\hat{\boldsymbol{U}}^{d_{i},\tilde{p}_{i}+1}-\hat{\boldsymbol{U}}^{d_{i},\tilde{p}_{i}}\right)\boldsymbol{\psi}\left(\boldsymbol{\xi}\right)\right\Vert _{\boldsymbol{G}}^{2}\mu\left(\boldsymbol{\xi}\right)d\boldsymbol{\xi}}+\chi^{*}\rho^{-p}\nonumber \\
 & =\tilde{\chi}\sqrt{\int_{\Xi}\left\Vert \boldsymbol{G}^{\frac{1}{2}}\left(\hat{\boldsymbol{U}}^{d_{i},\tilde{p}_{i}+1}-\hat{\boldsymbol{U}}^{d_{i},\tilde{p}_{i}}\right)\boldsymbol{\psi}\left(\boldsymbol{\xi}\right)\right\Vert _{2}^{2}\mu\left(\boldsymbol{\xi}\right)d\boldsymbol{\xi}}+\chi^{*}\rho^{-p}\nonumber \\
 & \leq\tilde{\chi}\left\Vert \boldsymbol{G}^{\frac{1}{2}}\left(\hat{\boldsymbol{U}}^{d_{i},\tilde{p}_{i}+1}-\hat{\boldsymbol{U}}^{d_{i},\tilde{p}_{i}}\right)\right\Vert _{2}\sqrt{\int_{\Xi}\left\Vert \boldsymbol{\psi}\left(\boldsymbol{\xi}\right)\right\Vert _{2}^{2}\mu\left(\boldsymbol{\xi}\right)d\boldsymbol{\xi}}+\chi^{*}\rho^{-p}\nonumber \\
 & \leq\chi\left\Vert \left(\hat{\boldsymbol{U}}^{d_{i},\tilde{p}_{i}+1}-\hat{\boldsymbol{U}}^{d_{i},\tilde{p}_{i}}\right)\right\Vert _{\boldsymbol{G}}+\chi^{*}\rho^{-p}
\end{align}
for some $\chi>0,\rho\geq\rho^{*}$. By substituting \hyperlink{eq357}{Eq. 3.57} into
\hyperlink{eq355}{Eq. 3.55}, we arrive at \hyperlink{eq356}{Eq. 3.56}.$\square$

\subsection{Algorithm and computational cost}

The proposed reduced ISP based uncertainty propagation method is described
in \hyperlink{alg2}{Algorithm 2}. Let $\tilde{Q}_{2}$ denote the size of the optimal
quadrature rule corresponding to the reduced dimension $d_{1}$ and
reduced order $\tilde{p}_{1}$, and let $\tilde{Q}_{1}$ denote the
size of the optimal quadrature rule corresponding to the reduced dimension
$d_{2}$ and reduced order $\tilde{p}_{2}$. Therefore, still assuming
that the computational costs are dominated by the repeated execution
of module operators $\tilde{\boldsymbol{M}}_{1}$ and $\tilde{\boldsymbol{M}}_{2}$,
the computational cost of the reduced ISP method
\begin{equation}
\mathcal{C}_{r}\approx\mathcal{O}\left(\bar{\mathcal{C}}_{1}P_{1}^{\alpha_{1}}\tilde{Q}_{2}+\bar{\mathcal{C}}_{2}P_{2}^{\alpha_{2}}\tilde{Q}_{1}\right)
\end{equation}

would grow exponentially with respect to the reduced dimensions $s_{1}+d_{1}$,
$s_{2}+d_{2}$ and reduced orders $\tilde{p}_{1},\tilde{p}_{2}$.
Therefore, the proposed reduced ISP method indeed mitigates the curse
of dimensionality associated with the standard ISP method. 

\begin{algorithm}
\hypertarget{alg51}{}
\caption{Reduced ISP based uncertainty propagation for a two-module multi-physics system}
\SetKwInOut{Input}{inputs}\SetKwInOut{Output}{outputs}
\SetKwRepeat{Repeat}{repeat}{}
\DontPrintSemicolon
\Input{order $p\geq0$, level $q\geq p$, $\epsilon_{1,\mathrm{dim}}$,
$\epsilon_{1,\mathrm{dim}}$, $\epsilon_{2,\mathrm{ord}}$, $\epsilon_{2,\mathrm{ord}}$,
$\hat{\boldsymbol{U}}_{1}^{0}$, $\hat{\boldsymbol{U}}_{2}^{0}$} 
\Output{$\hat{\boldsymbol{U}}_{1}$, $\hat{\boldsymbol{U}}_{2}$}
\textbf{precompute}: $\left\{ \left(w_{1}^{\left(j\right)},\boldsymbol{\Pi}_{2}\left(\boldsymbol{\xi}_{1}^{\left(j\right)}\right)\right)\right\} _{j=1}^{Q_{1}}$,
$\left\{ \left(w_{2}^{\left(j\right)},\boldsymbol{\Pi}_{1}\left(\boldsymbol{\xi}_{2}^{\left(j\right)}\right)\right)\right\} _{j=1}^{Q_{2}}$

$\ell\leftarrow0$, $\tilde{p}_{1}\leftarrow0$, $\tilde{p}_{2}\leftarrow0$

\Repeat{(contd.)}{

\SetKwBlock{Begin}{dimension reduction}{end}
\Begin{
\Input{$\hat{\boldsymbol{U}}_{1}^{\ell}$, $\hat{\boldsymbol{U}}_{2}^{\ell}$,
$\epsilon_{1,\mathrm{dim}}$}

\Output{$\left\{ \left(\tilde{\boldsymbol{U}}_{1,1,j},\tilde{\boldsymbol{U}}_{2,1,j}\right)\right\} _{j=0}^{d_{1}}$,
$\left\{ \boldsymbol{\theta}_{1}^{\left(j\right)}=\left[\begin{array}{ccc}
\theta_{1,1}^{\left(j\right)} & \cdots & \theta_{1,d_{1}}^{\left(j\right)}\end{array}\right]^{\mathbf{T}}\right\} _{j=1}^{Q_{2}}$}
}
\SetKwBlock{Begin}{reduced basis/quadrature construction}{end}
\Begin{
\Input{$\left\{ \left(\boldsymbol{\theta}_{1}^{\left(j\right)},w_{2}^{\left(j\right)}\right)\right\} _{j=1}^{Q_{2}}$,
$\tilde{p}_{1}$}

\Output{$\left\{ \left(\boldsymbol{\phi}_{1}^{\ell,\tilde{p}_{1}}\left(\boldsymbol{\theta}_{1}^{\left(j\right)}\right),\boldsymbol{\phi}_{1}^{\ell,\tilde{p}_{1}+1}\left(\boldsymbol{\theta}_{1}^{\left(j\right)}\right),\tilde{w}_{2}^{\left(j\right)}\right)\right\} _{j=1}^{Q_{2}}$,
$\boldsymbol{S}_{1}^{\ell,\tilde{p}_{1}}$, $\boldsymbol{S}_{1}^{\ell,\tilde{p}_{1}+1}$,
$\mathcal{Z}_{2}$}
}

$\hat{\boldsymbol{U}}_{1}^{d_{1},\tilde{p}_{1}}\leftarrow\boldsymbol{0}$,
$\hat{\boldsymbol{U}}_{1}^{d_{1},\tilde{p}_{1}+1}\leftarrow\boldsymbol{0}$

\For{$ j\in\tilde{\mathcal{Z}}_{2}$}{

$\tilde{\boldsymbol{U}}_{1,1}\leftarrow\tilde{\boldsymbol{U}}_{1,1,0}$,
$\tilde{\boldsymbol{U}}_{2,1}\leftarrow\tilde{\boldsymbol{U}}_{2,1,0}$

\For{$ k\leftarrow 1$ \KwTo $d_{1}$}{

$\tilde{\boldsymbol{U}}_{1,1}\leftarrow\tilde{\boldsymbol{U}}_{1,1}+\theta_{1,k}^{\left(j\right)}\tilde{\boldsymbol{U}}_{1,1,k}$,
$\tilde{\boldsymbol{U}}_{2,1}\leftarrow\tilde{\boldsymbol{U}}_{2,1}+\theta_{1,k}^{\left(j\right)}\tilde{\boldsymbol{U}}_{2,1,k}$

}
$\tilde{\boldsymbol{U}}_{1}\leftarrow\tilde{\boldsymbol{M}}_{1}\left(\tilde{\boldsymbol{U}}_{1,1},\tilde{\boldsymbol{U}}_{2,1}\right)$

$\hat{\boldsymbol{U}}_{1}^{d_{1},\tilde{p}_{1}}\leftarrow\hat{\boldsymbol{U}}_{1}^{d_{1},\tilde{p}_{1}}+\tilde{w}_{2}^{\left(j\right)}\tilde{\boldsymbol{U}}_{1}\left(\left(\boldsymbol{\phi}_{1}^{\ell,\tilde{p}_{1}}\left(\boldsymbol{\theta}_{1}^{\left(j\right)}\right)^{\mathbf{T}}\boldsymbol{S}_{1}^{\ell,\tilde{p}_{1}}\right)\otimes\boldsymbol{I}_{P_{1}+1}\right)$

$\hat{\boldsymbol{U}}_{1}^{d_{1},\tilde{p}_{1}+1}\leftarrow\hat{\boldsymbol{U}}_{1}^{d_{1},\tilde{p}_{1}+1}+\tilde{w}_{2}^{\left(j\right)}\tilde{\boldsymbol{U}}_{1}\left(\left(\boldsymbol{\phi}_{1}^{\ell,\tilde{p}_{1}+1}\left(\boldsymbol{\theta}_{1}^{\left(j\right)}\right)^{\mathbf{T}}\boldsymbol{S}_{1}^{\ell,\tilde{p}_{1}+1}\right)\otimes\boldsymbol{I}_{P_{1}+1}\right)$

}

$\hat{\boldsymbol{U}}_{1}^{\tilde{p}_{1}}\leftarrow\boldsymbol{0}$,
$\hat{\boldsymbol{U}}_{1}^{\tilde{p}_{1}+1}\leftarrow\boldsymbol{0}$

$\hat{\boldsymbol{U}}_{1}^{\tilde{p}_{1}}\leftarrow\boldsymbol{0}$,
$\hat{\boldsymbol{U}}_{1}^{\tilde{p}_{1}+1}\leftarrow\boldsymbol{0}$

\For{$ j\leftarrow 1$ \KwTo $Q_{2}$}{

$\hat{\boldsymbol{U}}_{1}^{\tilde{p}_{1}}\leftarrow\hat{\boldsymbol{U}}_{1}^{\tilde{p}_{1}}+w_{2}^{\left(j\right)}\hat{\boldsymbol{U}}_{1}^{d_{1},\tilde{p}_{1}}\left(\boldsymbol{\phi}_{1}^{\ell,\tilde{p}_{1}}\left(\boldsymbol{\theta}_{1}^{\left(j\right)}\right)\otimes\boldsymbol{I}_{P_{1}+1}\right)\boldsymbol{\Pi}_{1}\left(\boldsymbol{\xi}_{2}^{\left(j\right)}\right)^{\mathbf{T}}$

$\hat{\boldsymbol{U}}_{1}^{\tilde{p}_{1}+1}\leftarrow\hat{\boldsymbol{U}}_{1}^{\tilde{p}_{1}+1}+w_{2}^{\left(j\right)}\hat{\boldsymbol{U}}_{1}^{d_{1},\tilde{p}_{1}+1}\left(\boldsymbol{\phi}_{1}^{\ell,\tilde{p}_{1}+1}\left(\boldsymbol{\theta}_{1}^{\left(j\right)}\right)\otimes\boldsymbol{I}_{P_{1}+1}\right)\boldsymbol{\Pi}_{1}\left(\boldsymbol{\xi}_{2}^{\left(j\right)}\right)^{\mathbf{T}}$

}
\If{$\left\Vert \hat{\boldsymbol{U}}_{1}^{\tilde{p}_{1}+1}-\hat{\boldsymbol{U}}_{1}^{\tilde{p}_{1}}\right\Vert _{\boldsymbol{\boldsymbol{G}}_{1}}>\epsilon_{1,\mathrm{ord}}\left\Vert \hat{\boldsymbol{U}}_{1}^{\tilde{p}_{1}+1}\right\Vert _{\boldsymbol{\boldsymbol{G}}_{1}}$}
{$\tilde{p}_{1}\leftarrow\tilde{p}_{1}+1$}

$\hat{\boldsymbol{U}}_{1}^{\ell+1}\leftarrow\hat{\boldsymbol{U}}_{1}^{\tilde{p}_{1}}$

}
\end{algorithm}
\RestyleAlgo{boxed}
\begin{algorithm}
\SetKwInOut{Input}{inputs}\SetKwInOut{Output}{outputs}
\SetKwRepeat{Repeat}{}{until}
\DontPrintSemicolon
 \emph{(contd.)}\Repeat{$\hat{\boldsymbol{U}}_{1}^{\ell}$, $\hat{\boldsymbol{U}}_{2}^{\ell}$ $\mathrm{converge}$}{

\SetKwBlock{Begin}{dimension reduction}{end}
\Begin{

\Input{$\hat{\boldsymbol{U}}_{2}^{\ell}$, $\hat{\boldsymbol{U}}_{1}^{\ell+1}$,
$\epsilon_{2,\mathrm{dim}}$}

\Output{$\left\{ \left(\tilde{\boldsymbol{U}}_{2,2,j},\tilde{\boldsymbol{U}}_{1,2,j}\right)\right\} _{j=0}^{d_{2}}$,
$\left\{ \boldsymbol{\theta}_{2}^{\left(j\right)}=\left[\begin{array}{ccc}
\theta_{2,1}^{\left(j\right)} & \cdots & \theta_{2,d_{1}}^{\left(j\right)}\end{array}\right]^{\mathbf{T}}\right\} _{j=1}^{Q_{1}}$}
}

\SetKwBlock{Begin}{reduced basis/quadrature construction}{end}
\Begin{

\Input{$\left\{ \left(\boldsymbol{\theta}_{2}^{\left(j\right)},w_{1}^{\left(j\right)}\right)\right\} _{j=1}^{Q_{2}}$,
$\tilde{p}_{2}$
}

\Output{$\left\{ \left(\boldsymbol{\phi}_{2}^{\ell,\tilde{p}_{2}}\left(\boldsymbol{\theta}_{2}^{\left(j\right)}\right),\boldsymbol{\phi}_{2}^{\ell,\tilde{p}_{2}+1}\left(\boldsymbol{\theta}_{2}^{\left(j\right)}\right),\tilde{w}_{1}^{\left(j\right)}\right)\right\} _{j=1}^{Q_{1}}$,
$\boldsymbol{S}_{2}^{\ell,\tilde{p}_{2}}$, $\boldsymbol{S}_{2}^{\ell,\tilde{p}_{2}+1}$,
$\mathcal{Z}_{1}$}
}
$\hat{\boldsymbol{U}}_{2}^{d_{2},\tilde{p}_{2}}\leftarrow\boldsymbol{0}$,
$\hat{\boldsymbol{U}}_{2}^{d_{2},\tilde{p}_{2}+1}\leftarrow\boldsymbol{0}$

\For{$ j\in\tilde{\mathcal{Z}}_{1}$}{

$\tilde{\boldsymbol{U}}_{2,2}\leftarrow\tilde{\boldsymbol{U}}_{2,2,0}$,
$\tilde{\boldsymbol{U}}_{1,2}\leftarrow\tilde{\boldsymbol{U}}_{1,2,0}$

\For {$ k\leftarrow 1$ \KwTo $d_{2}$}{

$\tilde{\boldsymbol{U}}_{2,2}\leftarrow\tilde{\boldsymbol{U}}_{2,2}+\theta_{2,k}^{\left(j\right)}\tilde{\boldsymbol{U}}_{2,2,k}$,
$\tilde{\boldsymbol{U}}_{1,2}\leftarrow\tilde{\boldsymbol{U}}_{1,2}+\theta_{2,k}^{\left(j\right)}\tilde{\boldsymbol{U}}_{1,2,k}$

}
$\tilde{\boldsymbol{U}}_{2}\leftarrow\tilde{\boldsymbol{M}}_{2}\left(\tilde{\boldsymbol{U}}_{2,2},\tilde{\boldsymbol{U}}_{1,2}\right)$

$\hat{\boldsymbol{U}}_{2}^{d_{2},\tilde{p}_{2}}\leftarrow\hat{\boldsymbol{U}}_{2}^{d_{2},\tilde{p}_{2}}+\tilde{w}_{1}^{\left(j\right)}\tilde{\boldsymbol{U}}_{2}\left(\left(\boldsymbol{\phi}_{2}^{\ell,\tilde{p}_{2}}\left(\boldsymbol{\theta}_{2}^{\left(j\right)}\right)^{\mathbf{T}}\boldsymbol{S}_{2}^{\ell,\tilde{p}_{2}}\right)\otimes\boldsymbol{I}_{P_{2}+1}\right)$

$\hat{\boldsymbol{U}}_{2}^{d_{2},\tilde{p}_{2}+1}\leftarrow\hat{\boldsymbol{U}}_{2}^{d_{2},\tilde{p}_{2}+1}+\tilde{w}_{1}^{\left(j\right)}\tilde{\boldsymbol{U}}_{2}\left(\left(\boldsymbol{\phi}_{2}^{\ell,\tilde{p}_{2}+1}\left(\boldsymbol{\theta}_{2}^{\left(j\right)}\right)^{\mathbf{T}}\boldsymbol{S}_{2}^{\ell,\tilde{p}_{1}+1}\right)\otimes\boldsymbol{I}_{P_{2}+1}\right)$

}
\For {$ j \leftarrow 1$ \KwTo $Q_{1}$}{

$\hat{\boldsymbol{U}}_{2}^{\tilde{p}_{2}}\leftarrow\hat{\boldsymbol{U}}_{2}^{\tilde{p}_{2}}+w_{1}^{\left(j\right)}\hat{\boldsymbol{U}}_{2}^{d_{2},\tilde{p}_{2}}\left(\boldsymbol{\phi}_{2}^{\ell,\tilde{p}_{2}}\left(\boldsymbol{\theta}_{2}^{\left(j\right)}\right)\otimes\boldsymbol{I}_{P_{2}+1}\right)\boldsymbol{\Pi}_{2}\left(\boldsymbol{\xi}_{1}^{\left(j\right)}\right)^{\mathbf{T}}$

$\hat{\boldsymbol{U}}_{2}^{\tilde{p}_{2}+1}\leftarrow\hat{\boldsymbol{U}}_{2}^{\tilde{p}_{2}+1}+w_{1}^{\left(j\right)}\hat{\boldsymbol{U}}_{2}^{d_{2},\tilde{p}_{2}+1}\left(\boldsymbol{\phi}_{2}^{\ell,\tilde{p}_{2}+1}\left(\boldsymbol{\theta}_{2}^{\left(j\right)}\right)\otimes\boldsymbol{I}_{P_{2}+1}\right)\boldsymbol{\Pi}_{2}\left(\boldsymbol{\xi}_{1}^{\left(j\right)}\right)^{\mathbf{T}}$

}

\If{$\left\Vert \hat{\boldsymbol{U}}_{2}^{\tilde{p}_{2}+1}-\hat{\boldsymbol{U}}_{2}^{\tilde{p}_{2}}\right\Vert _{\boldsymbol{\boldsymbol{G}}_{2}}>\epsilon_{2,\mathrm{ord}}\left\Vert \hat{\boldsymbol{U}}_{2}^{\tilde{p}_{2}+1}\right\Vert _{\boldsymbol{\boldsymbol{G}}_{2}}$}
{$\tilde{p}_{2}\leftarrow\tilde{p}_{2}+1$}

$\hat{\boldsymbol{U}}_{2}^{\ell+1}\leftarrow\hat{\boldsymbol{U}}_{2}^{\tilde{p}_{2}}$

$\ell\leftarrow\ell+1$

}
\end{algorithm}
\subsection{Error analysis}

$\forall i\in\left\{ 1,2\right\} $, let $\varepsilon_{i}:\forall\boldsymbol{\xi}\in\Xi$,
\begin{equation}
\varepsilon_{i}\left(\boldsymbol{\xi}\right)=\left\Vert \boldsymbol{u}_{i}\left(\boldsymbol{\xi}\right)-\boldsymbol{u}_{i}^{\ell,p,d_{i},\tilde{p}_{i}}\left(\boldsymbol{\xi}\right)\right\Vert _{\boldsymbol{G}_{i}}
\end{equation}
denote the mean-square error between the component solution $\boldsymbol{u}_{i}$
and its corresponding reduced gPC approximation $\boldsymbol{u}_{i}^{\ell,p,d_{i},\tilde{p}_{i}}\left(\boldsymbol{\xi}\right)$. 

Subsequently, $\varepsilon_{i}$ can be decomposed as a sum of individual
error terms as follows.
\begin{equation}
\varepsilon_{i}=\varepsilon_{i,\mathrm{BGS}}+\varepsilon_{i,\mathrm{gPC}}+\varepsilon_{i,\mathrm{dim}}+\varepsilon_{i,\mathrm{ord}},
\end{equation}
where $\varepsilon_{i,\mathrm{BGS}}$ denotes the convergence error,
$\varepsilon_{i,\mathrm{gPC}}$ denotes the gPC truncation error,
$\varepsilon_{i,\mathrm{dim}}$ denotes the reduced dimension approximation
error and $\varepsilon_{i,\mathrm{ord}}$ denotes the reduced order
approximation error. 

Using the triangle inequality property of norms,
an asymptotic upper bound for each constituent error term can be formulated
as follows. $\forall\boldsymbol{\xi}\in\Xi$,
\hypertarget{eq361}{}
\begin{align}
\varepsilon_{i}\left(\boldsymbol{\xi}\right)= & \left\Vert \boldsymbol{u}_{i}\left(\boldsymbol{\xi}\right)-\boldsymbol{u}_{i}^{\ell}\left(\boldsymbol{\xi}\right)+\boldsymbol{u}_{i}^{\ell}\left(\boldsymbol{\xi}\right)-\boldsymbol{u}_{i}^{\ell,p}\left(\boldsymbol{\xi}\right)+\boldsymbol{u}_{i}^{\ell,p}\left(\boldsymbol{\xi}\right)-\boldsymbol{u}_{i}^{\ell,p,d_{i}}\left(\boldsymbol{\xi}\right)+\boldsymbol{u}_{i}^{\ell,p,d_{i}}\left(\boldsymbol{\xi}\right)\right.\nonumber \\
 & \left.-\boldsymbol{u}_{i}^{\ell,p,d_{i},\tilde{p}_{i}}\left(\boldsymbol{\xi}\right)\right\Vert _{\boldsymbol{G}_{i}}\nonumber \\
\leq & \underset{\varepsilon_{i,\mathrm{BGS}}\left(\boldsymbol{\xi}\right)\leq\mathcal{O}\left(\eta^{-\ell}\right)}{\underbrace{\left\Vert \boldsymbol{u}_{i}\left(\boldsymbol{\xi}\right)-\boldsymbol{u}_{i}^{\ell}\left(\boldsymbol{\xi}\right)\right\Vert _{\boldsymbol{G}_{i}}}}+\underset{\varepsilon_{i,\mathrm{gPC}}\left(\boldsymbol{\xi}\right)\leq\mathcal{O}\left(\rho^{-p}\right)}{\underbrace{\left\Vert \boldsymbol{u}_{i}^{\ell}\left(\boldsymbol{\xi}\right)-\boldsymbol{u}_{i}^{\ell,p}\left(\boldsymbol{\xi}\right)\right\Vert _{\boldsymbol{G}_{i}}}}+\underset{\varepsilon_{i,\mathrm{dim}}\left(\boldsymbol{\xi}\right)\leq\mathcal{O}\left(\epsilon_{i,\mathrm{dim}}\right)}{\underbrace{\left\Vert \boldsymbol{u}_{i}^{\ell,p}\left(\boldsymbol{\xi}\right)-\boldsymbol{u}_{i}^{\ell,p,d_{i}}\left(\boldsymbol{\xi}\right)\right\Vert _{\boldsymbol{G}_{i}}}}\nonumber \\
 & +\underset{\varepsilon_{i,\mathrm{ord}}\left(\boldsymbol{\xi}\right)\leq\mathcal{O}\left(\epsilon_{i,\mathrm{ord}}\right)+\mathcal{O}\left(\tilde{\rho}^{-p}\right)}{\underbrace{\left\Vert \boldsymbol{u}_{i}^{\ell,p,d_{i}}\left(\boldsymbol{\xi}\right)-\boldsymbol{u}_{i}^{\ell,p,d_{i},\tilde{p}_{i}}\left(\boldsymbol{\xi}\right)\right\Vert _{\boldsymbol{G}_{i}}}}.
\end{align}

In the standard ISP method, the asymptotic upper bound on the approximationta
error would simply be the sum of the first two terms on the right
hand side of \hyperlink{eq361}{Eq. 3.61}, implying that $\boldsymbol{u}_{i}^{\ell,p}$
would converge to $\boldsymbol{u}_{i}$ as $\ell,p\rightarrow\infty$.
However, in the reduced ISP method, $\varepsilon_{i}$ would converge
to a non-zero quantity, and have an asymptotic upper bound of $\mathcal{O}\left(\varepsilon_{i,\mathrm{dim}}\right)+\mathcal{O}\left(\varepsilon_{i,\mathrm{ord}}\right)$. 

This analysis suggests that the approximation error in $\boldsymbol{u}_{i}^{\ell,p,d_{i},\tilde{p}_{i}}$
can be explicitly controlled by choosing the tolerances $\varepsilon_{i,\mathrm{dim}}$
and $\varepsilon_{i,\mathrm{ord}}$ appropriately.

\subsection{Selecting the tolerance values}
In practice, since the exact bounds on the approximation error are not known apriori, the tolerance values are selected based on the results computed in a preliminary study, wherein a lower fidelity multi-physics model is used and the various tolerance values can be tested. In our case, lower fidelity translates to a coarser spatial discretization of the original steady state coupled PDE system to formulate the coupled algebraic system. An illustration of this is provided in the numerical experiments in \hyperlink{sec4}{\S4}

\section{Numerical examples}

We will now demonstrate and compare the performance of the standard
and reduced ISP methods using two numerical examples.

\subsection{Thermal-Neutronics problem}

The first numerical example we consider is inspired from recent benchmark
studies using nuclear reactor models {[}\hyperlink{ref14}{14}, \hyperlink{ref15}{15}{]}. Here, we study the
steady-state neutron transport and heat transfer in a two-dimensional
reactor {[}\hyperlink{ref33}{33}{]}, with uncertain thermal and reactive properties.

\subsubsection{Model setup}

\begin{figure}
    \hypertarget{fig1}{}
    \centering
    \includegraphics[bb=125bp 500bp 475bp 720bp,clip,scale=0.65]{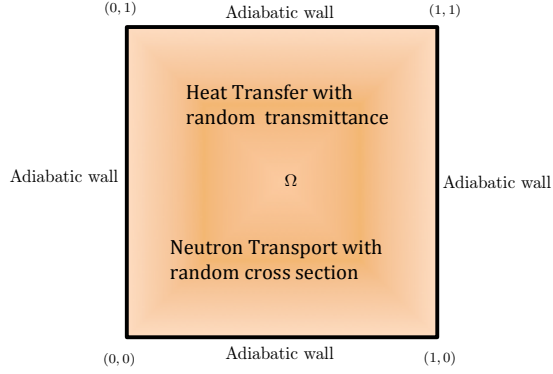}  
    \caption{Computational domain for the Thermal-Neutronics problem.}
    \label{fig:fig1}
\end{figure}

Let $\Omega\equiv\left(0,1\right)_{x_{1}}\times\left(0,1\right)_{x_{2}}$
denote the spatial domain and $T,\varphi$ denote the respective non-dimensional
temperature, neutron flux in a reactor. The governing equations for
the fluid variables are as follows. $\forall\boldsymbol{\xi}\in\Xi$,
\begin{align}
\boldsymbol{\nabla}^{\mathbf{T}}\boldsymbol{\nabla}T\left(\boldsymbol{x},\boldsymbol{\xi}\right)-H_{T}\left(\boldsymbol{x},\boldsymbol{\xi}_{1}\left(\boldsymbol{\xi}\right)\right)T\left(\boldsymbol{x},\boldsymbol{\xi}\right)+\frac{E_{T}\varphi\left(\boldsymbol{x},\boldsymbol{\xi}\right)}{\sqrt{T\left(\boldsymbol{x},\boldsymbol{\xi}\right)+1}} & =0,\nonumber \\
\boldsymbol{\nabla}^{\mathbf{T}}\left(\sqrt{T\left(\boldsymbol{x},\boldsymbol{\xi}\right)+1}\boldsymbol{\nabla}\varphi\left(\boldsymbol{x},\boldsymbol{\xi}\right)\right)-\frac{\varSigma_{\varphi}\left(\boldsymbol{x},\boldsymbol{\xi}_{2}\left(\boldsymbol{\xi}\right)\right)\varphi\left(\boldsymbol{x},\boldsymbol{\xi}\right)}{\sqrt{T\left(\boldsymbol{x},\boldsymbol{\xi}\right)+1}}+S_{\varphi} & =0, & \boldsymbol{x}\in\Omega,
\end{align}
with homogenous Neumann boundary conditions for $T$ and $\varphi$
at all boundaries. \hyperlink{fig1}{Figure 1} illustrates the multi-physics setup and
computational domain. 

Here, $H_{T}$, $E_{T}$, $\varSigma_{\varphi}$ and $S_{\varphi}$
denote the (non-dimensional) thermal transmittance, fission energy,
reaction cross section and neutron source strength respectively. In
this study, $H_{T}$ and $\varSigma_{\varphi}$ are assumed to be
independent random fields, and modeled using the following KL expansions.
$\forall\boldsymbol{x}\in\Omega,\boldsymbol{\xi}_{1}\in\Xi_{1}$,
\begin{equation}
H_{T}\left(\boldsymbol{x},\boldsymbol{\xi}_{1}\right)=\bar{H}_{T}+\sqrt{3}\delta_{H_{T}}\sum_{j=1}^{s_{1}}\gamma_{H_{T},j}\left(\boldsymbol{x}\right)\xi_{1j},
\end{equation}
where $\bar{H}_{T}$ denotes the mean of $H_{T}$ and $\left\{ \xi_{1j}\sim U\left[-1.1\right]\right\} _{j=1}^{s_{1}}$
are $i.i.d.$ random variables. Similarly, $\forall\boldsymbol{x}\in\Omega,\boldsymbol{\xi}_{2}\in\Xi_{2}$,
\begin{equation}
\varSigma_{\varphi}\left(\boldsymbol{x},\boldsymbol{\xi}_{2}\right)=\bar{\varSigma}_{\varphi}+\sqrt{3}\delta_{\varSigma_{\varphi}}\sum_{j=1}^{s_{2}}\gamma_{\varSigma_{\varphi},j}\left(\boldsymbol{x}\right)\xi_{2j},
\end{equation}
where $\bar{\varSigma}_{\varphi}$ denotes the mean of $\varSigma_{\varphi}$
and $\left\{ \xi_{2j}\sim U\left[-1.1\right]\right\} _{j=1}^{s_{2}}$
are $i.i.d.$ random variables. Moreover, we assume that both $H_{T}$
and $\varSigma_{\varphi}$ have exponential covariance kernels
\begin{align}
C_{H_{T}}\left(\boldsymbol{x},\boldsymbol{y}\right) & =\delta_{H_{T}}^{2}\exp\left(-\frac{\left\Vert \boldsymbol{x}-\boldsymbol{y}\right\Vert _{1}}{l_{H_{T}}}\right),\nonumber \\
C_{\varSigma_{\varphi}}\left(\boldsymbol{x},\boldsymbol{y}\right) & =\delta_{\varSigma_{\varphi}}^{2}\exp\left(-\frac{\left\Vert \boldsymbol{x}-\boldsymbol{y}\right\Vert _{1}}{l_{\varSigma_{\varphi}}}\right), & \boldsymbol{x},\boldsymbol{y}\in\Omega
\end{align}
where $\delta_{H_{T}},\delta_{\varSigma_{\varphi}}$ denote the respective
coefficients of variation, and $l_{H_{T}},l_{\varSigma_{\varphi}}$
denote the respective correlation lengths. The analytical expressions
for $\left\{ \gamma_{H_{T},j}\right\} _{j>0}$ and $\left\{ \gamma_{\varSigma_{\varphi},j}\right\} _{j>0}$
are provided in Appendix A. 

The governing PDE system is spatially discretized using bilinear finite
elements {[}\hyperlink{ref34}{34}{]} and $m\times m$ equispaced nodes in $\Omega$. Let
$\boldsymbol{u}_{1},\boldsymbol{u}_{2}\in\mathbb{R}^{m^{2}}$ denote
the respective vectors of nodal temperature and neutron flux values,
which solve the nonlinear system
\hypertarget{eq45}{}
\begin{align}
\left(\boldsymbol{K}_{T}+\boldsymbol{H}_{T}\left(\boldsymbol{\xi}_{1}\right)\right)\boldsymbol{u}_{1}-\boldsymbol{E}_{T}\left(\boldsymbol{u}_{1}\right)\boldsymbol{u}_{2} & =\boldsymbol{0},\nonumber \\
\left(\boldsymbol{K}_{\varphi}\left(\boldsymbol{u}_{1}\right)+\boldsymbol{\varSigma}_{\varphi}\left(\boldsymbol{u}_{1},\boldsymbol{\xi}_{2}\right)\right)\boldsymbol{u}_{2}-\boldsymbol{s}_{\varphi} & =\boldsymbol{0},
\end{align}
where each term in \hyperlink{eq45}{Eq. 4.5} denotes its respective discretized operator
in the coupled PDE system. Subsequently, a multi-physics setup is
formulated by separating the thermal and neutronic components of the
coupled algebraic system, wherein the component residuals, as per
\hyperlink{eq21}{Eq. 2.1}, are formulated as follows. 
\hypertarget{eq46}{}
\begin{align}
\boldsymbol{f}_{1}\left(\boldsymbol{u}_{1};\boldsymbol{u}_{2},\boldsymbol{\xi}_{1}\right) & =\left(\boldsymbol{K}_{T}+\boldsymbol{H}_{T}\left(\boldsymbol{\xi}_{1}\right)\right)\boldsymbol{u}_{1}-\boldsymbol{E}_{T}\left(\boldsymbol{u}_{1}\right)\boldsymbol{u}_{2},\nonumber \\
\boldsymbol{f}_{2}\left(\boldsymbol{u}_{2};\boldsymbol{u}_{1},\boldsymbol{\xi}_{2}\right) & =\left(\boldsymbol{K}_{\varphi}\left(\boldsymbol{u}_{1}\right)+\boldsymbol{\varSigma}_{\varphi}\left(\boldsymbol{u}_{1},\boldsymbol{\xi}_{2}\right)\right)\boldsymbol{u}_{2}-\boldsymbol{s}_{\varphi}.
\end{align}

The quantities of interest in this study are the first two moments
and sensitivity indices of the temperature and neutron flux fields.
\hyperlink{tab1}{Table 1} lists the numerical values of the deterministic parameters
used in this study.

\begin{table}[htbp]
\hypertarget{tab1}{}
\caption{Deterministic parameter values in the Thermal-Neutronics
problem.}
\begin{center}\scriptsize
\renewcommand{\arraystretch}{1.3}
\begin{tabular}{cccccccc}
\toprule 
$\bar{H}_{T}$ & $E_{T}$ & $\bar{\varSigma}_{\varphi}$ & $S_{\varphi}$ & $\delta_{H_{T}}$ & $\delta_{\varSigma_{\varphi}}$ & $l_{H_{T}}$ & $l_{\varSigma_{\varphi}}$\tabularnewline
\midrule
\midrule 
17 & 5.8 & 13.6 & 22.7 & 9.0 & 9.0 & 0.15 & 0.5\tabularnewline
\bottomrule
\end{tabular}
\par\end{center}
\end{table}

\subsubsection{Modular deterministic solver: Setup and verification}

To compare the results of the gPC solution approximations with Monte-Carlo
solution samples, we implemented a modular deterministic solver for
the coupled algebraic system, where a BGS method is used to converge
to the solution. The corresponding iterations are computed as follows.
\begin{align}
\boldsymbol{u}_{1}^{\ell+1} & =\left(\boldsymbol{K}_{T}+\boldsymbol{H}_{T}\left(\boldsymbol{\xi}_{1}\right)\right)^{-1}\boldsymbol{E}_{T}\left(\boldsymbol{u}_{1}^{\ell}\right)\boldsymbol{u}_{2}^{\ell},\nonumber \\
\boldsymbol{u}_{2}^{\ell+1} & =\left(\boldsymbol{K}_{\varphi}\left(\boldsymbol{u}_{1}^{\ell+1}\right)+\boldsymbol{\varSigma}_{\varphi}\left(\boldsymbol{u}_{1}^{\ell+1},\boldsymbol{\xi}_{2}\right)\right)^{-1}\boldsymbol{s}_{\varphi}.
\end{align}

Subsequently, a verification study was carried out on the modular
deterministic solver, using the method of manufactured solutions (MMS)
{[}\hyperlink{ref35}{35}{]}. Further details are provided in Appendix B. Moreover, similar
formulations define the BGS iterations in the standard and reduced
ISP method implementations, which respectively solve the global and
modular Galerkin forms of the coupled algebraic system formulated
in \hyperlink{eq46}{Eq. 4.6}. In each implementation, the the components of the solver
were developed as $\mathtt{MATLAB}^{\text{\texttrademark}}$ function
modules,. Furthermore, the corresponding linear systems are solved
using LSQR {[}\hyperlink{ref36}{36}{]}, with the mean of the stochastic left hand side
matrices used in defining preconditioners {[}\hyperlink{ref37}{37}{]}.

\subsubsection{ISP based uncertainty propagation}

Here, we compare the results and performance of the standard and reduced
ISP method implementations. For $m=21$, $s_{1}=s_{2}=3$ and $p=4$, \hyperlink{fig2}{Figure 2} and \hyperlink{fig3}{Figure 3} respectively compare the mean and standard deviation
of the temperature and neutron flux fields computed from their converged
gPC coefficients. Moreover, \hyperlink{tab2}{Table 2} compares the sensitivity indices
computed using the ANOVA decomposition technique [\hyperlink{ref38}{38}]. Here, $\forall i\in\left\{1,2\right\}$, $\hat{\mathcal{V}}_{i,1}$ and $\hat{\mathcal{V}}_{i,2}$ denote the main effects of $\boldsymbol{\xi_{1}}$ and $\boldsymbol{\xi_{2}}$, while $\hat{\mathcal{V}}_{i,12}$ denotes the interaction effects on the variance of $\boldsymbol{u}_{i}$. Therefore, $\hat{\mathcal{V}}_{i,1}+\hat{\mathcal{V}}_{i,2}+\hat{\mathcal{V}}_{i,12}=1:$

\begin{equation}
\hat{\mathcal{V}}_{i,1}=\frac{\displaystyle\sum_{j=1}^{P_1} \left\Vert  \hat{\boldsymbol{u}}_{i,j,0} \right\Vert_{\boldsymbol{G}}^{2}}{\displaystyle\sum_{j=1}^{P} \left\Vert  \hat{\boldsymbol{u}}_{i,\jmath_{1}\left(j\right),\jmath_{2}\left(j\right)} \right\Vert_{\boldsymbol{G}}^{2}},\ \hat{\mathcal{V}}_{i,2}=\frac{\displaystyle\sum_{j=1}^{P_2} \left\Vert  \hat{\boldsymbol{u}}_{i,0,j} \right\Vert_{\boldsymbol{G}}^{2}}{\displaystyle\sum_{j=1}^{P} \left\Vert  \hat{\boldsymbol{u}}_{i,\jmath_{1}\left(j\right),\jmath_{2}\left(j\right)} \right\Vert_{\boldsymbol{G}}^{2}}.
\end{equation}

Furthermore, for various values of
$s_{1}$, $s_{2}$ and $p$, a comparison of the respective errors
$\varepsilon_{s},\varepsilon_{r}$ and overall computational costs
$\mathcal{C}_{s}$, $\mathcal{C}_{r}$ (wall-times) is provided in
\hyperlink{tab3}{Table 3}. The tolerance values used in the reduced ISP method implementation
are $\epsilon_{1,\mathrm{dim}}=0.02$, $\epsilon_{2,\mathrm{dim}}=0.05$,
and $\epsilon_{1,\mathrm{ord}}=\epsilon_{2,\mathrm{ord}}=10^{-4}$.
While the choice of these tolerances seems arbitrary, they were selected
based on preliminary experiments using a coarser mesh, with $m=11$,
$s_{1}=s_{2}=3$ and $p=4$, to study the effects of changing $\epsilon_{1,\mathrm{dim}}$
and $\epsilon_{2,\mathrm{dim}}$ on the computed gPC coefficients
of the solutions. \hyperlink{fig4}{Figure 4} illustrates a comparison of the standard
deviations of the random neutron flux.

In each case, both the standard and reduced ISP algorithms
invariantly converged in $6$ iterations for the chosen tolerance $\epsilon_{\mathrm{BGS}}=10^{-6}$. This observation indicates that the convergence rate the BGS
iterations in both implementations were independent of the stochastic dimensions $s_{1},s_{2}$ and
the gPC order $p$.

\begin{figure}
    \hypertarget{fig2}{}
    \centering
    \includegraphics[bb=80bp 240bp 1100bp 905bp,scale=0.4]{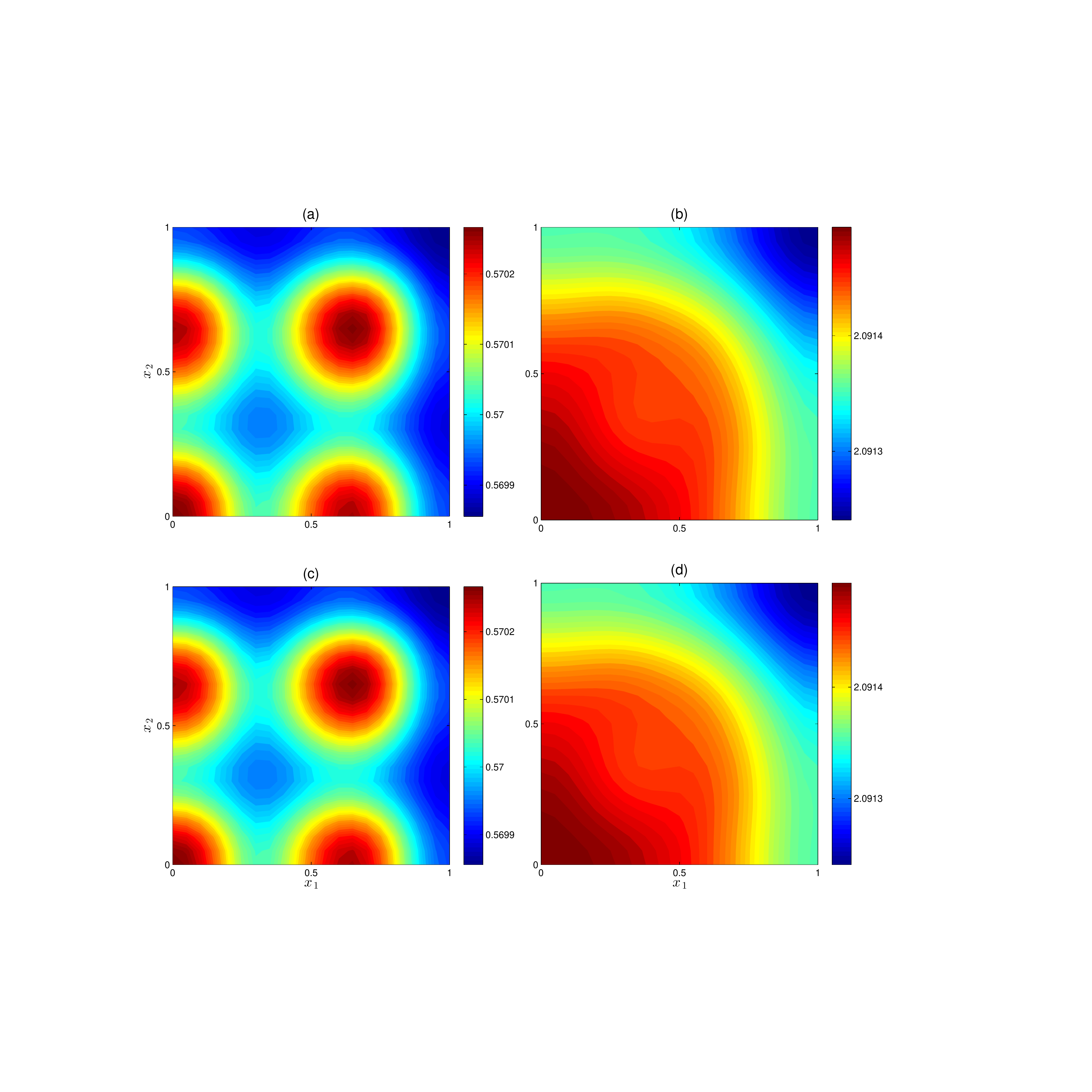}
    \caption{Mean of solution fields obtained using both ISP method implementations. Subfigures (a, c) and and (b, d) correspond to $T$ and $\varphi$ respectively, while (a, b) and (c, d) correspond to the reduced and standard ISP method implementations respectively.}
\end{figure}

\begin{figure}
    \hypertarget{fig3}{}
    \centering
    \includegraphics[bb=80bp 240bp 1100bp 905bp,scale=0.4]{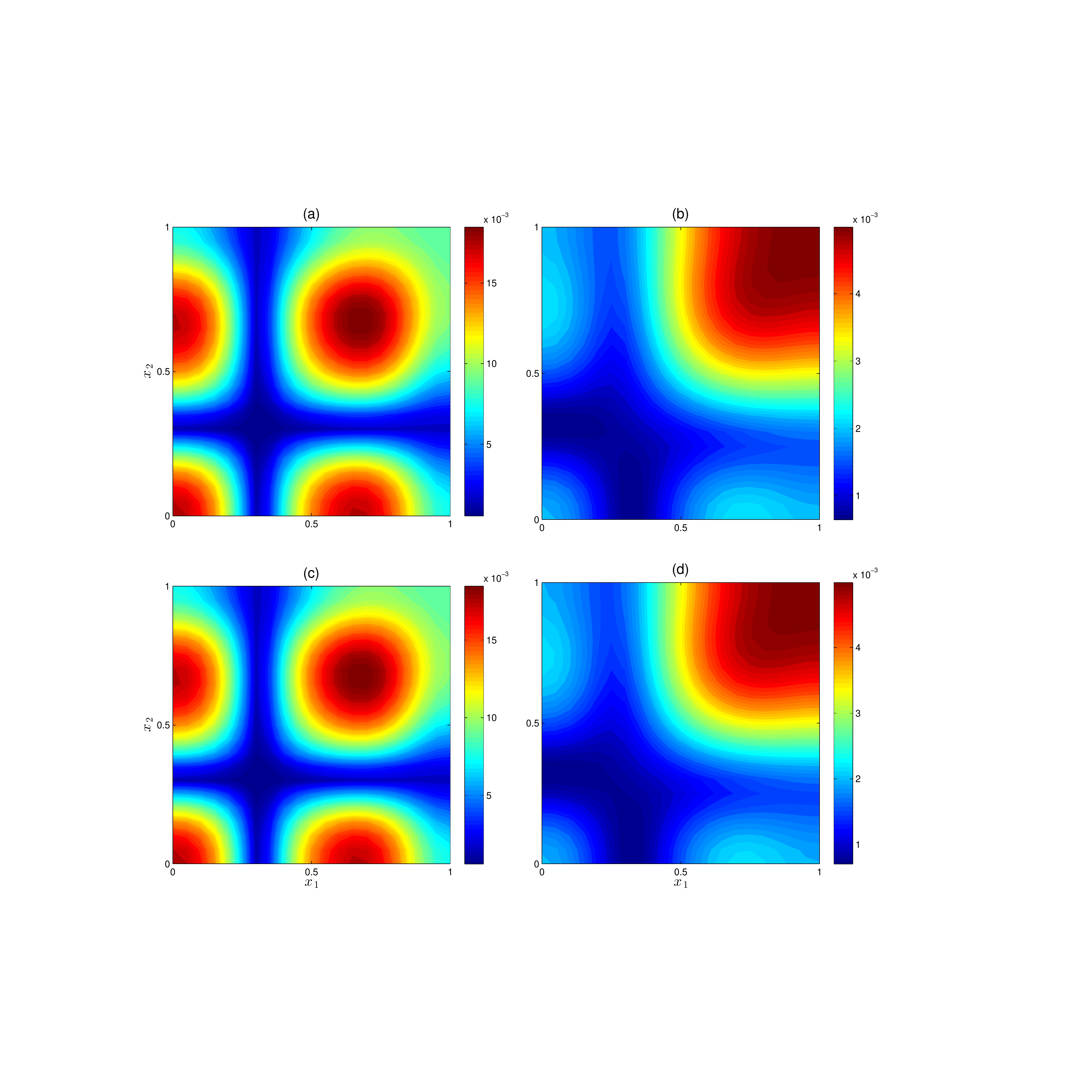}
    \caption{Standard deviation of solution fields obtained using both ISP method implementations. Subfigures (a, c) and and (b, d) correspond to $T$ and $\varphi$ respectively, while (a, b) and (c, d) correspond to the reduced and standard ISP method implementations respectively.}
\end{figure}

\begin{table}[htbp]
\hypertarget{tab2}{}
\caption{Comparison of the main effect and interaction effect sensitivities (percentage) obtained in the standard and reduced ISP method implementations for the Thermal-Neutronics problem.}
\begin{center}\scriptsize
\renewcommand{\arraystretch}{1.3}
\begin{tabular}{cccccc|cccccc}
\toprule 
\multicolumn{6}{c}{Standard ISP} \vline & \multicolumn{6}{c}{Reduced ISP} \tabularnewline
{$\hat{\mathcal{V}}_{1,1}$} & {$\hat{\mathcal{V}}_{1,2}$} & {$\hat{\mathcal{V}}_{1,12}$} & {$\hat{\mathcal{V}}_{2,1}$} & {$\hat{\mathcal{V}}_{2,2}$} & {$\hat{\mathcal{V}}_{2,12}$}\vline & {$\hat{\mathcal{V}}_{1,1}$} & {$\hat{\mathcal{V}}_{1,2}$} & {$\hat{\mathcal{V}}_{1,12}$} & {$\hat{\mathcal{V}}_{2,1}$} & {$\hat{\mathcal{V}}_{2,2}$} & {$\hat{\mathcal{V}}_{2,12}$} \tabularnewline
\midrule
\midrule 
{$91.1$} & {$8.9$} & {$0.0$} & {$97.6$} & {$2.7$} & {$0.1$} & {$92.6$} & {$7.4$} & {$0.0$} & {$97.8$} & {$2.2$} & {$0.0$} \tabularnewline
\bottomrule
\end{tabular}\par\end{center}
\end{table} 

\begin{table}[htbp]
\hypertarget{tab3}{}
\caption{Comparison of the approximation errors and computational costs
(seconds) obtained in the standard and reduced ISP method implementations for the Thermal-Neutronics problem. The
reduced dimensions and orders are also listed here.}
\begin{center}\scriptsize
\renewcommand{\arraystretch}{1.3}
\begin{tabular}{cc|cc|cccccccc|c}
\toprule 
 &  & \multicolumn{2}{c}{Standard ISP} \vline & \multicolumn{8}{c}{Reduced ISP}\vline & \tabularnewline
{$s_{1}$, $s_{2}$} & {$p$} & {$\varepsilon_{s}$} & {$\mathcal{C}_{s}$ } & {$d_{1}$} & {$\tilde{p}_{1}$} & {$\tilde{Q}_{1}$} & {$d_{2}$} & {$\tilde{p}_{2}$} & {$\tilde{Q}_{2}$} & {$\varepsilon_{r}$} & {$\mathcal{C}_{r}$} & {$\mathcal{C}_{s}/ \mathcal{C}_{r}$} \tabularnewline
\midrule
\midrule 
 & {$2$}  & {$7.6\times10^{-4}$} & {$34$} & {$2$} & {$1$} & {$15$} & {$1$} & {$1$} & {$5$} & {$3.3\times10^{-2}$} & {$19$} & {$1.8$}\tabularnewline
\cmidrule{2-13} 
{$3$} & {$3$} & {$1.8\times10^{-4}$} & {$193$} & {$2$} & {$1$} & {$15$} & {$1$} & {$1$} & {$5$} & {$9.7\times10^{-3}$} & {$47$} & {$4.1$}\tabularnewline
\cmidrule{2-13}  
 & {$4$} & {$8.7\times10^{-5}$} & {$1263$} & {$2$} & {$1$} & {$15$} & {$1$} & {$1$} & {$5$} & {$8.5\times10^{-3}$} & {$130$} & {$9.7$}\tabularnewline
\midrule 
 & {$2$} & {$1.9\times10^{-3}$} & {$65$} & {$2$} & {$1$} & {$15$} & {$1$} & {$1$} & {$5$} & {$3.4\times10^{-2}$} & {$25$} & {$2.6$}\tabularnewline
\cmidrule{2-13}  
{$4$} & {$3$} & {$3.7\times10^{-4}$} & {$878$} & {$2$} & {$1$} & {$15$} & {$1$} & {$1$} & {$5$} & {$1.0\times10^{-2}$} & {$117$} & {$7.5$}\tabularnewline
\cmidrule{2-13}  
 & {$4$} & {$1.1\times10^{-4}$} & {$8110$} & {$2$} & {$1$} & {$15$} & {$1$} & {$1$} & {$5$} & {$9.4\times10^{-3}$} & {$529$} & {$15.3$}\tabularnewline
\bottomrule
\end{tabular}
\par\end{center}
\end{table}

\begin{figure}
    \hypertarget{fig4}{}
    \centering
    \includegraphics[bb=80bp 240bp 1100bp 905bp,scale=0.4]{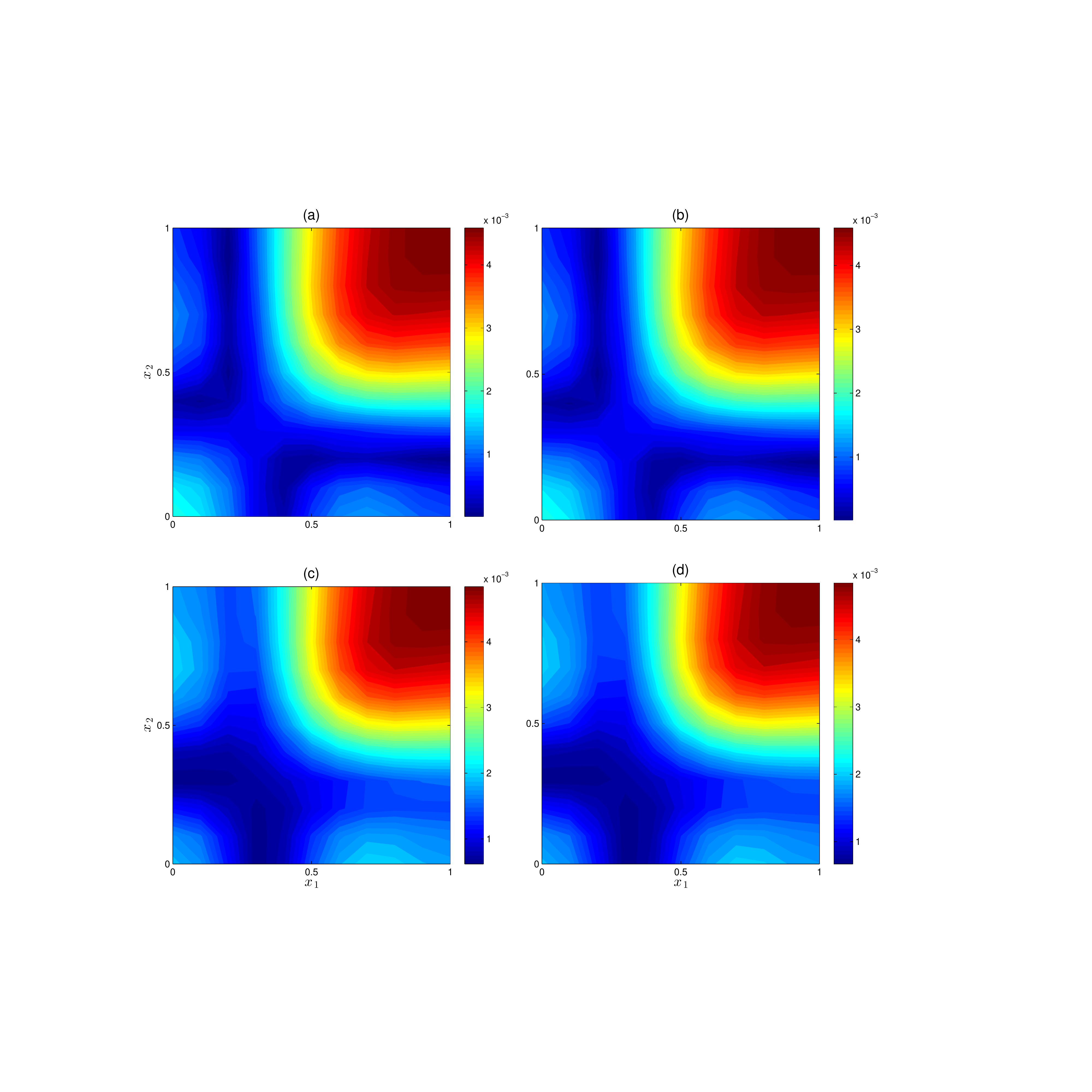}
    \caption{Comparison of computed standard deviations on coarser mesh. Subfigures (a, b, c) correspond to the reduced ISP method implementation while (d) corresponds to the standard ISP method implementation. The tolerance values are as follows. (a) $\epsilon_{1,\mathrm{dim}}=0.02$, $\epsilon_{2,\mathrm{dim}}=0.1$, (b) $\epsilon_{1,\mathrm{dim}}=0.01$, $\epsilon_{2,\mathrm{dim}}=0.1$, (c) $\epsilon_{1,\mathrm{dim}}=0.02$, $\epsilon_{2,\mathrm{dim}}=0.05$.}
\end{figure}

The highest speedup factor observed is $\approx 14.3$. Moreover, the approximation errors in the standard ISP method implementation
are observed to decay exponentially, which indicates a high degree
of regularity in the stochastic solutions. Furthermore, for the reduced ISP method implementation, the asymptotic upper bound 
predicted in \hyperlink{sec36}{\S 3.6}, is observed in its approximation
errors.

\subsection{Boussinesq flow problem.}

The Boussinesq model describes thermally driven, incompressible flows
and is widely used in oceanic and atmospheric modeling {[}\hyperlink{ref39}{39}, \hyperlink{ref40}{40}{]}.
Here, we consider a multi-physics setup with uncertain fluid properties
and boundary conditions.

\subsubsection{Model setup}

Let $\Omega\equiv\left(0,1\right)_{x_{1}}\times\left(0,1\right)_{x_{2}}$
denote the spatial domain and \textbf{$\boldsymbol{u}=\left[\begin{array}{cc}
u_{1} & u_{2}\end{array}\right]^{\mathbf{T}}$} , $p$ $,T$ denote the non-dimensional {[}\hyperlink{ref40}{40}{]} fluid velocity,
pressure and temperature respectively. The governing equations for
the fluid variables are as follows. $\forall\boldsymbol{\xi}\in\Xi$,

\begin{align}
\boldsymbol{\nabla}^{\mathbf{T}}\boldsymbol{u}\left(\boldsymbol{x},\boldsymbol{\xi}\right) & =0,\nonumber \\
\left(\boldsymbol{u}\left(\boldsymbol{x},\boldsymbol{\xi}\right)^{\mathbf{T}}\boldsymbol{\nabla}\right)\boldsymbol{u}\left(\boldsymbol{x},\boldsymbol{\xi}\right)+\boldsymbol{\nabla}p\left(\boldsymbol{x},\boldsymbol{\xi}\right)\nonumber \\
-\mathrm{Pr}\boldsymbol{\nabla}^{\mathbf{T}}\boldsymbol{\nabla}\boldsymbol{u}\left(\boldsymbol{x},\boldsymbol{\xi}\right)-\mathrm{Pr}\mathrm{Ra}\left(\boldsymbol{x},\boldsymbol{\xi}_{1}\left(\boldsymbol{\xi}\right)\right)T\left(\boldsymbol{x},\boldsymbol{\xi}\right)\boldsymbol{e}_{2} & =\boldsymbol{0},\nonumber \\
\left(\boldsymbol{u}\left(\boldsymbol{x},\boldsymbol{\xi}\right)^{\mathbf{T}}\boldsymbol{\nabla}\right)T\left(\boldsymbol{x},\boldsymbol{\xi}\right)-\boldsymbol{\nabla}^{\mathbf{T}}\boldsymbol{\nabla}T\left(\boldsymbol{x},\boldsymbol{\xi}\right) & =0, & \boldsymbol{x}\in\Omega,
\end{align}
with homogenous Dirichlet boundary conditions for $\boldsymbol{u}$
and Neumann boundary conditions for $p$ at all boundaries. The boundary
conditions for temperature, as shown in \hyperlink{fig5}{Figure 5}, are as follows.
\begin{align}
\frac{\partial T}{\partial x_{2}}\left(x_{1},0,\boldsymbol{\xi}\right)=\frac{\partial T}{\partial x_{2}}\left(x_{1},1,\boldsymbol{\xi}\right) & =0, & x_{1}\in\left[0,1\right],\nonumber \\
T\left(0,x_{2},\boldsymbol{\xi}\right)-T_{h}\left(x_{2},\boldsymbol{\xi}_{2}\right)=T\left(1,x_{2},\boldsymbol{\xi}\right) & =0, & x_{2}\in\left[0,1\right].
\end{align}

\begin{figure}
    \hypertarget{fig5}{}
    \centering
    \includegraphics[bb=125bp 500bp 475bp 720bp,clip,scale=0.65]{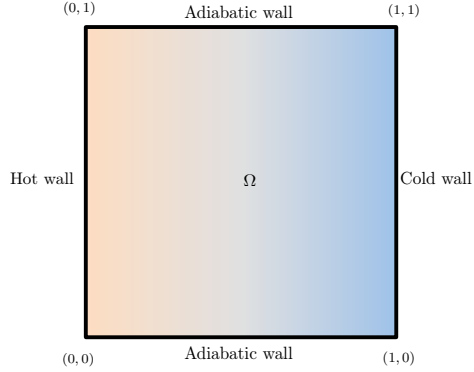}  
    \caption{Computational domain for the Boussinesq flow problem.}
    \label{fig:fig5}
\end{figure}

Here, $\boldsymbol{e}_{2}$ denotes $\left[\begin{array}{cc}
0 & 1\end{array}\right]^{\mathbf{T}}$ while $\mathrm{Pr}$ and $\mathrm{Ra}$ denote the Prandtl and Rayleigh
numbers respectively. $T_{h}$ denotes the hot-wall temperature such
that $\forall x_{2}\in\left[0,1\right],\boldsymbol{\xi}_{2}\in\Xi_{2}$,
\begin{equation}
T_{h}\left(x_{2},\boldsymbol{\xi}_{2}\right)=\bar{T}_{h}+h\left(x_{2},\boldsymbol{\xi}_{2}\right)\sin^{2}\left(\pi x_{2}\right),
\end{equation}
where $\bar{T}_{h}$ is the mean hot-wall temperature and $h$ denotes
the perturbation amplitude.

In this study, $\mathrm{Ra}$ and $h$ are assumed to be independent
random fields, and modeled using the following KL expansions. $\forall\boldsymbol{x}\in\Omega,\boldsymbol{\xi}_{1}\in\Xi_{1}$,
\begin{equation}
\mathrm{Ra}\left(\boldsymbol{x},\boldsymbol{\xi}_{1}\right)=\bar{\mathrm{Ra}}+\sqrt{3}\delta_{\mathrm{Ra}}\sum_{j=1}^{s_{1}}\gamma_{\mathrm{Ra},j}\left(\boldsymbol{x}\right)\xi_{1j},
\end{equation}
where $\bar{\mathrm{Ra}}$ denotes the mean of $\mathrm{Ra}$ and
$\left\{ \xi_{1j}\sim U\left[-1,1\right]\right\} _{j=1}^{s_{1}}$
are $i.i.d.$ random variables. Similarly, $\forall x_{2}\in\left(0,1\right),\boldsymbol{\xi}_{2}\in\Xi_{2}$,
\begin{equation}
h\left(x_{2},\boldsymbol{\xi}_{2}\right)=\sqrt{3}\delta_{h}\sum_{j=1}^{s_{2}}\gamma_{h,j}\left(\boldsymbol{x}\right)\xi_{2j},
\end{equation}
where $\left\{ \xi_{2j}\sim U\left[-1,1\right]\right\} _{j=1}^{s_{2}}$
are $i.i.d.$ random variables. Moreover, we assume that both $\mathrm{Ra}$
and $h$ have exponential covariance kernels 
\begin{align}
C_{\mathrm{Ra}}\left(\boldsymbol{x},\boldsymbol{y}\right) & =\delta_{\mathrm{Ra}}^{2}\exp\left(-\frac{\left\Vert \boldsymbol{x}-\boldsymbol{y}\right\Vert _{1}}{l_{\mathrm{Ra}}}\right), & \boldsymbol{x},\boldsymbol{y}\in\Omega,\nonumber \\
C_{h}\left(x_{2},y_{2}\right) & =\delta_{h}^{2}\exp\left(-\frac{\left|x_{2}-y_{2}\right|}{l_{h}}\right), & x_{2},y_{2}\in\left[0,1\right],
\end{align}
where $\delta_{\mathrm{Ra}}$, $\delta_{h}$ denote the respective
coefficients of variations, and $l_{\mathrm{Ra}}$, $l_{h}$ denote
the respective correlation lengths. The analytic expressions for $\left\{ \gamma_{\mathrm{Ra},j}\right\} _{j>0}$
and $\left\{ \gamma_{h,j}\right\} _{j>0}$ are provided in Appendix
A. 

The pressure Poisson equation 
\begin{equation}
\boldsymbol{\nabla}^{\mathbf{T}}\boldsymbol{\nabla}p\left(\boldsymbol{x},\boldsymbol{\xi}\right)+\boldsymbol{\nabla}^{\mathbf{T}}\left(\left(\boldsymbol{u}\left(\boldsymbol{x},\boldsymbol{\xi}\right)^{\mathbf{T}}\boldsymbol{\nabla}\right)\boldsymbol{u}\left(\boldsymbol{x},\boldsymbol{\xi}\right)-\mathrm{Pr}\mathrm{Ra}\left(\boldsymbol{x},\boldsymbol{\xi}_{1}\right)T\left(\boldsymbol{x},\boldsymbol{\xi}\right)\boldsymbol{e}_{2}\right)=0
\end{equation}
is solved in place of the continuity equation, to close the momentum
component of the PDE system.

Each component of the PDE system is spatially discretized using a
finite volume method, with linear central-differencing schemes [\hyperlink{ref42}{42}],
on a uniform grid with $m\times m$ cells. Let $\boldsymbol{u}_{1}^{\prime},\boldsymbol{u}_{2}^{\prime},\boldsymbol{p}^{\prime},\boldsymbol{t}^{\prime}\in\mathbb{R}^{m^{2}}$
denote the respective vectors of cell-centroidal horizontal velocity,
vertical velocity, pressure and temperature, which solve the nonlinear
system 
\hypertarget{eq416}{}
\begin{align}
\left(\boldsymbol{K}_{u}+\boldsymbol{A}\left(\boldsymbol{u}_{1}^{\prime},\boldsymbol{u}_{2}^{\prime}\right)\right)\boldsymbol{u}_{1}^{\prime}+\boldsymbol{B}_{1}\boldsymbol{p}^{\prime} & =\boldsymbol{0},\nonumber \\
\left(\boldsymbol{K}_{u}+\boldsymbol{A}\left(\boldsymbol{u}_{1}^{\prime},\boldsymbol{u}_{2}^{\prime}\right)\right)\boldsymbol{u}_{2}^{\prime}+\boldsymbol{B}_{2}\boldsymbol{p}^{\prime}-\boldsymbol{R}\left(\boldsymbol{\xi}_{1}\right)\boldsymbol{t}^{\prime} & =\boldsymbol{0},\nonumber \\
\boldsymbol{K}_{p}\boldsymbol{p}^{\prime}+\boldsymbol{C}_{1}\left(\boldsymbol{u}_{1}^{\prime},\boldsymbol{u}_{2}^{\prime}\right)\boldsymbol{u}_{1}^{\prime}+\boldsymbol{C}_{2}\left(\boldsymbol{u}_{1}^{\prime},\boldsymbol{u}_{2}^{\prime}\right)\boldsymbol{u}_{2}^{\prime}-\boldsymbol{S}\left(\boldsymbol{\xi}_{1}\right)\boldsymbol{t}^{\prime} & =\boldsymbol{0},\nonumber \\
\left(\boldsymbol{K}_{T}+\boldsymbol{A}\left(\boldsymbol{u}_{1}^{\prime},\boldsymbol{u}_{2}^{\prime}\right)\right)\boldsymbol{t}^{\prime}-\boldsymbol{h}\left(\boldsymbol{\xi}_{2}\right) & =\boldsymbol{0}.
\end{align}
where each term in \hyperlink{eq416}{Eq. 4.16} denotes its respective discretized operator
in the coupled PDE system. Subsequently, we formulate a modular multi-physics
setup, as per \hyperlink{eq21}{Eq. 2.1}, by separating the momentum and energy components
of the coupled algebraic system. Let $\boldsymbol{u}_{1}=\left[\boldsymbol{u}_{1}^{\prime};\boldsymbol{u}_{2}^{\prime};\boldsymbol{p}^{\prime}\right]\in\mathbb{R}^{n_{1}}\equiv\mathbb{R}^{3m^{2}}$,
$\boldsymbol{u}_{2}=\boldsymbol{t}^{\prime}\in\mathbb{R}^{n_{2}}\equiv\mathbb{R}^{m^{2}}$
denote the respective solution variables in the modular algebraic
system. Therefore, the component residuals are formulated as follows.
\begin{align}
\boldsymbol{f}_{1}\left(\boldsymbol{u}_{1};\boldsymbol{u}_{2},\boldsymbol{\xi}_{1}\right) & =\left[\begin{array}{cc}
\boldsymbol{K}_{u}+\boldsymbol{A}\left(\boldsymbol{u}_{1}^{\prime}\left(\boldsymbol{u}_{1}\right),\boldsymbol{u}_{2}^{\prime}\left(\boldsymbol{u}_{1}\right)\right) & \boldsymbol{0}\\
\boldsymbol{0} & \boldsymbol{K}_{u}+\boldsymbol{A}\left(\boldsymbol{u}_{1}^{\prime}\left(\boldsymbol{u}_{1}\right),\boldsymbol{u}_{2}^{\prime}\left(\boldsymbol{u}_{1}\right)\right)\\
\boldsymbol{C}_{1}\left(\boldsymbol{u}_{1}^{\prime}\left(\boldsymbol{u}_{1}\right),\boldsymbol{u}_{2}^{\prime}\left(\boldsymbol{u}_{1}\right)\right) & \boldsymbol{C}_{2}\left(\boldsymbol{u}_{1}^{\prime}\left(\boldsymbol{u}_{1}\right),\boldsymbol{u}_{2}^{\prime}\left(\boldsymbol{u}_{1}\right)\right)
\end{array}\right.\nonumber \\
 & \left.\begin{array}{c}
\boldsymbol{B}_{1}\\
\boldsymbol{B}_{2}\\
\boldsymbol{K}_{p}
\end{array}\right]\boldsymbol{u}_{1}-\left[\begin{array}{c}
\boldsymbol{0}\\
\boldsymbol{R}\left(\boldsymbol{\xi}_{1}\right)\\
\boldsymbol{S}\left(\boldsymbol{\xi}_{1}\right)
\end{array}\right]\boldsymbol{u}_{2},\nonumber \\
\boldsymbol{f}_{2}\left(\boldsymbol{u}_{2};\boldsymbol{u}_{1},\boldsymbol{\xi}_{2}\right) & =\left(\boldsymbol{K}_{T}+\boldsymbol{A}\left(\boldsymbol{u}_{1}^{\prime}\left(\boldsymbol{u}_{1}\right),\boldsymbol{u}_{2}^{\prime}\left(\boldsymbol{u}_{1}\right)\right)\right)\boldsymbol{u}_{2}-\boldsymbol{h}\left(\boldsymbol{\xi}_{2}\right).
\end{align}

\begin{table}[htbp]
\hypertarget{tab4}{}
\caption{Deterministic parameter values in the Boussinesq flow problem.}
\begin{center}\scriptsize
\renewcommand{\arraystretch}{1.3}
\begin{tabular}{ccccccc}
\toprule 
$\mathrm{Pr}$ & $\bar{\mathrm{Ra}}$ & $\bar{T}_{h}$ & $\delta_{\mathrm{Ra}}$ & $\delta_{h}$ & $l_{\mathrm{Ra}}$ & $l_{h}$\tabularnewline
\midrule
\midrule 
$0.71$ & $1000$ & $1$ & $200$ & $0.5$ & $0.5$ & $0.5$\tabularnewline
\bottomrule
\end{tabular}
\par\end{center}
\end{table}

\hyperlink{tab4}{Table 4} lists the numerical values of the deterministic parameters
used in this study. The quantities of interest are the probability
density functions of the (scaled) kinetic energy $K$ and thermal
energy $E$, which are defined as follows. $\forall\boldsymbol{\xi}\in\Xi$,

\begin{equation}
K\left(\boldsymbol{\xi}\right)=\frac{1}{2}\left(\int_{\Omega}u_{1}\left(\boldsymbol{x},\boldsymbol{\xi}\right)^{2}d\boldsymbol{x}+\int_{\Omega}u_{2}\left(\boldsymbol{x},\boldsymbol{\xi}\right)^{2}d\boldsymbol{x}\right),\ E\left(\boldsymbol{\xi}\right)=\int_{\Omega}T\left(\boldsymbol{x},\boldsymbol{\xi}\right)d\boldsymbol{x},
\end{equation}
and the statistics of the fluid velocity and temperature.

\subsubsection{Modular deterministic solver: Setup and verification}

Following the same procedure as in the previous numerical example,
a modular deterministic solver was implemented, wherein the BGS iterations
are formulated as follows.
\begin{align}
\boldsymbol{u}_{1}^{\ell+1} & =\boldsymbol{u}_{1}^{\ell}-\left(\frac{\partial\boldsymbol{f}_{1}}{\partial\boldsymbol{u}_{1}}\left(\boldsymbol{u}_{1}^{\ell};\boldsymbol{u}_{2}^{\ell},\boldsymbol{\xi}_{1}\right)\right)^{-1}\boldsymbol{f}_{1}\left(\boldsymbol{u}_{1}^{\ell};\boldsymbol{u}_{2}^{\ell},\boldsymbol{\xi}_{1}\right),\nonumber \\
\boldsymbol{u}_{2}^{\ell+1} & =\boldsymbol{u}_{2}^{\ell}-\left(\frac{\partial\boldsymbol{f}_{2}}{\partial\boldsymbol{u}_{2}}\left(\boldsymbol{u}_{2}^{\ell};\boldsymbol{u}_{1}^{\ell+1},\boldsymbol{\xi}_{2}\right)\right)^{-1}\boldsymbol{f}_{2}\left(\boldsymbol{u}_{2}^{\ell};\boldsymbol{u}_{1}^{\ell+1},\boldsymbol{\xi}_{2}\right).
\end{align}
Similar formulations are used in the standard and reduced ISP method
implementations, wherein the corresponding linear systems are solved
using LSQR, and mean based preconditioners. Once again, in each implementation,
the components of the solver were developed as $\mathtt{MATLAB}^{\text{\texttrademark}}$
function modules.

\subsubsection{ISP based uncertainty propagation}

We implemented both ISP based uncertainty propagation methods to compute
the quantities of interest. Subsequently, for $m=25$, $s_{1}=s_{2}=3$
and $p=4$, the probability density function of $K$ and $E$, and
the first two moments of $u_{1},u_{2},T$, were computed and compared.
The results are shown in \hyperlink{fig6}{Figure 6}, \hyperlink{fig7}{Figure 7} and \hyperlink{fig8}{Figure 8} respectively.
Moreover, in the reduced ISP method implementation, the joint probability density functions of the reduced random
variables was computed at the last iteration. The results are illustrated in \hyperlink{fig9}{Figure 9}.
The tolerance values used in this study are $\epsilon_{1,\mathrm{dim}}=0.01,\epsilon_{2,\mathrm{dim}}=0.02$
and $\epsilon_{1,\mathrm{ord}}=\epsilon_{2,\mathrm{ord}}=10^{-4}$. 

For various values of $s_{1}$, $s_{2}$ and $p$, the errors
$\varepsilon_{s},\varepsilon_{r}$ and overall costs $\mathcal{C}_{s}$,
$\mathcal{C}_{r}$ (wall-times) of both methods are compared in \hyperlink{tab5}{Table 5}. In both ISP method implementations, the rate of convergence of the BGS iterations is observed to be invariant with respect to the stochastic dimensions $s_{1},s_{2}$ and gPC order $p$. For $\epsilon_{\mathrm{BGS}}=10^{-6}$, convergence is achieved in $9$ iterations.

\begin{figure}
    \hypertarget{fig6}{}
    \centering
    \includegraphics[bb=185bp 400bp 1050bp 750bp,clip,scale=0.5]{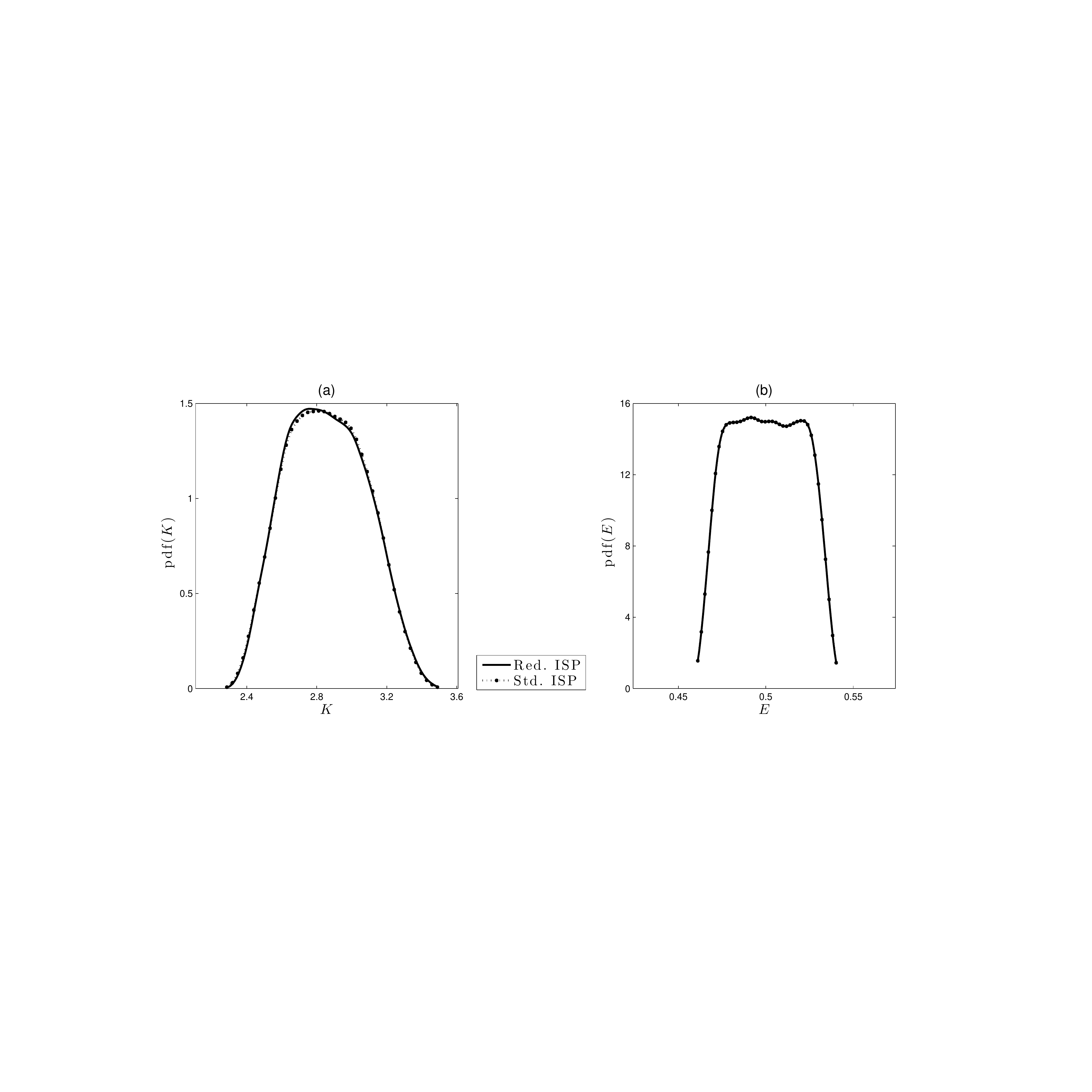} 
    \caption{Probability density function of the fluid energies computed
using both ISP method implementations. Subfigure (a) corresponds to the kinetic energy $K$, while (b) corresponds to the thermal energy $E$. The densities were computed using the KDE method with $10^{5}$ samples.}
    \label{fig:fig6}
\end{figure}

\begin{figure}
    \hypertarget{fig7}{}
    \centering
    \includegraphics[bb=140bp 360bp 1100bp 810bp,scale=0.47]{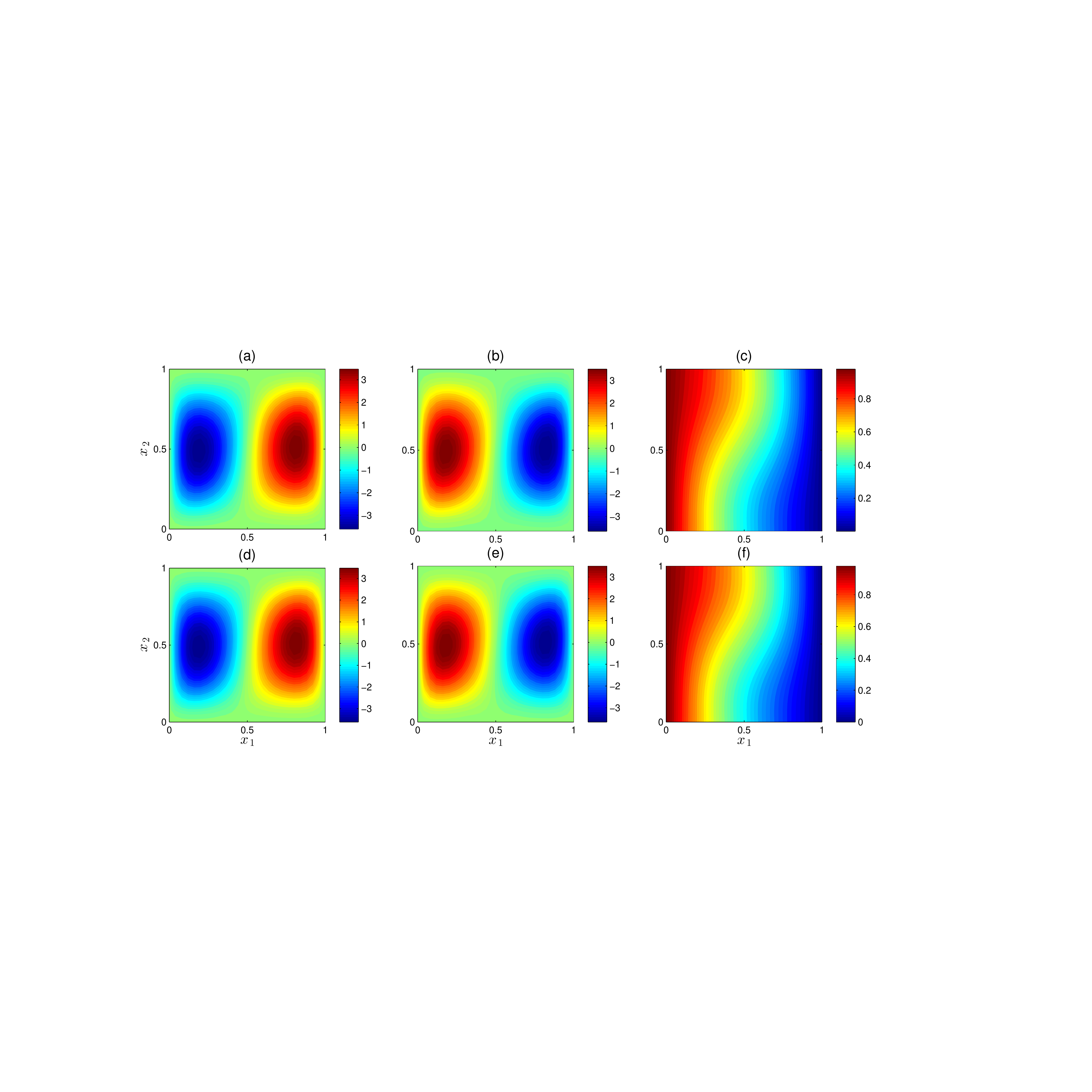}
    \caption{Mean of solution fields obtained using both ISP method implementations. Subfigures (a, d), (b, e) and (c, f) correspond to $u_1$, $u_2$ and $T$ respectively, while (a, b, c) and (d, e, f) correspond to the reduced and standard ISP method implementations respectively.}
    \label{fig:fig7}

    \hypertarget{fig8}{}
    \centering
    \includegraphics[bb=140bp 360bp 1100bp 810bp,scale=0.47]{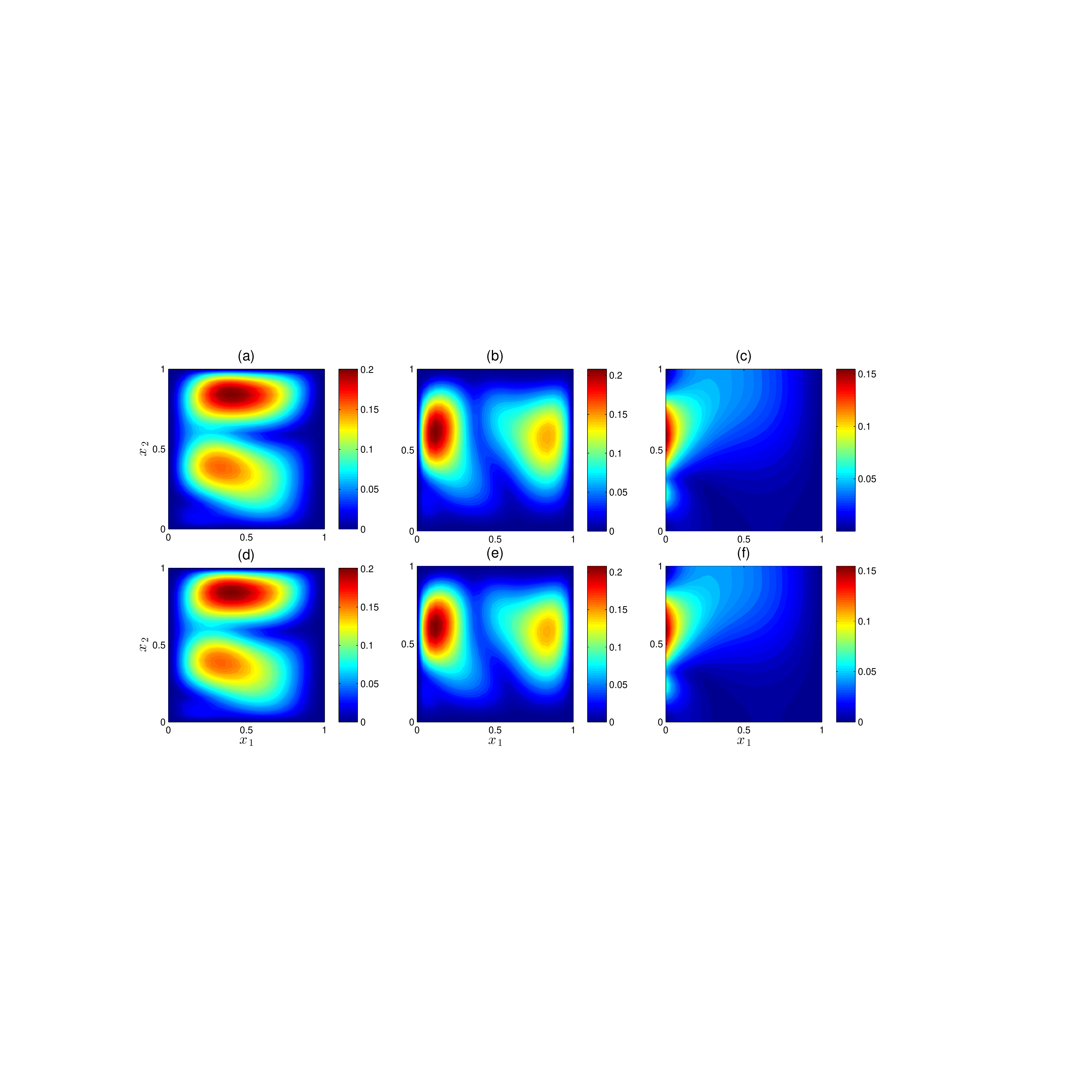}
    \caption{Standard deviation of solution fields obtained using both ISP method implementations. Subfigures (a, d), (b, e) and (c, f) correspond to $u_1$, $u_2$ and $T$ respectively, while (a, b, c) and (d, e, f) correspond to the reduced and standard ISP method implementations respectively.}
    \label{fig:fig8}
\end{figure}

\begin{figure}
    \hypertarget{fig9}{}
    \centering
    \includegraphics[bb=185bp 390bp 1050bp 750bp,clip,scale=0.5]{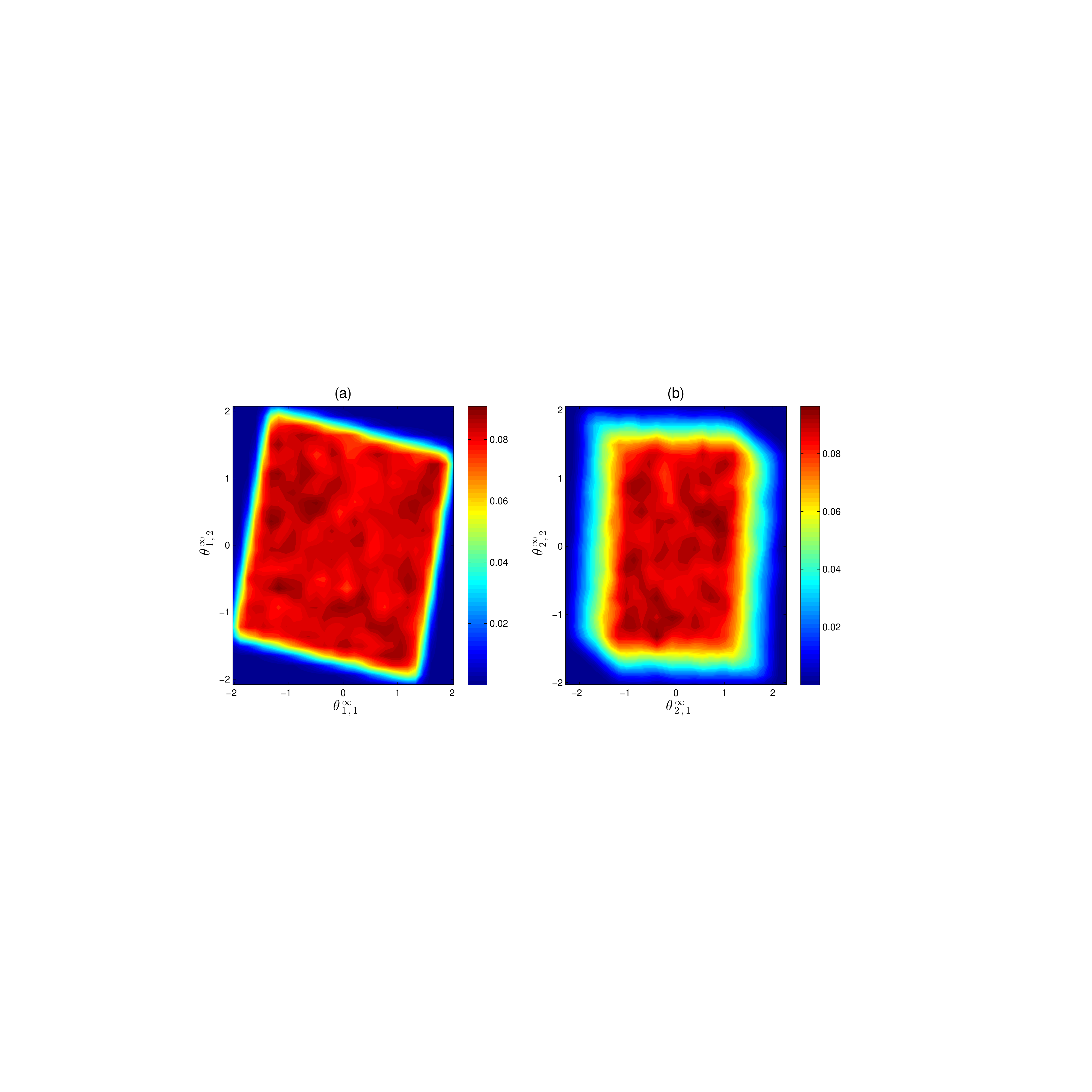} 
    \caption{Probability density function of the fluid energies computed
using both ISP method implementations. Subfigure (a) corresponds to the kinetic energy $K$, while (b) corresponds to the thermal energy $E$. The densities were computed using the KDE method with $10^{5}$ samples.}
    \label{fig:fig6}
\end{figure}

\begin{table}[htbp]
\hypertarget{tab5}{}
\caption{Comparison of the approximation errors and computational costs
(seconds) obtained in the standard and reduced ISP method implementations for the Boussinesq flow problem. The
reduced dimensions and orders are also listed here.}
\begin{center}\scriptsize
\renewcommand{\arraystretch}{1.3}
\begin{tabular}{cc|cc|cccccccc|c}
\toprule 
 &  & \multicolumn{2}{c}{Standard ISP} \vline & \multicolumn{8}{c}{Reduced ISP}\vline & \tabularnewline
{$s_{1}$, $s_{2}$} & {$p$} & {$\varepsilon_{s}$} & {$\mathcal{C}_{s}$ } & {$d_{1}$} & {$\tilde{p}_{1}$} & {$\tilde{Q}_{1}$} & {$d_{2}$} & {$\tilde{p}_{2}$} & {$\tilde{Q}_{2}$} & {$\varepsilon_{r}$} & {$\mathcal{C}_{r}$} & {$\mathcal{C}_{s}/ \mathcal{C}_{r}$} \tabularnewline
\midrule
\midrule 
 & {$2$}  & {$8.6\times10^{-4	}$} & {$46$} & {$2$} & {$1$} & {$15$} & {$2$} & {$1$} & {$15$} & {$2.3\times10^{-2}$} & {$32$} & {$1.4$}\tabularnewline
\cmidrule{2-13} 
{$3$} & {$3$} & {$2.3\times10^{-4}$} & {$334$} & {$2$} & {$1$} & {$15$} & {$2$} & {$1$} & {$15$} & {$5.7\times10^{-3}$} & {$93$} & {$3.6$}\tabularnewline
\cmidrule{2-13}  
 & {$4$} & {$8.7\times10^{-5}$} & {$2713$} & {$2$} & {$1$} & {$15$} & {$2$} & {$1$} & {$15$} & {$3.6\times10^{-3}$} & {$326$} & {$8.3$}\tabularnewline
\midrule 
 & {$2$} & {$1.9\times10^{-3}$} & {$98$} & {$2$} & {$1$} & {$15$} & {$2$} & {$1$} & {$15$} & {$2.4\times10^{-2}$} & {$46$} & {$2.1$}\tabularnewline
\cmidrule{2-13}  
{$4$} & {$3$} & {$4.7\times10^{-4}$} & {$1710$} & {$2$} & {$1$} & {$15$} & {$2$} & {$1$} & {$15$} & {$6.2\times10^{-3}$} & {$298$} & {$5.7$}\tabularnewline
\cmidrule{2-13}  
 & {$4$} & {$1.3\times10^{-4}$} & {$16487$} & {$2$} & {$1$} & {$15$} & {$2$} & {$1$} & {$15$} & {$3.8\times10^{-3}$} & {$1259$} & {$13.1$}\tabularnewline
\bottomrule
\end{tabular}
\par\end{center}
\end{table}

The highest speedup factor observed is $\approx 12.1$. The convergence rate of the approximation 
error in the standard ISP method decays exponentially, which indicates a high degree of regularity in 
the solutions. Moreover, the predicted asymptotic upper bound in the approximation error in the
reduced ISP method is observed as well.

\section{Conclusions and future work}
A reduced ISP based uncertainty propagation method for stochastic
multi-physics models is presented, and compared against the standard ISP method using two numerical examples. 
At the expense of relatively small approximation errors, the reduced ISP method exhibited
significant speedup over the standard ISP method.
This can be primarily attributed to the dimension
and order reduction steps for constructing an approximation of the input data in each module,
which mitigates the rate of exponential growth of overall computational costs.

The two most significant overheads observed while
implementing the reduced ISP algorithm
were the intermediate SVD computation towards constructing the reduced dimensional approximation space, and the computation of the global
gPC coefficient matrices from the reduced and modular gPC coefficient matrices. Both these overheads can be
easily eliminated due to the embarrassingly parallel nature of the computations involved.
Moreover, alternatives, for example, Schur complement based elimination, to the BGS partitioning approach can be explored
for reducing the overall computational costs.

Since the global gPC coefficients still need to be computed, manipulated and stored, the applicability of our proposed method is limited to models in which the global stochastic dimension is manageably low ($<20$). If a particular module contributes a large number of uncertainties to the global set,
a Monte-Carlo based sampling approach can be employed in that module alone, while
other modules could still afford the use of spectral methods. This method
has been recently demonstrated in [\hyperlink{ref43}{43}]. 
Moreover, for tackling models in which the solution regularity is low,
the proposed method can also be easily adapted towards multi-element gPC
[\hyperlink{ref44}{44}] and discontinuous wavelet [\hyperlink{ref45}{45}] based spectral representations. 

Other active areas being currently investigated also include derivative or active-subspace based dimension reduction
methods [\hyperlink{ref45}{45}] as alternative reduction strategies. Moreover, studying the effects of uncertainties on the approximation
errors in a multi-physics simulation context is a subject of future research.

\section*{Acknowledgement}
This research was funded by the US Department of Energy, Office of
Advanced Computing Research and Applied Mathematics Program and partially
funded by the US Department of Energy NNSA ASC Program.

\newpage

\section*{Appendix A - Karhunen-Loeve expansion}

\subsection*{Exponential covariance kernel}

Given an $n$-dimensional spatial domain $\Omega\subseteq\mathbb{R}^{n}$,
let $C_{u}:\Omega\times\Omega\rightarrow\mathbb{R}^{+}$ denote the exponential
covariance kernel of a spatially varying random field $u\in\Omega\rightarrow\mathbb{R}$.
Therefore, $C_{u}:\forall\boldsymbol{x}=\left[\begin{array}{ccc}
x_{1} & \cdots & x_{n}\end{array}\right]^{\mathbf{T}},\boldsymbol{y}=\left[\begin{array}{ccc}
y_{1} & \cdots & y_{n}\end{array}\right]^{\mathbf{T}}\in\Omega$, 
\begin{align*}
C_{u}\left(\boldsymbol{x},\boldsymbol{y}\right)=\exp\left(-\frac{\left\Vert \boldsymbol{x}-\boldsymbol{y}\right\Vert _{1}}{l}\right)=\prod_{j=1}^{n}\exp\left(-\frac{\left|x_{j}-y_{j}\right|}{l}\right),\tag{A.1}
\end{align*}
where $l$ denotes the correlation length. Subsequently, we can define
the KL expansion of $u$ using an infinite set of random variables
$\left\{ \xi_{\boldsymbol{j}}:\boldsymbol{j}\in\mathbb{N}^{n}\right\} $,
each having zero mean and unit variance, as follows. $\forall\boldsymbol{x}\in\Omega$,
\hypertarget{eqA2}{}
\begin{align*}
u\left(\boldsymbol{x}\right)-\bar{u}\left(\boldsymbol{x}\right) & =\sum_{\boldsymbol{j}\in\mathbb{R}^{n}}\gamma_{\boldsymbol{j}}\left(\boldsymbol{x}\right)\xi_{\boldsymbol{j}}\\
 & =\sum_{j_{1}\in\mathbb{R}}\cdots\sum_{j_{n}\in\mathbb{R}}\gamma_{j_{1}\ldots j_{n}}\left(\boldsymbol{x}\right)\xi_{j_{1}\ldots j_{n}}\\
 & =\sum_{j_{1}\in\mathbb{R}}\cdots\sum_{j_{n}\in\mathbb{R}}\prod_{k=1}^{n}g_{j_{k}}\left(x_{k}\right)\xi_{j_{1}\ldots j_{n}}\tag{A.2}
\end{align*}
where $\forall j>0$, if $\zeta_{j}$ solves 
\[
l\zeta_{j}+\tan\left(\frac{\zeta_{j}}{2}\right)=0,\tag{A.3}
\]
and $\zeta_{j+1}>\zeta_{j}>0$, then $\forall x\in\mathbb{R}$, 
\[
g_{j}\left(x\right)=\begin{cases}
2{\displaystyle \sqrt{\frac{l\zeta_{j}}{1+l^{2}\zeta_{j}^{2}}}\frac{\cos\left(\zeta_{j}x\right)}{\sqrt{\zeta_{j}+\sin\left(\zeta_{j}\right)}}} & j\ \mathrm{is\ odd},\\
2{\displaystyle \sqrt{\frac{l\zeta_{j}}{1+l^{2}\zeta_{j}^{2}}}\frac{\sin\left(\zeta_{j}x\right)}{\sqrt{\zeta_{j}-\sin\left(\zeta_{j}\right)}}} & j\ \mathrm{is\ even}.
\end{cases}\tag{A.4}
\]
Therefore, as is required in \S4.1 and \S4.2, a truncated KL expansion
can be easily obtained from the single index form of the expansion
in \hyperlink{eqA2}{Eq. A.2}.

\renewcommand\thefigure{B\arabic{figure}}    
\setcounter{figure}{0}    
\newpage
\section*{Appendix B - Verification of modular deterministic solvers}

\subsection*{B1 - Thermal-Neutronics problem}

Given the stochastic transmittance $H_{T}$ and cross section $\varSigma_{\varphi}$,
we choose analytical functions $T^{*},\varphi^{*}:\forall\boldsymbol{x}\in\Omega,\boldsymbol{\xi}\in\Xi$,
\begin{align*}
T^{*}\left(\boldsymbol{x},\boldsymbol{\xi}\right) & =\left(2+\cos\left(\pi x_{1}\right)\right)^{2}-1,\\
\varphi^{*}\left(\boldsymbol{x},\boldsymbol{\xi}\right) & =3+\cos\left(2\pi x_{2}\right),\tag{B.1}
\end{align*}
which solve the modified system of equations: $\forall\boldsymbol{x}\in\Omega,\boldsymbol{\xi}\in\Xi$,
\begin{align*}
\boldsymbol{\nabla}^{\mathbf{T}}\boldsymbol{\nabla}T\left(\boldsymbol{x},\boldsymbol{\xi}\right)-H_{T}\left(\boldsymbol{x},\boldsymbol{\xi}_{1}\left(\boldsymbol{\xi}\right)\right)T\left(\boldsymbol{x},\boldsymbol{\xi}\right)+\frac{E_{T}\varphi\left(\boldsymbol{x},\boldsymbol{\xi}\right)}{\sqrt{T\left(\boldsymbol{x},\boldsymbol{\xi}\right)+1}}+f_{T}^{*}\left(\boldsymbol{x},\boldsymbol{\xi}_{1}\left(\boldsymbol{\xi}\right)\right) & =0,\\
\boldsymbol{\nabla}^{\mathbf{T}}\left(\sqrt{T\left(\boldsymbol{x},\boldsymbol{\xi}\right)+1}\boldsymbol{\nabla}\varphi\left(\boldsymbol{x},\boldsymbol{\xi}\right)\right)-\frac{\varSigma_{\varphi}\left(\boldsymbol{x},\boldsymbol{\xi}_{2}\left(\boldsymbol{\xi}\right)\right)\varphi\left(\boldsymbol{x},\boldsymbol{\xi}\right)}{\sqrt{T\left(\boldsymbol{x},\boldsymbol{\xi}\right)+1}}+f_{\varphi}^{*}\left(\boldsymbol{x},\boldsymbol{\xi}_{2}\left(\boldsymbol{\xi}\right)\right) & =0.\tag{B.2}
\end{align*}

Here, $\forall\boldsymbol{x}\in\Omega,\boldsymbol{\xi}_{1}\in\Xi_{1},\boldsymbol{\xi}_{2}\in\Xi_{2}$,

\begin{align*}
f_{T}^{*}\left(\boldsymbol{x},\boldsymbol{\xi}_{1}\right) & =2\pi^{2}\left(2\cos\left(\pi x_{1}\right)+\cos\left(2\pi x_{1}\right)\right)+H_{T}\left(\boldsymbol{x},\boldsymbol{\xi}_{1}\right)\left(\left(2+\cos\left(\pi x_{1}\right)\right)^{2}-1\right)\\
 & -E_{T}\frac{3+\cos\left(2\pi x_{2}\right)}{2+\cos\left(\pi x_{1}\right)},\\
f_{\varphi}^{*}\left(\boldsymbol{x},\boldsymbol{\xi}_{2}\right) & =4\pi^{2}\left(2+\cos\left(\pi x_{1}\right)\right)\cos\left(2\pi x_{2}\right)+\varSigma_{\varphi}\left(\boldsymbol{x},\boldsymbol{\xi}_{2}\right)\frac{3+\cos\left(2\pi x_{2}\right)}{2+\cos\left(\pi x_{1}\right)}.\tag{B.3}
\end{align*}

Let $\varepsilon^{m}\left(\boldsymbol{\xi}\right):\forall\boldsymbol{\xi}\in\Xi$,
\[
\varepsilon^{m}\left(\boldsymbol{\xi}\right)=\sqrt{\frac{\int_{\Omega}\left\Vert \left[\begin{array}{c}
T^{*}\left(\boldsymbol{x},\boldsymbol{\xi}\right)-T^{m}\left(\boldsymbol{x},\boldsymbol{\xi}\right)\\
\varphi^{*}\left(\boldsymbol{x},\boldsymbol{\xi}\right)-\varphi^{m}\left(\boldsymbol{x},\boldsymbol{\xi}\right)
\end{array}\right]\right\Vert _{2}^{2}d\boldsymbol{x}}{\int_{\Omega}\left\Vert \left[\begin{array}{c}
T^{*}\left(\boldsymbol{x},\boldsymbol{\xi}\right)\\
\varphi^{*}\left(\boldsymbol{x},\boldsymbol{\xi}\right)
\end{array}\right]\right\Vert _{2}^{2}d\boldsymbol{x}}}\tag{B.4}
\]
denote the mean-square error between the exact and approximate solutions
$T^{m}$ and $\varphi^{m}$ computed using $m\times m$ nodes. Subsequently,
with a convergence tolerance value $\epsilon_{\mathrm{BGS}}=10^{-6}$
and stochastic dimensions to $s_{1}=s_{2}=4$, the sample average
of $\varepsilon^{m}$, denoted here as $\bar{\varepsilon}^{m}$ was
computed using $100$ Monte-Carlo solution samples, for various values
of $m$. \hyperlink{figB1}{Figure B1} illustrates the expected second order rate of decay
in $\bar{\varepsilon}^{m}$.

\begin{figure}
    \hypertarget{figB1}{}
    \centering
    \includegraphics[bb=40bp 220bp 550bp 590bp,scale=0.55]{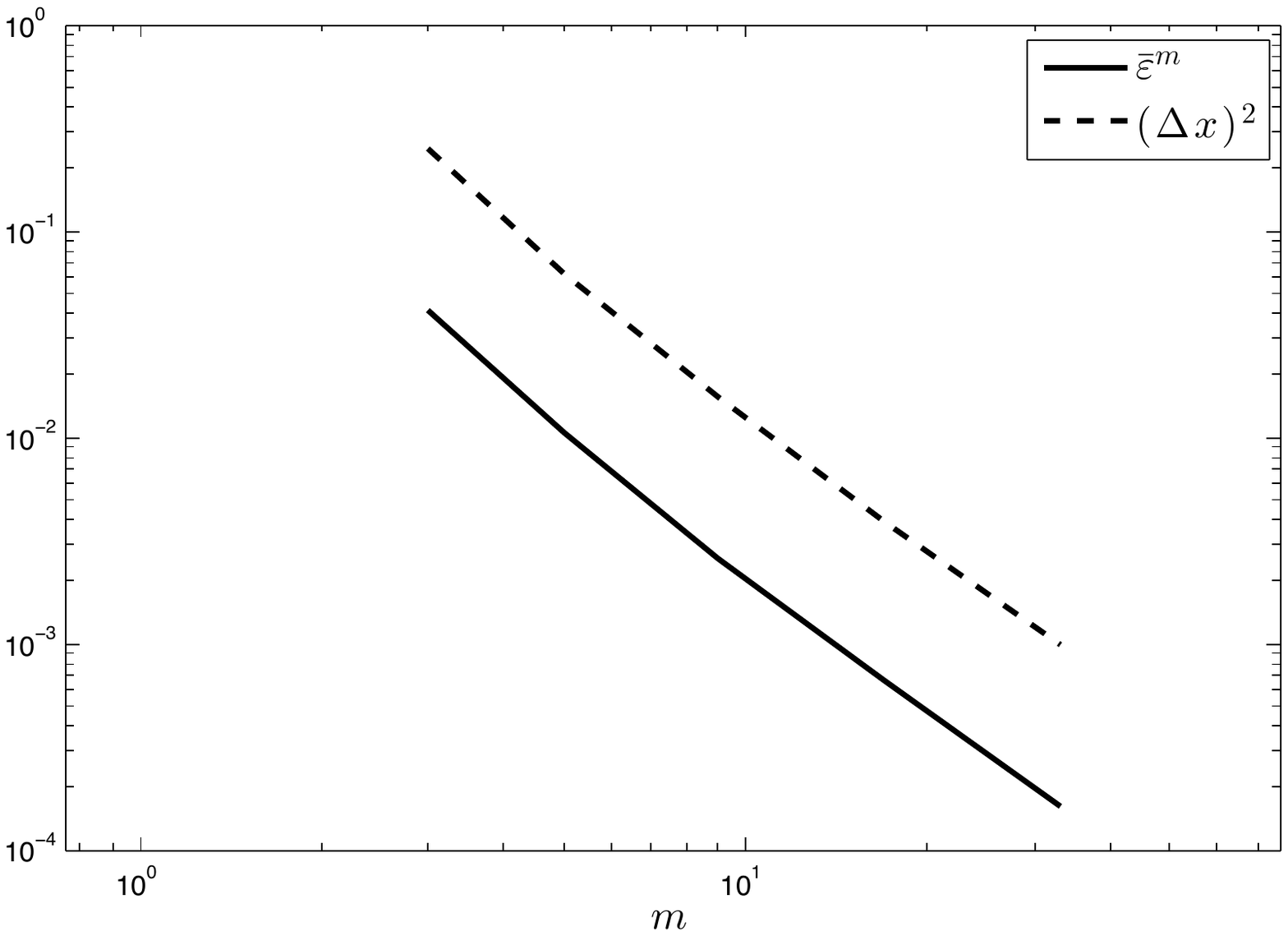} 
    \caption{Average error $\bar{\varepsilon}^{m}$ v.s.
$m$. $\Delta x$ denotes the node spacing.}
\end{figure}

\subsection*{B2 - Boussinesq flow problem}

Following the same procedure in Appendix B1, we choose analytical
functions $u_{1}^{*},u_{2}^{*},p^{*},t^{*}:\forall\boldsymbol{x}\in\Omega,\boldsymbol{\xi}\in\Xi$,
\begin{align*}
u_{1}^{*}\left(\boldsymbol{x},\boldsymbol{\xi}\right) & =-\sin^{2}\left(\pi x_{1}\right)\sin\left(2\pi x_{2}\right),\\
u_{2}^{*}\left(\boldsymbol{x},\boldsymbol{\xi}\right) & =\sin\left(2\pi x_{1}\right)\sin^{2}\left(\pi x_{2}\right),\\
p^{*}\left(\boldsymbol{x},\boldsymbol{\xi}\right) & =\cos\left(\pi x_{1}\right)\cos\left(\pi x_{2}\right),\\
t^{*}\left(\boldsymbol{x},\boldsymbol{\xi}\right) & =\cos\left(\frac{\pi}{2}x_{1}\right)T_{h}\left(x_{2},\boldsymbol{\xi}_{2}\left(\boldsymbol{\xi}\right)\right),\tag{B.5}
\end{align*}
which solve the modified Boussinesq equations: $\forall\boldsymbol{x}\in\Omega,\boldsymbol{\xi}\in\Xi$,
\begin{align*}
\boldsymbol{\nabla}^{\mathbf{T}}\boldsymbol{u}\left(\boldsymbol{x},\boldsymbol{\xi}\right) & =0,\\
\left(\boldsymbol{u}\left(\boldsymbol{x},\boldsymbol{\xi}\right)^{\mathbf{T}}\boldsymbol{\nabla}\right)\boldsymbol{u}\left(\boldsymbol{x},\boldsymbol{\xi}\right)+\boldsymbol{\nabla}p\left(\boldsymbol{x},\boldsymbol{\xi}\right)\\
-\mathrm{Pr}\boldsymbol{\nabla}^{\mathbf{T}}\boldsymbol{\nabla}\boldsymbol{u}\left(\boldsymbol{x},\boldsymbol{\xi}\right)-\mathrm{Pr}\mathrm{Ra}\left(\boldsymbol{x},\boldsymbol{\xi}_{1}\left(\boldsymbol{\xi}\right)\right)T\left(\boldsymbol{x},\boldsymbol{\xi}\right)\boldsymbol{e}_{2}+\boldsymbol{f}_{u}^{*}\left(\boldsymbol{x},\boldsymbol{\xi}\right) & =\boldsymbol{0},\\
\left(\boldsymbol{u}\left(\boldsymbol{x},\boldsymbol{\xi}\right)^{\mathbf{T}}\boldsymbol{\nabla}\right)T\left(\boldsymbol{x},\boldsymbol{\xi}\right)-\boldsymbol{\nabla}^{\mathbf{T}}\boldsymbol{\nabla}T\left(\boldsymbol{x},\boldsymbol{\xi}\right)+f_{T}^{*}\left(\boldsymbol{x},\boldsymbol{\xi}\right) & =0,\tag{B.6}
\end{align*}

Here, $\forall\boldsymbol{x}\in\Omega,\boldsymbol{\xi}\in\Xi$, 
\begin{align*}
\boldsymbol{f}_{u}^{*}\left(\boldsymbol{x},\boldsymbol{\xi}\right) & =2\pi\left[\begin{array}{c}
\sin\left(2\pi x_{1}\right)\sin^{2}\left(\pi x_{1}\right)\sin^{2}\left(2\pi x_{2}\right)\\
\sin\left(2\pi x_{2}\right)\sin^{2}\left(\pi x_{2}\right)\sin^{2}\left(2\pi x_{1}\right)
\end{array}\right]\\
 & -\pi\left[\begin{array}{c}
\sin\left(2\pi x_{1}\right)\sin^{2}\left(\pi x_{1}\right)\sin^{2}\left(2\pi x_{2}\right)-\sin\left(\pi x_{1}\right)\cos\left(\pi x_{2}\right)\\
\sin\left(2\pi x_{2}\right)\sin^{2}\left(\pi x_{2}\right)\sin^{2}\left(2\pi x_{1}\right)-\sin\left(\pi x_{2}\right)\cos\left(\pi x_{1}\right)
\end{array}\right]\\
 & +2\pi^{2}\mathrm{Pr}\left[\begin{array}{c}
\cos\left(2\pi x_{1}\right)\sin\left(2\pi x_{2}\right)-2\sin^{2}\left(\pi x_{1}\right)\sin\left(2\pi x_{2}\right)\\
-\cos\left(2\pi x_{2}\right)\sin\left(2\pi x_{1}\right)+2\sin^{2}\left(\pi x_{2}\right)\sin\left(2\pi x_{1}\right)
\end{array}\right]\\
 & +\Pr\mathrm{Ra}\left(\boldsymbol{x},\boldsymbol{\xi}_{1}\left(\boldsymbol{\xi}\right)\right)\cos\left(\frac{\pi}{2}x_{1}\right)T_{h}\left(x_{2},\boldsymbol{\xi}_{2}\left(\boldsymbol{\xi}\right)\right)\boldsymbol{e}_{2},\\
f_{T}^{*}\left(\boldsymbol{x},\boldsymbol{\xi}\right) & =-\frac{\pi}{2}\left(\sin^{2}\left(\pi x_{1}\right)\sin\left(2\pi x_{2}\right)\sin\left(\frac{\pi}{2}x_{1}\right)+\frac{\pi}{2}\cos\left(\frac{\pi}{2}x_{1}\right)\right)T_{h}\left(x_{2},\boldsymbol{\xi}_{2}\left(\boldsymbol{\xi}\right)\right)\\
 & -\sin\left(2\pi x_{1}\right)\cos\left(\frac{\pi}{2}x_{1}\right)\sin^{2}\left(\pi x_{2}\right)\frac{\partial T_{h}}{\partial x_{2}}\left(x_{2},\boldsymbol{\xi}_{2}\left(\boldsymbol{\xi}\right)\right)\\
 & +\cos\left(\frac{\pi}{2}x_{1}\right)\frac{\partial^{2}T_{h}}{\partial x_{2}^{2}}\left(x_{2},\boldsymbol{\xi}_{2}\left(\boldsymbol{\xi}\right)\right).\tag{B.7}
\end{align*}

Let $\varepsilon^{m}\left(\boldsymbol{\xi}\right):\forall\boldsymbol{\xi}\in\Xi$,
\begin{equation*}
\varepsilon^{m}\left(\boldsymbol{\xi}\right)=\sqrt{\frac{\int_{\Omega}\left\Vert \left[\begin{array}{c}
u_{1}^{*}\left(\boldsymbol{x},\boldsymbol{\xi}\right)-u_{1}^{m}\left(\boldsymbol{x},\boldsymbol{\xi}\right)\\
u_{2}^{*}\left(\boldsymbol{x},\boldsymbol{\xi}\right)-u_{2}^{m}\left(\boldsymbol{x},\boldsymbol{\xi}\right)\\
p^{*}\left(\boldsymbol{x},\boldsymbol{\xi}\right)-p^{m}\left(\boldsymbol{x},\boldsymbol{\xi}\right)\\
t^{*}\left(\boldsymbol{x},\boldsymbol{\xi}\right)-t^{m}\left(\boldsymbol{x},\boldsymbol{\xi}\right)
\end{array}\right]\right\Vert _{2}^{2}d\boldsymbol{x}}{\int_{\Omega}\left\Vert \left[\begin{array}{c}
u_{1}^{*}\left(\boldsymbol{x},\boldsymbol{\xi}\right)\\
u_{2}^{*}\left(\boldsymbol{x},\boldsymbol{\xi}\right)\\
p^{*}\left(\boldsymbol{x},\boldsymbol{\xi}\right)\\
t^{*}\left(\boldsymbol{x},\boldsymbol{\xi}\right)
\end{array}\right]\right\Vert _{2}^{2}d\boldsymbol{x}}}\tag{B.8}
\end{equation*}
denote the mean-square error between the exact and approximate solutions
$u_{1}^{m}$, $u_{2}^{m}$, $p^{m}$, $t^{m}$, computed using $m\times m$ cells.

Subsequently, the average error $\bar{\varepsilon}^{m}$, was computed
using 100 Monte-Carlo samples for various values of $m$, shown in
\hyperlink{figB2}{Figure B2}, keeping the tolerance at $\epsilon_{\mathrm{BGS}}=10^{-6}$
and stochastic dimensions to $s_{1}=s_{2}=4$. As expected, a second
order rate of decay rate is observed in $\bar{\varepsilon}^{m}$.

\begin{figure}
    \hypertarget{figB2}{}
    \centering
    \includegraphics[bb=40bp 220bp 550bp 590bp,scale=0.55]{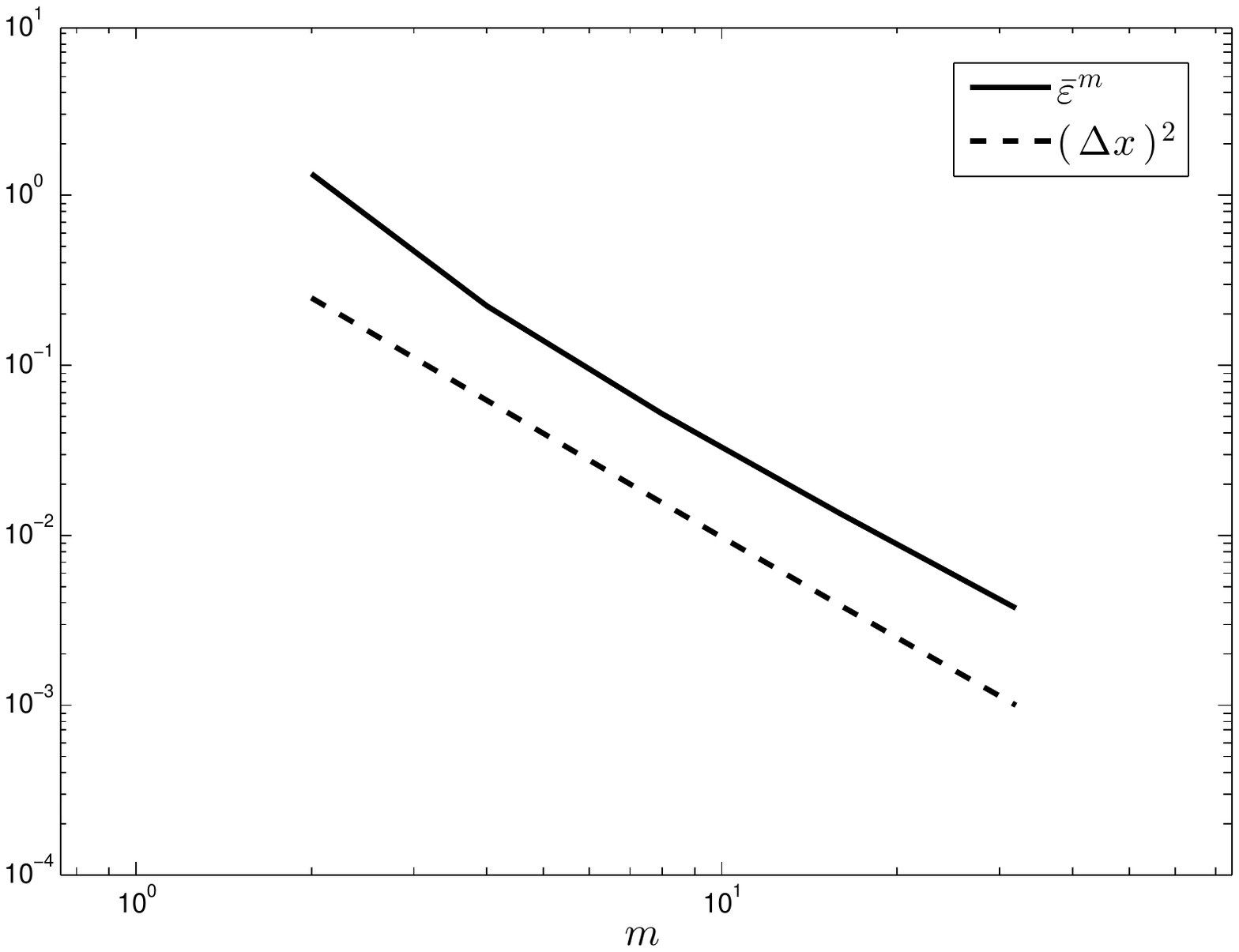} 
    \caption{Average error $\bar{\varepsilon}^{m}$ v.s.
$m$. $\Delta x$ denotes the cell side length.}
\end{figure}


\begin{thebibliography}{10}
\bibitem{key-1}\textsc{G. Casella and C. Robert}, \emph{Monte
Carlo statistical methods}, Springer, New York, 1999.
\hypertarget{ref1}{}
\bibitem{key-2}\textsc{P. Glasserman}, \emph{Quasi-Monte Carlo},
Monte Carlo Methods in Financial Engineering, Springer, New York, 2003,
pp. 281-337.
\hypertarget{ref2}{}
\bibitem{key-3}\textsc{H. Rauhut and R. Ward}, \emph{Sparse
Legendre expansions via $\ell$1-minimization},
J. Approx. Theory, 164.5 (2012), pp. 517-533.
\hypertarget{ref3}{}
\bibitem{key-4}\textsc{J. Peng, J. Hampton, and A. Doostan,} \emph{A
weighted $\ell$1-minimization approach for sparse polynomial
chaos expansions,} J. Comput. Phys.,
267 (2014), pp. 92-111.
\hypertarget{ref4}{}
\bibitem{key-5}\textsc{D. Xiu and G. Karniadakis,} \emph{Modeling
uncertainty in flow simulations via generalized polynomial chaos,}
J. Comput. Phys., 187.1 (2003), pp. 137-167.
\hypertarget{ref5}{}
\bibitem{key-6}\textsc{H. Matthies, R. Niekamp and J. Steindorf,}
 \emph{Algorithms for strong coupling procedures.}
Computer Methods Appl. Mech. Engrg., 195 (2006), pp. 2028-2049.
\hypertarget{ref6}{}
\bibitem{key-7}\textsc{C. Felippa, K. Park and C. Farhat,}
 \emph{Partitioned analysis of coupled mechanical systems.}
Computer Methods Appl. Mech. Engrg., 190 (2006), pp. 3247-3270.
\hypertarget{ref7}{}
\bibitem{key-8}\textsc{P. Constantine, A. Doostan and G.
Iaccarino,} \emph{A hybrid collocation/Galerkin scheme for
convective heat transfer problems with stochastic boundary conditions,}
International Journal for Numerical Methods in Engineering, 80.7
(2009), pp. 868-880.
\hypertarget{ref8}{}
\bibitem{key-9}\textsc{X. Chen, B. Ng, Y. Sun and C. Tong,}
 \emph{A flexible uncertainty quantification method for linearly
coupled multi-physics systems,} J. Comput. Phys.,
Physics 248 (2013), pp. 383-401.\hypertarget{ref9}{}
\bibitem{key-10}\textsc{X. Chen, B. Ng, Y. Sun and C. Tong,} \emph{A
computational method for simulating subsurface flow and reactive transport
in heterogeneous porous media embedded with flexible uncertainty quantification,}
Water Resources Research, 49 (2013), pp. 5740-5755.
\hypertarget{ref10}{}
\bibitem{key-11}\textsc{M. Hadigol, A. Doostan, H. Matthies 
and R. Niekamp,} \emph{Partitioned treatment of uncertainty
in coupled domain problems: A separated representation approach,}
Computer Methods Appl. Mech. Engrg., 274 (2014), pp.
103-124.
\hypertarget{ref11}{}
\bibitem{key-12}\textsc{P. Constantine, E. Phipps and T. Wildey}, \emph{Efficient uncertainty propagation for network
multi-physics systems,} International Journal for Numerical
Methods in Engineering, (2014).\hypertarget{ref12}{}
\bibitem{key-13}\textsc{M. Arnst, R. Ghanem, E. Phipps and J.
Red-Horse,} \emph{Dimension reduction in stochastic modeling of coupled problems,}
International Journal for Numerical Methods in Engineering, 92 (2012), pp. 940-968.
\hypertarget{ref13}{}
\bibitem{key-14}\textsc{M. Arnst, R. Ghanem, E. Phipps and J.
Red-Horse,} \emph{Measure transformation and efficient quadrature
in reduced-dimensional stochastic modeling of coupled problems,}
International Journal for Numerical Methods in Engineering, 92 (2012), pp. 1044-1080.
\hypertarget{ref14}{}
\bibitem{key-15}\textsc{M. Arnst, R. Ghanem, E. Phipps and J. Red-Horse,} 
 \emph{Reduced chaos expansions with random coefficientsin
reduced-dimensional stochastic modeling of coupled problems,}
International Journal for Numerical Methods in Engineering, 97 (2014), pp. 352-376.
\hypertarget{ref15}{}
\bibitem{key-16}\textsc{C. Soize  and R. Ghanem,} \emph{Reduced
chaos decomposition with random coefficients of vector-valued random
variables and random fields,} Computer Methods Appl. Mech. Engrg., 198  (2009), pp. 1926-1934.\hypertarget{ref16}{}
\bibitem{key-17}\textsc{P. Davis,} \emph{Interpolation and approximation,} 
Dover Publications, 1975.
\hypertarget{ref17}{}
\bibitem{key-18}\textsc{W. Gautschi}, \emph{Orthogonal polynomials:
applications and computation,} Acta Numer., 5 (1996), pp.
45-119.\hypertarget{ref18}{}
\bibitem{key-19}\textsc{W. Gautschi,} \emph{On generating orthogonal
polynomials,} SIAM J. Sci. Comp., 3 (1982), pp. 289-317.\hypertarget{ref19}{}
\bibitem{key-20}\textsc{R. Cameron and W. Martin,} \emph{The
orthogonal development of non-linear functionals in series of Fourier-Hermite
functionals,} Ann. of Math., (1947), pp. 385-392.
45-119.\hypertarget{ref20}{}
\bibitem{key-21}\textsc{B. Silverman,} \emph{Density estimation for statistics
and data analysis,}, CRC press, 1986.
\hypertarget{ref21}{}
\bibitem{key-22}\textsc{I. Babsuka, R. Tempone and G. Zouraris,} \emph{Galerkin finite element approximations of stochastic elliptic partial differential equations}, SIAM J. Numer. Anal., 42 (2004), pp. 800-825.
\hypertarget{ref22}{}
\bibitem{key-23}\textsc{O. Ernst and U. Ellmann,} \emph{Stochastic galerkin matrices.}, SIAM J. Matrix Anal. App., 31 (2010), pp. 1848-1872.
\hypertarget{ref23}{}
\bibitem{key-24}\textsc{K. Park, C. Carlos, A. Felippa and R. Ohayon,}   
 \emph{Partitioned formulation of internal fluid\textendash structure
interaction problems by localized Lagrange multipliers,}
Computer Methods Appl. Mech. Engrg., 190.24 (2001), pp.
2989-3007.
\hypertarget{ref24}{}
\bibitem{key-25}\textsc{H. Hartley,} \emph{The modified Gauss-Newton
method for the fitting of non-linear regression functions by least
squares,} Technometrics, 3 (1961), pp. 269-280.
\hypertarget{ref25}{}
\bibitem{key-26}\textsc{J. More}, \emph{The Levenberg-Marquardt
algorithm: implementation and theory,} in Numerical
analysis, Springer, Heidelberg, 1978, pp. 105-116.
\hypertarget{ref26}{}
\bibitem{key-27}\textsc{M. Joosten, W. Dettmer and D. Peric,} \emph{Analysis of the block Gauss-Seidel solution procedure for a strongly coupled model problem with reference to fluid-structure interaction,} International Journal for Numerical Methods in Engineering, 78.7
(2009), pp. 757-778.
\hypertarget{ref27}{}
\bibitem{key-28}\textsc{F. Nobile, Fabio, R. Tempone and C. Webster,}
 \emph{A sparse grid stochastic collocation method for partial
differential equations with random input data,} SIAM
J. Numer. Anal., 46 (2008), pp. 2309-2345.
\hypertarget{ref28}{}
\bibitem{key-29}\textsc{M. Loeve,} \emph{Probability Theory. Foundations. Random
Sequences,}, D. Van Nostrand Company, New York, 1955.
\hypertarget{ref29}{}
\bibitem{key-30}\textsc{K. Petersen and M. Pedersen,}
 \emph{The matrix cookbook,} Technical University
of Denmark, (2008), pp. 7-15.
\hypertarget{ref30}{}
\bibitem{key-31}\textsc{J. Nelder and R. Mead,} \emph{A
simplex method for function minimization,} The Computer
Journal,  7 (1965), pp. 308-313.
\hypertarget{ref31}{}
\bibitem{key-32}\textsc{V. Tchakaloff,} \emph{Formules de cubatures mécaniquesa coefficients non négatifs,} Bull.
Sci. Math., 81 (1957), pp. 123-134.
\hypertarget{ref32}{}
\bibitem{key-33}\textsc{J. Lamarsh,} \emph{Introduction to nuclear reactor theory,} Addison-Wesley, Reading, 1966.\hypertarget{ref33}{}
\bibitem{key-34}\textsc{T. Hughes,} \emph{The finite element method:
linear static and dynamic finite element analysis,}, Dover
Publications, 2012.
\hypertarget{ref34}{}
\bibitem{key-35}\textsc{P. Roache,} \emph{Code verification
by the method of manufactured solutions,} Journal of
Fluids Engineering, 124 (2002), pp. 4-10.
\hypertarget{ref35}{}
\bibitem{key-36}\textsc{C. Paige and M. Saunders,} \emph{LSQR: An algorithm for sparse linear equations and sparse least squares,} ACM Transactions on Mathematical Software, 8.1 (1982), pp. 43-71.
\hypertarget{ref36}{}
\bibitem{key-37}\textsc{P. Constantine, D. Gleich, and G. Iaccarino,} \emph{A Factorization of the Spectral Galerkin System for Parameterized Matrix Equations: Derivation and Applications,} SIAM J. Matrix Anal. App., 33.5 (2011), pp. 2995-3009.
\hypertarget{ref37}{}
\bibitem{key-38}\textsc{A. Saltelli, K. Chan and E. Scott,} \emph{Sensitivity analysis,} Wiley, New York, 2000.
\hypertarget{ref38}{}
\bibitem{key-39}\textsc{L. Kantha and C. Clayson,} \emph{Numerical
models of oceans and oceanic processes,} Academic press, 2000.\hypertarget{ref39}{}
\bibitem{key-40}\textsc{J. Massaguer and J. Zahn,} \emph{Cellular
convection in a stratified atmosphere,} Astronom. and
Astrophys. Lib., 87 (1980), pp. 315-327.
\hypertarget{ref40}{}
\bibitem{key-41}\textsc{D. de Vahl,} \emph{Laminar natural
convection in an enclosed rectangular cavity,} International
Journal of Heat and Mass Transfer, 11 (1968), pp. 1675-1693.
\hypertarget{ref41}{}
\bibitem{key-42}\textsc{R. LeVeque}, \emph{Finite volume methods for
hyperbolic problems,} Cambridge University Press, 2002.
\hypertarget{ref42}{}
\bibitem{key-43}\textsc{M. Arnst, C. Soize and R. Ghanem,}
 \emph{A Hybrid Sampling/Spectral Method for Solving Stochastic
Coupled Problems,} SIAM Journal on Uncertainty Quantification, (2013), pp. 218-243.
\hypertarget{ref43}{}
\bibitem{key-44}\textsc{X. Wan and G. Karniadakis,} \emph{Multi-element
generalized polynomial chaos for arbitrary probability measures,}
SIAM J. Sci. Comp., 28 (2006), pp. 901-928.
\hypertarget{ref44}{}
\bibitem{key-44}\textsc{O Le Maitre, O. Knio, H. Najm, and R. Ghanem,} \emph{Uncertainty propagation using Wiener-Haar
expansions,} J. Comput. Phys., 197 (2004), pp. 28-57.
\hypertarget{ref44}{}
\bibitem{key-46}\textsc{T. Lukaczyk, F. Palacios, J. Alonso and P. Constantine,} \emph{Active Subspaces for Shape
Optimization,} 10th AIAA Multidisciplinary Design Optimization Conference, Maryland, 2014.
\hypertarget{ref45}{}

\end{thebibliography}
\end{document}